\documentclass[10pt]{amsart}
\usepackage{amssymb,xypic}
\newtheorem{theorem}{Theorem}
\newtheorem{proposition}{Proposition}
\newtheorem{lemma}[proposition]{Lemma}
\newtheorem{cor}[proposition]{Corollary}

\theoremstyle{definition}
\newtheorem*{remark}{Remark}
\newcommand{\Q}{{\mathbb Q}}
\newcommand{\Z}{{\mathbb Z}}
\newcommand{\N}{{\mathbb N}}
\newcommand{\R}{\mathbb{R}}

\newcommand{\tensor}{\otimes}
\newcommand{\coker}{\mathrm{coker}}
\newcommand{\sD}{\mathcal{D}}
\newcommand{\PP}{\mathbb{P}}

\newcommand{\sC}{\mathcal{C}}
\newcommand{\elem}{\mathrm{Elem}}
\newcommand{\sS}{\mathcal{S}}

\newcommand{\sF}{\mathcal{F}}
\newcommand{\sV}{{\mathcal V}}
\newcommand{\sG}{\mathcal{G}}
\newcommand{\sE}{{\mathcal E}}

\newcommand{\sP}{\mathcal{P}}
\newcommand{\sL}{{\mathcal L}}
\newcommand{\oC}{\overline{C}}
\newcommand{\oL}{\overline{L}}
\newcommand{\sH}{{\mathcal H}}
\newcommand{\sAut}{{\mathcal{A}ut}}
\newcommand{\semidirect}{\ltimes}
\newcommand{\Maps}{\mathrm{Maps}}
\newcommand{\gr}{{\rm gr}\,}
\newcommand{\FL}{{\mathbb{FL}}}

\newcommand{\SPL}{{\mathbb{SPL}}}
\newcommand{\sSPL}{{\mathcal{SPL}}}
\newcommand{\ET}{{\mathbb{ET}}}
\newcommand{\GL}{{\rm GL}\,}
\newcommand{\sta}{{\rm star}\,}
\newcommand{\rank}{{\rm rank}\,}

\newcommand{\limdir}[1]{{\displaystyle{\varinjlim_{#1}}}}
\renewcommand{\tilde}{\widetilde}
\def\YEAR{\year}\newcount\VOL\VOL=\YEAR\advance\VOL by-1995
\def\firstpage{1}\def\lastpage{1000}
\def\received{}\def\revised{}
\def\communicated{}

\makeatletter
\def\magnification{\afterassignment\m@g\count@}
\def\m@g{\mag=\count@\hsize6.5truein\vsize8.9truein\dimen\footins8truein}
\makeatother

\font\eightrm=cmr8
\font\caps=cmcsc10                    
\font\Caps=cmcsc10 scaled \magstep1   

\makeatletter
\@specialpagefalse\headheight=8.5pt
\def\DocMath{}
\renewcommand{\@evenhead}{%
    \ifnum\thepage>\lastpage\rlap{\thepage}\hfill%
    \else\rlap{\thepage}\slshape\leftmark\hfill{\caps\SAuthor}\hfill\fi}%
\renewcommand{\@oddhead}{%
    \ifnum\thepage=\firstpage{\DocMath\hfill\llap{\thepage}}%
    \else{\slshape\rightmark}\hfill{\caps\STitle}\hfill\llap{\thepage}\fi}%
\makeatother

\def\TSkip{\bigskip}
\newbox\TheTitle{\obeylines\gdef\GetTitle #1
\ShortTitle  #2
\SubTitle    #3
\Author      #4
\ShortAuthor #5
\EndTitle
{\setbox\TheTitle=\vbox{\baselineskip=20pt\let\par=\cr\obeylines%
\halign{\centerline{\Caps##}\cr\noalign{\medskip}\cr#1\cr}}%
	\copy\TheTitle\TSkip\TSkip%
\def\next{#2}\ifx\next\empty\gdef\STitle{#1}\else\gdef\STitle{#2}\fi%
\def\next{#3}\ifx\next\empty%
    \else\setbox\TheTitle=\vbox{\baselineskip=20pt\let\par=\cr\obeylines%
    \halign{\centerline{\caps##} #3\cr}}\copy\TheTitle\TSkip\TSkip\fi%
\centerline{\caps #4}\TSkip\TSkip%
\def\next{#5}\ifx\next\empty\gdef\SAuthor{#4}\else\gdef\SAuthor{#5}\fi%
\ifx\received\empty\relax
    \else\centerline{\eightrm Received: \received}\fi%
\ifx\revised\empty\TSkip%
    \else\centerline{\eightrm Revised: \revised}\TSkip\fi%
\ifx\communicated\empty\relax
    \else\centerline{\eightrm Communicated by \communicated}\fi\TSkip\TSkip%
\catcode'015=5}}\def\Title{\obeylines\GetTitle}
\def\Abstract{\begingroup\narrower
    \parskip=\medskipamount\parindent=0pt{\caps Abstract. }}
\def\EndAbstract{\par\endgroup\TSkip}

\long\def\MSC#1\EndMSC{\def\arg{#1}\ifx\arg\empty\relax\else
     {\par\narrower\noindent%
     2000 Mathematics Subject Classification: #1\par}\fi}

\long\def\KEY#1\EndKEY{\def\arg{#1}\ifx\arg\empty\relax\else
	{\par\narrower\noindent Keywords and Phrases: #1\par}\fi\TSkip}

\newbox\TheAdd\def\Addresses{\vfill\copy\TheAdd\vfill
    \ifodd\number\lastpage\vfill\eject\phantom{.}\vfill\eject\fi}
{\obeylines\gdef\GetAddress #1
\Address #2 
\Address #3
\Address #4
\EndAddress
{\def\xs{4.3truecm}\parindent=0pt
\setbox0=\vtop{{\obeylines\hsize=\xs#1\par}}\def\next{#2}
\ifx\next\empty 
     \setbox\TheAdd=\hbox to\hsize{\hfill\copy0\hfill}
\else\setbox1=\vtop{{\obeylines\hsize=\xs#2\par}}\def\next{#3}
\ifx\next\empty 
     \setbox\TheAdd=\hbox to\hsize{\hfill\copy0\hfill\copy1\hfill}
\else\setbox2=\vtop{{\obeylines\hsize=\xs#3\par}}\def\next{#4}
\ifx\next\empty\ 
     \setbox\TheAdd=\vtop{\hbox to\hsize{\hfill\copy0\hfill\copy1\hfill}
                \vskip20pt\hbox to\hsize{\hfill\copy2\hfill}}
\else\setbox3=\vtop{{\obeylines\hsize=\xs#4\par}}
     \setbox\TheAdd=\vtop{\hbox to\hsize{\hfill\copy0\hfill\copy1\hfill}
	        \vskip20pt\hbox to\hsize{\hfill\copy2\hfill\copy3\hfill}}
\fi\fi\fi\catcode'015=5}}\gdef\Address{\obeylines\GetAddress}

\begin{document}
\Title{K-theory and the enriched Tits building}
\ShortTitle
\SubTitle   
\Author 
Madhav V. Nori\\
Vasudevan Srinivas 
\ShortAuthor 
Nori and Srinivas 
\EndTitle
{\it To A. A. Suslin with admiration, on his sixtieth birthday.}\\[5mm]
\Abstract 
Motivated by the splitting principle, we define certain simplicial complexes associated to an associative ring $A$, 
which have an action of the general linear group $GL(A)$. This leads to an exact sequence, involving Quillen's algebraic 
K-groups of $A$ and the symbol map. Computations in low degrees lead to another view on Suslin's theorem on the Bloch group,
 and perhaps show a way towards possible generalizations.
\EndAbstract
\MSC 
\EndMSC
\KEY 
\EndKEY
\Address 
Dept. of Mathematics
The University of Chicago
5734 University Ave.
Chicago, IL-60637
USA
\Address
School of Mathematics
Tata Institute of Fundamental Research, 
Homi Bhabha Road
Colaba
Mumbai-400005
India
\Address 
\email{nori@math.uchicago.edu}
\email{srinivas@math.tifr.res.in}
\Address
\EndAddress
The homology of $GL_n(A)$ has been studied in great depth by A.A. Suslin.  
In some of his works (\cite{SuslinMilnor} and \cite{SuslinBloch} for
 example), the action of $GL_n(A)$ on certain simplicial complexes
 facilitated his homology computations.

 We introduce  three simplicial complexes in this paper. They are motivated by the splitting principle. The description of these spaces is given below.
This is followed by the little information we possess on their homology.
After that comes the connection with K-theory.

These objects are defined 
quickly in the context of affine algebraic groups as follows. Let
$G$ be a connected algebraic group\footnotemark \footnotetext
{Gopal Prasad informed us
that we should take $G$ reductive or $k$ perfect.}
 defined over a field $k$. The collection
of minimal parabolic subgroups $P\subset G$ is denoted by $FL(G)$ and
the collection of maximal $k$-split tori $T\subset G$ is denoted by
$SPL(G)$. The simplicial complex $\FL(G)$ has $FL(G)$ as its set of vertices.
Minimal parabolics $P_0,P_1,...,P_r$ of $G$ form an $r$-simplex if 
their intersection
contains a maximal $k$-split torus\footnotemark\footnotetext{John Rognes has an analogous construction with maximal parabolics replacing minimal parabolics. His spaces, different homotopy
 types from ours,  are connected with K-theory as well, see \cite{Rognes} . }. The dimension of $\FL(G)$ is one less
than the order of the Weyl group of any $T\in SPL(G)$. 
Dually, we define $\SPL(G)$ as
the simplicial complex with $SPL(G)$ as its set of vertices, and
$T_0,T_1,...,T_r$ forming an $r$-simplex if they are all contained in
 a minimal parabolic. In general, $\SPL(G)$ is infinite dimensional. 
 
That both  $\SPL(G)$ and $\FL(G)$ have the same 
homotopy type can be deduced from corollary~\ref{comparison},
which is a general principle.  A third simplicial complex, denoted by 
$\ET(G)$, which we refer to as the {\em enriched Tits building},  is  
better suited for homology computations. This is the simplicial complex 
whose simplices are (nonempty) chains of the partially ordered set 
$\sE(G)$  whose definition follows.
For a parabolic subgroup $P\subset G$, we denote by $U(P)$ its
unipotent radical and by $j(P):P\to P/U(P)$ the given morphism. Then
$\sE(G)$ is the set of pairs $(P,T)$ where $P\subset G$ is a parabolic subgroup
and $T\subset P/U(P)$ is a maximal
$k$-split torus.   
We say $(P',T')\leq (P,T)$ in $\sE(G)$ if $P'\subset P$
and $j(P)^{-1}(T)\subset j(P')^{-1}T'$. Note that $\dim\sE(G)$ is
the split rank of the quotient of $G/U(G)$ by its center. 
Assume for the moment that this quotient is a simple algebraic group. 
Then $(P,T)\mapsto P$
gives a map to the cone of the Tits building. The topology of $\ET(G)$ 
is  more complex than the topology of the Tits building, which is well 
known to be a bouquet of spheres.

When $G=GL(V)$, we denote the above three simplicial complexes by
$\FL(V)$, $\SPL(V)$ and $\ET(V)$. These constructions have 
simple analogues even when one is working over an arbitrary associative 
ring $A$. Their precise definition is given with some motivation in 
section 2. Some  basic properties of these spaces are also established 
in section 2. Amongst them is Proposition ~\ref{mainflags} which shows 
that $\ET(V)$ and $\FL(V)$ have the same homotopy type. 

$\ET(A^n)$ has a polyhedral decomposition (see lemma~\ref{poly}). This 
produces a spectral sequence (see Theorem~\ref{ss}) that computes its
homology. The $E^1_{p,q}$ terms and the differentials $d^1_{p,q}$
are recognisable since they involve only the homology groups of 
$\ET(A^a)$ for $a<n$. The differentials $d^r_{p,q}$ for $r>1$ are not 
understood well enough, however.

There are natural inclusions $\ET(A^n)\hookrightarrow
 \ET(A^{d})$ for $d>n$, and the induced
map on homology factors through
\begin{center} $H_m(\ET(A^n))\to H_0(E_n(A),H_m(\ET(A^n)))
\to H_m(\ET(A^d))$.
\end{center}
where $E_n(A)$ is the group of elementary matrices
(see Corollary~\ref{elementary}). 

For the remaining statements on the homology of $\ET(A^n)$, we assume 
that $A$ is a commutative ring with many units, in the sense of Van der 
Kallen. See \cite{Mirzaii} for a nice exposition of the definition and 
its consequences. Commutative local rings $A$ with infinite residue 
fields are examples of such rings. Under this assumption, $E_n(A)$ can 
be replaced by $GL_n(A)$ in the above statement. 

We have observed that $\ET(A^{m+1})$ has dimension $m$. Thus it is  
natural question to ask whether
\begin{center} $H_0(GL_{m+1}(A),H_m(\ET(A^{m+1})))\otimes\Q
\to H_m(\ET(A^d))\otimes\Q$
\end{center}
is an isomorphism when $d>m+1$. Theorem~\ref{stab} asserts that this is true
for $m=1,2,3$. The statement is true in general (see 
Proposition~\ref{generalstab}) if a certain Compatible Homotopy Question 
has an affirmative answer. The higher differentials of the spectral 
sequence can be dealt with if this is true.  Proposition~\ref{three} 
shows that this holds in some limited situations.

The computation of $H_0(GL_3(A),H_2(\ET(A^3)))\otimes\Q$ is carried out 
at in the last lemma of the paper. This is intimately connected with 
Suslin's result (see \cite{SuslinBloch}) 
connecting $K_3$ and the Bloch group. A closed form for 
$H_0(GL_4(A),H_3(\ET(A^4))$ is awaited. This should impact on the study
of $K_4(A)$. 

We now come to the connection with the Quillen K-groups $K_i(A)$ as 
obtained by his plus construction.

\begin{sloppypar}
$GL(A)$ acts on the geometric realisation  $|\SPL(A^{\infty})|$ and thus we have the Borel construction,
namely the quotient of $|\SPL(A^{\infty})|\times EGL(A)$ by $GL(A)$, 
 a familiar object in the study of equivariant 
homotopy. We denote this space by  $\SPL(A^{\infty})//GL(A)$.
We apply Quillen's plus construction to $\SPL(A^{\infty})//GL(A)$ 
and a suitable perfect subgroup of its fundamental group to obtain a 
space $Y(A)$. Proposition~\ref{Hspace} shows  that $Y(A)$ is an H-space 
and that the natural map $Y(A)\to BGL(A)^+$ is an H-map. Its homotopy 
fiber, denoted by 
$\SPL(A^{\infty})^+$, is thus also a H-space. The $n$-th homotopy group of
$\SPL(A^{\infty})^+$ at its canonical base point is denoted by $L_n(A)$.
There is of course a natural map $\SPL(A^{\infty})\to 
\SPL(A^{\infty})^+$. That this map is a homology isomorphism is shown  
in lemma~\ref{asser(ii)}. This assertion is easy, but not tautological: 
it relies once again on the triviality of the action of $E(A)$ on 
$H_*(\SPL(A^{\infty}))$. As a consequence of this lemma, 
$L_n(A)\otimes\Q$ is identified with the primitive rational 
homology of $\SPL(A^{\infty})$, or equivalently, that of 
$\ET(A^{\infty})$.
\end{sloppypar}

We have the inclusion $N_n(A)\hookrightarrow GL_n(A)$, where 
$N_n(A)$ is the semidirect product of the permutation group
 $S_n$ with $(A^{\times})^n$. Taking direct limits over $n\in\N$,
 we obtain $N(A)\subset GL(A)$. Let $H'\subset N(A)$ be the infinite
alternating group and let $H$ be the normal subgroup generated by $H'$.
Applying Quillen's plus construction
to the space $BN(A)$ with respect to $H$, we obtain $BN(A)^+$. Its
$n$-th  homotopy group is defined to be $\sH_n(A^{\times})$. From the
 Dold-Thom theorem, it is easy to see that
$\mathcal{H}_n(A^{\times})\tensor\Q$ is isomorphic to
the group homology $H_n(A^{\times})\tensor\Q$. When $A$ is commutative,
this is simply $\wedge^n_{\Q}(A^{\times}\tensor\Q)$. 
Proposition~\ref{stable} identifies
the groups $\sH_n(A^{\times})$ with certain stable homotopy groups.
Its proof was shown to us by J. Peter May.
 It is sketched in the text of the paper after the proof of the Theorem below.

\begin{theorem}\label{theorem1}
Let $A$ be a Nesterenko-Suslin ring. 
Then there is a long exact sequence, functorial in $A$:
$$\displaylines{
\cdots \to L_2(A)\to \mathcal{H}_2(A^{\times})\to K_2(A)\to L_1(A)\to \mathcal{H}_1(A^{\times}) \hfill\cr\hfill
\to K_1(A)\to 
L_0(A)\to 0\cr}$$
\end{theorem}

We call a ring $A$ Nesterenko-Suslin if it satisfies the hypothesis
of Remark 1.13 of their paper \cite{NS}.  The  precise requirement is
that for every finite set $F$, there is a function $f_F:F\to  $ the 
\emph{center} of $A$
so that the sum $\Sigma\{f_F(s):s\in S\}$ is a 
\emph{unit} of $A$ for every nonempty $S\subset F$. If $k$ is an infinite field, every 
associative $k$-algebra
 is Nesterenko-Suslin, and so is every commutative ring with many units in
the sense of Van der Kallen.  Remark 1.13 of  Nesterenko-Suslin
\cite{NS} permits us
  to ignore unipotent radicals. This is used crucially in the proof
of Theorem 1 (see also Proposition~\ref{suslin-lemma2}). 

In the first draft of the paper, we conjectured that this theorem is true
without any hypothesis on $A$. Sasha Beilinson then brought to our 
attention Suslin's paper \cite{SuslinVolodin} on the equivalence of 
Volodin's K-groups and Quillen's. From Suslin's description of 
Volodin's spaces, it is possible to show that these spaces are homotopy 
equivalent to the total space of the $N_n(A)$-torsor on $\FL(A^n)$ given 
in section 2  of this paper. This requires proposition~\ref{covering} 
and a little organisation. Once this is done, Corollary~\ref{elementary} 
can also be obtained from Suslin's set-up. The statement  ``$X(R)$ is 
acyclic'' stated and proved by Suslin in \cite{SuslinVolodin} now 
validates Proposition~\ref{suslin-lemma2} at the infinite level, thus 
showing that Theorem~\ref{theorem1} is true without any hypothesis on 
$A$. The details have not been included here. 

R. Kottwitz informed us that the maximal simplices of $\FL(V)$ are referred to
as ``regular stars'' in the work of Langlands( see \cite{Langlands}).

We hope that this paper will eventually connect with mixed Tate motives 
(see\cite{BK},\cite{BGSV}). 

The arrangement of the paper is a follows. 
Section 1 has some topological preliminaries
 used through most of the paper. The proofs of Corollary~\ref{comparison}
and Proposition~\ref{mainflags} rely on Quillen's Theorem A. Alternatively,
they can both be proved directly by repeated applications of
 Proposition~\ref{covering}.
The definitions of $\SPL(A^n),\FL(A^n),\ET(A^n)$ and first properties
are given in section two. The next four sections are devoted to the 
proof of Theorem 1. The last four sections are concerned with the 
homology of $\ET(A^n)$.

The lemmas, corollaries and propositions 
 are labelled   sequentially. For instance, corollary 9 is followed by lemma 10
 and later by proposition 11; there is no proposition 10 or corollary 10. 
The other numbered statements are the three theorems. Theorems 2 and 3 are 
stated and proved in sections 7 and 9 respectively.  Section 0 records some assumptions
and notation, some perhaps non-standard, that are used in the paper. The reader might
find it helpful to glance at this section for notation  regarding  elementary matrices
 and the Borel construction  and the use of 
 ``simplicial complexes''.

\setcounter{section}{-1}
\section{Assumptions and Notation}

\begin{center} \emph{Rings, Elementary matrices, $\elem(W\hookrightarrow V),\elem(V,q)$, $\sL(V)$ and $\sL_p(V)$}
\end{center}
We are concerned with the Quillen K-groups of a ring $A$. 

\emph{We assume that $A$ has the
following property: if $A^m\cong A^n$ as left $A$-modules, then $m=n$.} The phrase ``$A$-module" 
always means left $A$-module.

For a finitely generated free $A$-module $V$, the collection of $A$-submodules $L\subset V$
so that (i) $V/L$ is free and (ii) $L$ is free of rank one , is denoted by $\sL(V)$.

$\sL_p(V)$ is the collection of subsets $q$ of cardinality $(p+1)$ of $\sL(V)$ so that 
\\ $\oplus\{L:L\in q \}\to V$ is a monomorphism whose cokernel is free.

Given an $A$-submodule $W$ of a $A$-module $V$ so that the
short exact sequence 
\begin{center}$0\to W\to V\to V/W\to 0$
\end{center}
is split, we have the subgroup $\elem(W\hookrightarrow V)\subset Aut_A(V)$, defined 
as follows. Let $H(W)$ be the group of automorphisms $h$  of $V$ so that 
\\ $(id_V-h)V\subset W\subset ker(id_V-h)$.  Let $W'\subset V$  be a submodule
that is complementary  to $W$. Define $H(W')$ in the same manner.  The subgroup
 of $Aut_A(V)$ generated by $H(W)$ and $H(W')$ is   $\elem(W\hookrightarrow V)$.
 It
 does not depend on the choice 
 of $W'$ because $H(W)$ acts transitively on the collection of such $W'$. 
 
 For example, if $V=A^n$, and $W$ is the $A$-submodule generated by any $r$ members
   of
 the given basis of $A^n$, then $\elem(W\hookrightarrow A^n)$ equals $E_n(A)$ , the
 subgroup of elementary matrices in $GL_n(A)$, provided of course that
 $0<r<n$. 
 
 If $V$ is finitely generated free and if  $q\in\sL_p(V)$, the above statement 
 implies that $\elem(L\hookrightarrow V)$ does not depend on the choice of $L\in q$. Thus
  we denote this subgroup by $\elem(V,q)\subset GL(V)$.

\begin{center} \emph{The Borel Construction}
\end{center}
Let $X$ be a topological space equipped with the action of group $G$.  Let $EG$
 be the principal  $G$-bundle on $BG$ (as in \cite{Quillen}). The Borel construction,
  namely the quotient of $X\times EG$ by the $G$-action, is denoted by $X//G$ 
  throughout the paper.

\begin{center} \emph{Categories, Geometric realisations, Posets}
\end{center}
Every category $\sC$ gives rise to a simplicial set, namely its nerve (see \cite{Quillen}). Its 
geometric realisation is denoted by $B\sC$. 

 A poset (partially ordered set)   $P$ gives rise to a 
category. The $B$-construction of this category, by abuse of notation, is denoted by $BP$.
 Associated to
$P$ is the simplicial complex with $P$ as its set of vertices; the simplices are finite non-empty 
chains in $P$. The geometric realisation of this simplicial complex coincides with $BP$.

\begin{center} \emph{Simplicial Complexes, Products and Internal Hom, Barycentric subdivision}
\end{center}
Simplicial complexes crop up throughout this paper.
We refer to Chapter 3,\cite{Spanier}, for the definition of a simplicial complex and its 
barycentric subdivision. $\sS(K)$ and
 $\sV(K)$ denote the sets of vertices and simplices respectively of a simplicial complex $K$.  
The geometric realisation of $K$ is denoted by $|K|$.  The set $\sS(K)$ is a partially ordered set
 (with respect to inclusion of subsets). Note that $B\sS(K)$ is simply the (geometric realisation of)
 the barycentric subdivision $sd(K)$.  The geometric  realisations of $K$ and $sd(K)$ are 
 canonically homeomorphic to each other, but not by a simplicial map. 
 
 Given simplicial complexes $K_1$ and $K_2$, the product $|K_1|\times |K_2|$ 
 (in the compactly generated topology) is canonically
 homeomorphic to $B(\sS(K_1)\times\sS(K_2))$.
\\ \\
The category of simplicial complexes and simplicial maps has a \emph{categorical product}:
\\ $\sV(K_1\times K_2)=\sV(K_1)\times \sV(K_2)$.  A non-empty subset of $\sV(K_1\times K_2)$
 is a simplex of $K_1\times K_2$ if and only if it is contained in $S_1\times S_2$ for
 some $S_i\in\sS(K_i)$  for $i=1,2$. The geometric realisation of the product is not 
 homeomorphic to the 
 product of the geometric realisations, but they do have the same homotopy type. In fact 
 Proposition ~\ref{covering} of section 1 provides a contractible collection of homotopy 
 equivalences $|K_1|\times|K_2|\to |K_1\times K_2|$.  For most purposes, it suffices 
 to note that there is a \emph{canonical} map  $P(K_1,K_2):|K_1|\times|K_2|\to |K_1\times K_2|$.  This is
  obtained in the following manner.  Let $C(K)$ denote the $\R$-vector space  with basis $\sV(K)$
  for a simplicial complex $K$. Recall that $|K|$ is a subset of $C(K)$.  For simplicial complexes
  $K_1$ and $K_2$, we have the evident isomorphism 
 \begin{center} $j:C(K_1)\otimes_{\R}C(K_2)\to C(K_1\times K_2).$
 \end{center}
 For  $c_i\in |K_i|$ for $i=1,2$ we put $P(K_1,K_2)(c_1,c_2)=
 j(c_1\otimes c_2)\in C(K_1\times K_2)$.  We note that   $j(c_1\otimes c_2)$ belongs to the subset
 $ |K_1\times K_2|$. This gives the  canonical  $P(K_1,K_2)$.
 \\ \\
 Given simplicial complexes $K,L$ there is a simplicial complex $\mathcal{H}om (K,L)$ with the following property: if $M$ is a simplicial complex, then the set of simplicial maps $K\times M\to L$ 
 is naturally identified with the set of simplicial maps $M\to\mathcal{H}om(K,L)$. This simple verification
  is left to the reader. 
\\ \\  
 Simplicial maps  $f:K_1\times K_2\to K_3$ occur in sections 2 and 5 of this paper.  
 \begin{center}$|f|\circ P(K_1,K_2):|K_1|\times |K_2|\to |K_3|$ 
 \end{center}
 is the map we employ on geometric realisations.  Maps $|K_1|\times|K_2|\to |K_3|$ associated to simplicial maps $f_1$
  and $f_2$ are  seen (by contiguity) to be homotopic to each other if $\{f_1,f_2\}$ is a simplex 
  of $\mathcal{H}om(K_1\times K_2,K_3)$. This fact is employed in Lemma~\ref{contiguity}. 
 
 Simplicial maps  $f:K_1\times K_2\to K_3$ 
 are in 
 reality maps $\sV(f):\sV(K_1)\times\sV(K_2)\to\sV(K_3)$  with the property that 
$\sV(f)(S_1\times S_2)$ is a simplex of $K_3$ whenever $S_1$ and $S_2$ are simplices 
of $K_1$ and $K_2$ respectively.   One should note that  such an $f$ induces a map
of posets $\sS(K_1)\times\sS(K_2)\to \sS(K_3)$, which in turn induces a continuous map
  $B(\sS(K_1)\times\sS(K_2))\to B\sS(K_3)$. In view of the natural identifications, 
  this is the same as giving a map $|K_1|\times|K_2|\to |K_3|$. This map coincides with
  the $|f|\circ P(K_1,K_2)$ considered above.  
  
   The homotopy assertion of maps
  $|K_1|\times|K_2|\to |K_3|$ associated to $f_1,f_2$ where $\{f_1,f_2\}$ is an edge of $\mathcal{H}om(K_1\times K_2, K_3)$  cannot be proved by the quick poset definition of the maps (for $|K_3|$ has been subdivided
   and contiguity is not available any more). This explains our preference for the longwinded 
   $|f|\circ P(K_1,K_2)$ definition.

\section{Some preliminaries from topology}

We work with the category of
compactly generated  weakly Hausdorff spaces. A good source is
Chapter 5 of \cite{Maybook}. 
This category possesses products.
It also possesses an internal
Hom in the following sense: for compactly generated Hausdorff $X,Y,Z$,  
continuous maps $Z\to\mathcal{H}om(X,Y)$ are the same as continuous maps
$Z\times X\to Y$, where $Z\times X$ denotes the product in this category.
This internal Hom property is required in the proof of Proposition~\ref{covering} stated below. 

$\mathcal{H}om(X,Y)$ is
the space of continuous maps from $X$ to $Y$. This space of maps
has the compact-open topology, which is then replaced by the inherited compactly generated
topology. This space $\mathcal{H}om (X,Y)$ is referred to frequently
 as $\mathrm{Map}(X,Y)$, and some times even as $\mathrm{Maps}(X,Y)$ ,                   in the text.

Now consider the following set-up. Let $\Lambda$ be a partially ordered 
set assumed to be Artinian: (i) every non-empty 
subset in $\Lambda$ has a minimal element with respect to the partial 
order, or equivalently (ii) there are no infinite strictly descending chains 
$\lambda_1>\lambda_2>\cdots$ in $\Lambda$. The poset  $\Lambda$
will remain fixed throughout the discussion below. 

We consider topological spaces $X$ equipped with a family of  closed subsets
 $X_{\lambda}, \lambda\in\Lambda$  with the property that 
 $X_{\mu}\subset X_{\lambda}$ whenever  $\mu\leq \lambda$.

Given another $Y,Y_{\lambda},\lambda\in\Lambda$ as above,
the collection of $\Lambda$-compatible continuous $f:X\to Y$  (i.e. 
satisfying $f(X_{\lambda}) \subset Y_{\lambda},\forall\lambda\in\Lambda$)
will be denoted by $Map_{\Lambda}(X,Y)$.  
$Map_{\Lambda}(X,Y)$ is a closed subset of $\mathcal{H}om(X,Y)$, and this
topologises $Map_{\Lambda}(X,Y)$. 

We say that $\{X_{\lambda}\}$ is a {\em weakly admissible covering} of 
$X$ if the three conditions listed below are satisfied.  It is an
{\em  admissible covering}
if in addition, each
$X_{\lambda}$ is contractible.
\begin{enumerate}
\item[(1)] For each pair of indices $\lambda,\mu\in \Lambda$, we have
\[X_{\lambda}\cap X_{\mu}=\mathop{\cup}_{\nu\leq\lambda,\nu\leq 
\mu} X_{\nu}\]
\item[(2)] If 
\[\partial X_{\lambda}=\mathop{\cup}_{\nu<\lambda}X_{\nu},\]
then 
$\partial X_{\lambda}\hookrightarrow X_{\lambda}$ 
is a cofibration
\item[(3)] 
 The topology on $X$ is coherent with respect to the family of 
subsets 
$\{X_{\lambda}\}_{\lambda\in \Lambda}$, 
that is, 
$X=\cup_{\lambda}X_{\lambda}$, 
and a subset $Z\subset X$ is 
closed precisely when $Z\cap X_{\lambda}$ is closed in $X_{\lambda}$ 
in the relative topology, for all $\lambda$.

\end{enumerate}

\begin{proposition}\label{covering}
Assume that $\{X_{\lambda}\}$ is a weakly admissible covering of $X$.
Assume also that each $Y_{\lambda}$ is contractible. 

Then the space ${\rm Map}_{\Lambda}(X,Y)$ of  
$\Lambda$-compatible maps $f:X\to Y$ is  
contractible.  In particular, it is non-empty and path-connected.

\end{proposition}
\begin{cor}\label{cor-covering} If both $\{X_{\lambda}\}$ and 
$\{Y_{\lambda}\}$ are admissible, then $X$ and $Y$ are homotopy 
equivalent.
\end{cor}
With assumptions as in the above corollary, the proposition yields the existence of
$\Lambda$-compatible maps
$f:X\to Y$ and $g:Y\to X$.  Because $g\circ f$ and $f\circ g$  are also $\Lambda$-compatible, 
that they are homotopic to $id_X$ and $id_Y$ respectively  is deduced from  the
path-connectivity of  ${\rm Map}_{\Lambda}(X,X)$  and ${\rm Map}_{\Lambda}(Y,Y)$.
 
\begin{cor}\label{cor-covering2} 
If $\{X_{\lambda}\}$ is admissible, then there is a  homotopy 
equivalence $X\to B\Lambda$.
\end{cor}
Here, recall that $B\Lambda$ is the geometric realization of the 
simplicial complex associated to the set of nonempty finite 
chains (totally ordered subsets) in $\Lambda$; equivalently, regarding 
$\Lambda$ as a category, $B\Lambda$ is the geometric realization of its 
nerve. We put $Y= B\Lambda$ and $Y_{\lambda}=B\{\mu\in\Lambda:\mu\leq \lambda\} $ 
in  Corollary~\ref{cor-covering} to deduce Corollary~ \ref{cor-covering2}. 

\emph{In both corollaries, what one obtains is a contractible collection of homotopy equivalences}; 
there is no preferred or `natural' choice. Naturally, this situation persists in all applications of the
above proposition and corollaries.

The proof of Proposition~\ref{covering} is easily reduced to the 
following extension lemma. 

\begin{lemma}\label{covering2}
Let $\{X_{\lambda}\},  \{Y_{\lambda}\}$ etc. be as in the above proposition.
Let $\Lambda'\subset 
\Lambda$ be a subset, with induced partial order, so that for any 
$\lambda\in \Lambda'$, $\mu\in \Lambda$  with $\mu\leq \lambda$, we have 
$\mu\in \Lambda'$. Let $X'=\cup_{\lambda\in \Lambda'}X_{\lambda}$, 
$Y'=\cup_{\lambda\in \Lambda'}Y_{\lambda}$. Assume given a continuous 
map $f':X'\to Y'$ with $f'(X_{\lambda})\subset Y_{\lambda}$ for all 
$\lambda\in \Lambda'$. Then $f'$ extends to a continuous map $f:X\to Y$ 
with $f(X_{\lambda})\subset Y_{\lambda}$ for all $\lambda\in \Lambda$.  
\end{lemma}
\begin{proof} Consider the collection of pairs $(\Lambda'',f'')$ satisfying:
\\(a) $\Lambda'\subset \Lambda''\subset \Lambda$
\\(b) $\mu\in \Lambda,\lambda\in\Lambda'',\mu\leq\lambda$ implies $\mu\in\Lambda''$
\\(c) $f'':\cup\{X_{\mu}|\mu\in\Lambda''\}\to Y$ is a continuous map
\\(d) $f''(X_{\mu})\subset Y_{\mu}$ for all $\mu\in\Lambda''$
\\(e) $f'|X_{\mu}=f''|X_{\mu}$ for all $\mu\in\Lambda'$  

This collection is partially ordered in a natural manner. 
The coherence condition on the topology of
$X$ ensures that every chain in this collection has an upper bound. 
The presence of  
$(\Lambda',f')$ shows that it is non-empty. By Zorn's lemma, there is a 
maximal element   $(\Lambda'',f'')$ in this collection. The Artinian hypothesis on
$\Lambda$ shows that  if  $\Lambda''\neq \Lambda$,  then 
its complement possesses a minimal element $\mu$.  Let $D''$ be the domain of
$f''$. The minimality of $\mu$ shows that  $D''\cap X_{\mu}=\partial X_{\mu}$.
By condition (d) above, we see that $f''(   \partial X_{\mu}   )$
 is contained in the contractible space $Y_{\mu}$. Because $ \partial X_{\mu}
 \hookrightarrow X_{\mu}$ is a cofibration, it follows that $f''| \partial X_{\mu}$ 
 extends to a map $g:X_{\mu}\to Y_{\mu}$. The $f''$ and $g$ patch together to
 give a continuous map $h:D''\cup X_{\mu}\to Y$. Since the pair
 $(\Lambda''\cup\{\mu\},h)$ evidently belongs to this collection, the maximality 
   of $(\Lambda'',f'')$ is contradicted. Thus $\Lambda''=\Lambda$ and this completes
    the proof.

\end{proof}
The proof of the Proposition follows in three standard steps. 
\\Step 1:
 Taking $\Lambda'=\emptyset$ in Lemma~\ref{covering2}
we deduce that  ${\rm Map}_{\Lambda}(X,Y)$  is nonempty.  
\\Step 2: For the 
path-connectivity of ${\rm Map}_{\Lambda}(X,Y)$ ,  we 
 replace $X$ by $X\times [0,1]$ and replace  
 the original poset $\Lambda$ by the product 
$\Lambda\times\{\{0\},\{1\},\{0,1\}\}$, with the product partial order, 
where the second factor is partially ordered by inclusion. The subsets of $X\times I$ 
 (resp.$Y$)
indexed by $(\lambda,0),(\lambda,1),(\lambda,\{0,1\})$ are $X_{\lambda}\times\{0\},  
X_{\lambda}\times\{1\}$ and $X_{\lambda}\times [0,1]$  (resp. $Y_{\lambda}$ in all three cases).

We then apply the lemma
to   the  sub-poset $\Lambda\times\{\{0\},\{1\}\}$. 
\\Step 3: Finally, for the contractibility of
${\rm Map}_{\Lambda}(X,Y)$, we first choose $f_0\in  {\rm Map}_{\Lambda}(X,Y)$ 
and then consider the two  maps 
 ${\rm Map}_{\Lambda}(X,Y)\times X\to Y$ given by $(f,x)\mapsto f(x)$
 and  $(f,x)\mapsto f_0(x)$.  Putting 
 $({\rm Map}_{\Lambda}(X,Y)\times X)_{\lambda}
 ={\rm Map}_{\Lambda}(X,Y)\times X_{\lambda}$
 for all $\lambda\in\Lambda$, we see that both the above maps are $\Lambda$-compatible.
 The path-connectivity assertion in Step 2 now gives a homotopy between the identity map
 of  ${\rm Map}_{\Lambda}(X,Y)$ and the constant map $f\mapsto f_0$. This
completes the proof of Proposition~\ref{covering}.

We now want to make some remarks about equivariant versions of the above 
statements. 

Given $X,\{X_{\lambda};\lambda\in\Lambda\}$ as above,
an action of a group $G$ on $X$ 
 is called $\Lambda$-compatible if $G$  also acts on the poset $\Lambda$
 so that  
for all $g\in G$, $\lambda\in \Lambda$, we have 
$g(X_{\lambda})=X_{g\lambda}$. 

\begin{sloppypar}
Under the conditions of Proposition~\ref{covering}, suppose $\{X_{\lambda}\}$, 
$\{Y_{\lambda}\}$ admit $\Lambda$-compatible $G$-actions. There is no
$G$-equivariant $f\in Map_{\Lambda}(X,Y)$ in general. However, if
 $f\in Map_{\Lambda}(X,Y)$ and $g_X,g_Y$ denote the actions
of $g\in G$ on $X$ and $Y$ respectively, we see that
, then $g_Y\circ f\circ g_X^{-1}$
is also a $\Lambda$-compatible map.
By  Proposition~\ref{covering},
we see that this map is homotopic to $f$. Thus $g_Y\circ f$ and
$f\circ g_X$ are homotopic to each other. \emph{In particular, 
$H_n(f):H_n(X)\to H_n(Y)$ is a homomorphism of $G$-modules.}
\end{sloppypar}
In the sequel a better version of this involving the Borel
construction is needed.

We recall the {\em Borel construction} of equivariant homotopy quotient 
spaces. Let $EG$ denote a contractible CW complex on which $G$ has a 
proper free cellular action; for our purposes, it suffices to fix a 
choice of this space $EG$ to be the geometric realization of the nerve 
of the translation category of $G$ (the category with vertices $[g]$ 
indexed by the elements of $G$, and unique morphisms between ordered 
pairs of vertices $([g],[h])$, thought of as given by the left action 
of  $hg^{-1}$). The classifying space $BG$ is the quotient space $EG/G$. 

If $X$ is any $G$-space, let $X//G$ denote the {\em homotopy quotient 
of $X$ by $G$}, obtained using the Borel construction, i.e., 
\begin{equation}\label{borel}
X//G=(X\times EG)/G,
\end{equation}
where $EG$ is as above, and $G$ acts diagonally. Note that the natural 
quotient map 
\[q_X:X\times EG\to X//G\] 
is a Galois covering space, with covering group $G$. 

If $X$ and $Y$ are $G$-spaces, then considering $G$-equivariant maps 
$\tilde{f}:X\times EG\to Y\times EG$ compatible with the projections to $EG$, 
giving a commutative diagram
$$\xymatrix{
X\times EG\ar[dr] \ar[rr]^{\tilde{f}}&& Y\times EG\ar[dl] \\
&EG&
}
$$
is equivalent to considering maps $\overline{f}:X//G\to Y//G$ 
compatible with the maps $q_X:X//G\to BG$, $q_Y:Y//G\to BG$, giving a 
commutative diagram
$$\xymatrix{
X//G\ar[dr]_{q_X} \ar[rr]^{\overline{f}}&& Y//G\ar[dl]^{q_Y} \\
&BG&
}
$$

\begin{proposition}\label{covering-equivariant}
Assume that, in the situation of proposition~\ref{covering}, there are 
$\Lambda$-compatible $G$-actions on $X$ and $Y$. Let $EG$ be as above, 
and consider the $\Lambda$-compatible families $\{X_{\lambda}\times 
EG\}$, which is a  weakly admissible covering family for $X\times EG$, 
and $\{Y_{\lambda}\times EG\}$, which is an admissible covering family 
for $Y\times EG$. Then there is a $G$-equivariant map 
$\tilde{f}:X\times EG\to Y\times EG$, compatible with the projections to 
$EG$, such that 
\begin{enumerate}
\item[(i)] $\tilde{f}(X_{\lambda}\times EG)\subset (Y_{\lambda}\times 
EG)$ for all $\lambda\in \Lambda$
\item[(ii)] if $\tilde{g}:X\times EG\to Y\times EG$ is another such 
equivariant map, then there is a $G$-equivariant homotopy between 
$\tilde{f}$ and $\tilde{g}$, compatible with the projections to $EG$
\item[(iii)] The space of such equivariant maps $X\times EG\to Y\times 
EG$, as in (i), is contractible.
\end{enumerate}
\end{proposition}
\begin{proof} We show the existence of the desired map, and leave the 
proof of other properties, by similar arguments, to the reader.

Let ${\rm Map}_{\Lambda}(X,Y)$ be the contractible space of
$\Lambda$-compatible 
maps from $X$ to $Y$; note that it comes equipped with a natural $G$-action, 
so that the canonical evaluation map $X\times {\rm Map}_{\Lambda}(X,Y)\to Y$ is 
equivariant. This induces $X\times {\rm Map}_{\Lambda}(X,Y)\times EG\to Y\times
EG$. 
There is also a natural $G$-equivariant map $\pi:X\times {\rm 
Map}_{\Lambda}(X,Y)\times EG\to X\times EG$. This map $\pi$ has 
equivariant sections, since the projection ${\rm 
Map}_{\Lambda}(X,Y)\times EG\to EG$ is a $G$-equivariant map between 
weakly contractible spaces, so that the map on quotients modulo $G$ is a 
weak homotopy equivalence (i.e., $({\rm Map}_{\Lambda}(X,Y)\times EG)/G$ 
is another ``model'' for the classifying space $BG=EG/G$). However $BG$ 
is a CW complex, so the map has a section.
\end{proof}

As another preliminary, we note some facts (see lemma~\ref{quillenA} 
below) which are essentially corollaries of Quillen's Theorem A (these 
are  presumably well-known to experts, though we do not have a specific 
reference). 

If $P$ is any poset, let $C(P)$ be the poset consisting of non-empty 
finite chains (totally ordered subsets) of $P$. If $f:P\to Q$ is a 
morphism between posets (an order preserving map) there is an induced 
morphism $C(f):C(P)\to C(Q)$. If $S$ is a simplicial complex 
(literally, a collection of finite non-empty subsets of the vertex set), 
we may regard $S$ as a poset, partially ordered with respect to inclusion; 
then the classifying space $BS$  is naturally homeomorphic to the 
geometric realisation $|S|$ (and gives the barycentric subdivision of 
$|S|$). A simplicial map $f:S\to T$ between simplicial complexes (that 
is, a map on vertex sets which sends simplices to simplices, not necessarily 
preserving dimension) is also then a morphism of posets. We say that a 
poset $P$ is contractible if its classifying space $BP$ is contractible.
\begin{lemma}\label{quillenA}
(i)\quad 
Let $f:P\to Q$ be a morphism between posets. Suppose that for each 
$X\in C(Q)$, the fiber poset $C(f)^{-1}(X)$ is contractible. Then $Bf:BP\to BQ$ 
is a homotopy equivalence.\\
(ii)\quad Let $f:S\to T$ be a simplicial map between simplicial 
complexes.  Suppose that for any simplex $\sigma\in T$, the fiber 
$f^{-1}(\sigma)$, considered as a poset, is contractible. Then 
$|f|:|S|\to |T|$ is a homotopy equivalence.
\end{lemma}
\begin{proof}
We first prove (i). For any poset $P$, there is morphism of posets 
$\varphi_P:C(P)\to P$, sending a chain to its first (smallest) element. 
If $a,b\in P$ with $a\leq b$, and $C$ is a chain in 
$\varphi_P^{-1}(b)$, then $\{a\}\cup C$ is a chain in 
$\varphi_P^{-1}(a)$. This gives an order preserving map of posets 
$\varphi_P^{-1}(b)\to \varphi_P^{-1}(a)$ (i.e., a ``base-change'' 
functor). This makes $C(P)$ prefibred over $P$, in the sense of Quillen 
(see page 96 in \cite{Srinivas}, for example). Also, $\varphi_P^{-1}(a)$ 
has the  minimal element (initial object) $\{a\}$, and so its 
classifying space is contractible. 

Hence Quillen's Theorem A (see \cite{Srinivas}, page 96) implies that 
$B(\varphi_P)$ is a homotopy equivalence, for any $P$.

Now let $f:P\to Q$ be a morphism between posets. Let $C(f):C(P)\to 
C(Q)$ be the corresponding morphism on the posets of (finite, nonempty) 
chains. If $A\subset B$ are two chains in $C(Q)$, there is an obvious 
order preserving map $C(f)^{-1}(B)\to C(f)^{-1}(A)$. Again, this makes 
$C(f):C(P)\to C(Q)$ prefibred. 

Since we assumed that $BC(f)^{-1}(A)$ is contractible, for all $A\in 
C(Q)$, Quillen's Theorem A implies that $BC(f)$ is a homotopy 
equivalence. 

We thus have a commutative diagram of posets and order preserving maps
\[\begin{array}{ccc}
C(P) & \stackrel{C(f)}{\to} & C(Q)\\
\varphi_P\downarrow\quad& &\quad\downarrow\varphi_Q\\
P& \stackrel{f}{\to} &Q
\end{array}
\]
where three of the four sides yield homotopy equivalences on passing to 
classifying spaces. Hence $Bf:BP\to BQ$ is a homotopy equivalence, 
proving (i). 

The proof of (ii) is similar. This is equivalent to showing that 
$Bf:BS\to BT$ is a homotopy equivalence. Since $f:S\to T$, regarded as a 
morphism of posets, is naturally prefibered, and by assumption,
$Bf^{-1}(\sigma)$ is contractible for each $\sigma\in T$, Quillen's 
Theorem A implies that $Bf$ is a homotopy equivalence.
\end{proof}

We make use of Propositions~\ref{covering} and 
\ref{covering-equivariant} in the following way. 

Let $A$, $B$ be sets, $Z\subset A\times B$  a subset such that the 
projections $p:Z\to A$, $q:Z\to B$ are both surjective.
 Consider simplicial complexes $S_Z(A)$, $S_Z(B)$ 
on vertex sets $A$, $B$ respectively, with simplices in $S_Z(A)$ being 
finite nonempty subsets of fibers $q^{-1}(b)$, for any $b\in B$, and 
simplices in $S_Z(B)$ being finite, nonempty subsets of 
fibers $p^{-1}(a)$, for any $a\in A$. 

Consider also a third simplicial complex $S_Z(A,B)$ with vertex set $Z$, 
where a finite non-empty subset $Z'\subset Z$ is a simplex if and only 
it satisfies the following condition:
\[(a_1,b_1),(a_2,b_2)\in Z' \;\Rightarrow (a_1,b_2)\in Z.\]
Note that the natural maps on vertex sets $p:Z\to A$, $q:Z\to B$ 
induce canonical simplicial maps on geometric realizations 
\[|p|:|S_Z(A,B)|\to |S_Z(A)|,\;\; |q|:|S_Z(A,B)|\to |S_Z(B)|.\]
 
\begin{cor}\label{comparison}\begin{enumerate}
\item[(1)] With the above notation, the simplicial maps  
\[|p|:|S_Z(A,B)|\to |S_Z(A)|,\;\; |q|:|S_Z(A,B)|\to |S_Z(B)|\]
are homotopy equivalences. In particular, $|S_Z(A)|$, $|S_Z(B)|$ are 
homotopy equivalent. 
\item[(2)]If a group $G$ acts on $A$ and on $B$, so that 
$Z$ is stable under the diagonal $G$ action on $A\times B$, then the 
homotopy equivalences $|p|$, $|q|$ are $G$-equivariant homotopy 
equivalences.   Hence there exists a  $G$-equivariant homotopy 
equivalence between $|S_Z(A)|\times EG$ and $|S_Z(B)|\times EG$.
\end{enumerate}
\end{cor}
\begin{proof} 
We first discuss (1). Since the situation is symmetric with respect 
to the sets $A$, $B$, it suffices to show $|p|$ is a homotopy equivalence. 

Let $\Lambda$ be the poset of all simplices of $S_Z(A)$, 
thought of as subsets of $A$, and ordered by inclusion. Clearly 
$\Lambda$ is Artinian.  

Apply Corollary~\ref{cor-covering} with $X=|S_Z(A,B)|$, $Y=|S_Z(A)|$, 
$\Lambda$ as above, and the following $\Lambda$-admissible coverings: 
for $\sigma\in \Lambda$, let $Y_{\sigma}$ be the (closed) simplex in 
$Y=|S_Z(A)|$ determined by $\sigma$ (clearly $\{Y_{\sigma}\}$ is 
admissible); take $X_{\sigma}=|p|^{-1}(Y_{\sigma})$ (this is evidently 
weakly admissible). For admissibility of $\{X_{\sigma}\}$, we need to 
show that each $X_{\sigma}$ is contractible. 

In fact, regarding the sets of simplices $S_Z(A,B)$ and $S_Z(A)$ as 
posets, and $S_Z(A,B)\to S_Z(A)$ as a morphism of posets, $X_{\sigma}$ 
is the geometric realization of the simplicial complex determined by 
$\cup_{\tau\leq \sigma} p^{-1}(\tau)$. 

The corresponding map of posets 
\[p^{-1}(\{\tau | \tau\leq \sigma\})\to \{\tau | 
\tau\leq \sigma\}\] 
has contractible fiber posets -- if we 
fix an element $x\in p^{-1}(\tau)$, and $p^{-1}(\tau)(\geq x)$ is the 
 sub-poset of elements bounded below by $x$, then $y\mapsto y\cup x$ is 
a morphism of posets $r_x:p^{-1}(\tau)\to p^{-1}(\tau)(\geq x)$ which 
gives a homotopy equivalence on geometric realizations (it is left 
adjoint to the inclusion of the sub-poset). But the sub-poset has a 
minimal element, and so its realization is contractible. 

The poset $\{\tau | \tau\leq \sigma\}$ is obviously contractible, since 
it has a maximal element. Hence, applying lemma~\ref{quillenA}(ii), 
$X_{\sigma}$ is contractible.

Since the map $|p|:X\to Y$ is $\Lambda$-compatible, corollary~\ref
{cor-covering} provides a contractible collection of $\Lambda$-compatible
 homotopy inverses of $|p|$. 

In the presence of a $G$-action, Proposition~\ref{covering-equivariant} 
provides a contractible family of $G$-equivariant maps from 
$|S_Z(A)|\times EG$ to $|S_Z(A,B)|\times EG$. This suffices to give (2).

\end{proof}

\section{Flag Spaces}
In this section, we discuss various constructions of spaces 
(generally simplicial complexes) defined using flags of free modules, 
and various maps, and homotopy equivalences, between these. These are 
used as building blocks in the proof of Theorem~\ref{theorem1}.

Let $A$ be a ring, and let $V$ be a free (left) $A$-module of rank 
$n$. Define a simplicial complex $\FL(V)$ as follows.  Its vertex set is 
$$\displaylines{
FL(V)=\{F=(F_0,F_1,\ldots,F_n)\mid 0=F_0\subset 
F_1\subset\cdots\subset F_n=V\mbox{ are}\hfill\cr\hfill\mbox{ $A$-submodules,
and each quotient $F_i/F_{i-1}$ is  $A$-free of rank 
1}\}.\cr}$$
We think of this vertex set as the set of ``full flags'' in $V$.

To describe the simplices in $\FL(V)$, we will need another definition. 
Let 
$$\displaylines{SPL(V)=\{\{L_1,\ldots,L_n\}\mid L_i\subset V\mbox{ is a 
free $A$-submodule of rank 1,}\cr\hfill\mbox{and the induced 
map $\oplus_{i=1}^nL_i\to V$ is an isomorphism}\}.\cr}$$
Note that $\{L_1,\ldots,L_n\}$ is regarded as an {\em unordered} set of 
free $A$-submodules of rank 1 of $L$ (i.e., as a subset of cardinality $n$ in 
the set of all free $A$-submodules of rank 1 of $V$). We think of 
$SPL(V)$ as the ``set of unordered splittings of $V$ into direct sums 
of free rank 1 modules''.

Given $\alpha\in SPL(V)$, say $\alpha=\{L_1,\ldots,L_n\}$, we may choose 
some ordering $(L_1,\ldots,L_n)$ of its elements, and thus obtain a full 
flag in $V$ (i.e., an element in $FL(V)$), given by 
\[(0,L_1,L_1\oplus L_2,\cdots,L_1\oplus\cdots\oplus L_n=V)\in 
FL(V).\]
Let 
\[[\alpha]\subset FL(V)\]
be the set of $n!$ such full flags obtained from $\alpha$. 

We now define a simplex in $\FL(V)$ to be any subset of such a set 
$[\alpha]$ of vertices, for any $\alpha\in SPL(V)$. Thus, $\FL(V)$ 
becomes a simplicial complex of dimension $n!-1$, with the sets 
$[\alpha]$ as above corresponding to maximal dimensional simplices.


Clearly ${\rm Aut}\,(V)\cong \GL_n(A)$ acts on the simplical 
complex $\FL(V)$ through simplicial automorphisms, and thus acts on the 
homology groups $H_*(\FL(V),\Z)$ (and other similar invariants of 
$\FL(V)$). 

Next, remark that if $F\in FL(V)$ is any vertex of $\FL(V)$, we may 
associate to it the free $A$-module $\gr_F(V)=\oplus_{i=1}^nF_i/F_{i-1}$. 
If $(F,F')$  is an ordered pair of distinct vertices, which are joined 
by an edge in $\FL(V)$, then we obtain a {\em canonical} isomorphism 
(determined by the edge) 
\[\varphi_{F,F'}:\gr_F(V)\stackrel{\cong}{\longrightarrow}\gr_{F'}(V).\]
One way to describe it is by considering the edge as lying in a simplex 
$[\alpha]$, for some $\alpha=\{L_1,\ldots,L_n\}\in SPL(V)$; this 
determines an identification of $\gr_F(V)$ with $\oplus_iL_i$, and a 
similar identification of $\gr_{F'}(V)$, and thereby an identification 
between $\gr_F(V)$ and $\gr_{F'}(V)$. Note that from this 
description of the maps $\varphi_{F.F'}$, it follows that if $F,F',F''$ 
form vertices of a 2-simplex in $\FL(V)$,  
i.e., there exists some $\alpha\in SPL(V)$ such that $F,F',F''\in 
[\alpha]$, then we also have 
\[\varphi_{F,F''}=\varphi_{F',F''}\circ \varphi_{F,F'}.\]

The isomorphism $\varphi_{F,F'}$ depends only on the (oriented) edge in 
$\FL(V)$ determined by $(F,F')$, and not on the choice of the simplex 
$[\alpha]$ in which it lies. One way to see this is to use that, for 
any two such filtrations $F$, $F'$ of $V$ 
there is a canonical isomorphism  $\gr^p_F\gr^q_{F'}(V)\cong 
\gr^q_{F'}\gr^p_F(V)$ (Schur-Zassenhaus 
lemma) for each $p,q$. But in case $F$, $F'$ are flags which are 
connected by an edge, then there is also a canonical isomorphism 
$\gr_F\gr_{F'}(V)\cong\gr_{F'}(V)$ (in fact the $F$-filtration induced 
on $\gr^p_{F'}(V)$ has only 1 non-trivial step, for each $p$), and 
similarly there is a canonical isomorphism 
$\gr_{F'}\gr_{F}(V)\cong\gr_F(V)$. These three canonical 
isomorphisms combine to give the isomorphism 
$\varphi_{F,F'}$. 
 
Hence there is a well-defined {\em local system} ${\bf gr}(V)$ of 
$A$-modules on the geometric realization $|\FL(V)|$ of the simplicial 
complex $\FL(V)$, whose fibre over a vertex $F$ is $\gr_F(V)$. 

Notice further that this local system ${\bf gr}(V)$ comes equipped with 
a natural ${\rm Aut}(V)$ action, compatible with the natural actions on 
$FL(V)$ and $\FL(V)$.  Indeed, any element $g\in {\rm Aut}(V)$ gives a 
bijection on the set of full flags $FL(V)$, with 
\[F=(0=F_0,F_1,\ldots,F_n=V)\in FL(V))\]
mapping to 
\[gF=(0=gF_0,gF_1,\ldots,gF_n=V).\]
This clearly gives an induced isomorphism $\oplus_i  F_i/F_{i-1}\cong 
\oplus_i gF_i/gF_{i-1}$, identifying the fibers of the local 
system over $F$ and $gF$ in a specific way. It is easy to see that 
if $\alpha=\{L_1,\ldots,L_n\}\in SPL(V)$, then 
$g\alpha=\{gL_1,\ldots,gL_n\}\in SPL(V)$, giving the action of ${\rm 
Aut}\,(V)$ on $SPL(V)$, so that if a pair of vertices $F,F'$ of $FL(V)$ 
lie on an edge contained in $[\alpha]$, then $gF,gF'$ lie on an edge 
contained in $[g\alpha]$, and so the induced identification 
$\varphi_{F,F'}$ is compatible with $\varphi_{gF,gF'}$.    
This induces the desired action of ${\rm Aut}(V)$ on the local system. 

Further, note that if $F,F'\in FL(V)$ are connected by an edge in 
$\FL(V)$, then $\varphi_{F,F'}$ is a direct sum of isomorphisms of the form 
\[{\rm gr}^iF(V)\to {\rm gr}^{\sigma(i)}_{F'}(V)\] 
between free modules of rank 1, 
where $\sigma$ is a permutation of $\{1,\ldots,n\}$.  
Thus, given any edge-path joining vertices 
$F,F'$ in $\FL(V)$, the induced composite isomorphism 
$\gr_F(V)\to \gr_{F'}(V)$ is again realized by such a direct sum of isomorphisms,
upto permuting the factors. In particular, given an edge-path loop based as 
$F\in FL(V)$, the induced automorphism of $\gr_F(V)$ is the composition of a 
``diagonal'' automorphism and a permutation.

Hence, the monodromy group of the local system ${\bf gr}(V)$ is clearly 
contained  in $N_n(A)$, defined as a semidirect product 
\begin{equation}\label{N_n(A)}
N_n(A)=(A^{\times}\times\cdots\times A^{\times})\semidirect S_n
\end{equation} 
where $S_n$ is the permutation group; we regard $N_n(A)$ as a subgroup 
of ${\rm Aut}\,(\oplus_{i}L_i)$ in an obvious way.

Now we make infinite versions of the above constructions.

Let $A^{\infty}$ be the set of sequences $(a_1,a_2,\ldots,a_n,\ldots)$ 
of elements of $A$, all but finitely many of which are 0, considered as 
a free $A$-module of countable rank. There is a standard inclusion 
$i_n:A^n\hookrightarrow A^{\infty}$ of the standard free $A$-module of 
rank $n$ as the submodule of sequences with $a_m=0$ for all $m>n$. The 
induced inclusion $i:A^n\to A^{n+1}$ is the usual one, given by 
$i(a_1,\ldots,a_n)=(a_1,\ldots,a_n,0)$. 

We may thus view $A^{\infty}$ as being given with a tautological flag, 
consisting of the $A$-submodules $i_n(A^n)$. We define a 
simplicial complex $\FL(A^{\infty})$, with vertex set 
$FL(A^{\infty})$ equal to the set of flags $0=V_0\subset V_1\subset 
V_2\subset \cdots\subset A^{\infty}$ where $V_i/V_{i-1}$ is a free 
$A$-module of rank $1$, for each 
$i\geq 1$, and with $V_n=i(A^n)$ for all sufficiently large 
$n$. Thus $FL(A^{\infty})$ is naturally the union of subsets 
bijective with $FL(A^n)$. To make $\FL(A^{\infty})$ into a 
simplicial complex, we define a simplex to be a finite set of vertices 
in some subset $FL(A^n)$ which determines a simplex in the simplicial 
complex $\FL(A^n)$; this property does not depend on the choice of $n$, 
since the natural inclusion $FL(A^n)\hookrightarrow 
FL(A^{n+1})$, regarded as a map on vertex sets, identifies 
$\FL(A^n)$ with a subcomplex of $\FL(A^{n+1})$, such that any simplex of 
$\FL(A^{n+1})$ with vertices in $FL(A^n)$ is already in the 
subcomplex $\FL(A^n)$. 

We consider $GL(A)\subset {\rm Aut}(A^{\infty})$ as the union of the 
images of the obvious maps $i_n:GL_n(A)\hookrightarrow {\rm 
Aut}\,(A^{\infty})$, obtained by automorphisms which fix all the basis 
elements of $A^{\infty}$ beyond the first $n$. We clearly have an 
induced action of $GL(A)$ on the simplicial complex $\FL(A^{\infty})$, 
and hence on its geometric realisation $|\FL(A^{\infty})|$ through 
homeomorphisms preserving the simplicial structure. The 
inclusion $\FL(A^n)\hookrightarrow \FL(A^{\infty})$ as a subcomplex
is clearly $GL_n(A)$-equivariant. 

Next, observe that there is a local system ${\bf gr}(A^{\infty})$ on 
$\FL(A^{\infty})$ whose fiber over a vertex 
$F=(F_0=0,F_1,\ldots,F_n,\dots)$ is $\gr_F(V)=\oplus_iF_i/F_{i-1}$. 
This has monodromy contained in 
\[N(A)=\cup_nN_n(A)\subset GL(A),\]
where we may also view $N(A)$ as the semidirect product of 
\[(A^{\times})^{\infty}=\mbox{diagonal matrices in $GL(A)$}\]
by the infinite permutation group $S_{\infty}$. This local system also 
carries a natural $GL(A)$-action, compatible with the $GL(A)$-action on 
$\FL(A^{\infty})$.

Next, we prove a property (Corollary~\ref{elementary}) about the 
action of elementary matrices on homology, which is needed later. The 
corollary follows immediately from the lemma below. 

For the statement 
and proof of  the lemma, we suggest that the reader browse the remarks on
$\mathcal{H}om(K_1\times K_2, K_3)$  in section 0, given simplicial complexes
$K_i$ for $i=1,2,3$. The notation $\elem(V'\hookrightarrow V'\oplus V'')$ that appears in the lemma has also been introduced in
section 0 under the heading ``elementary matrices".

\begin{lemma}\label{contiguity}
Let $V'$, $V''$ be two free $A$-modules of finite rank, $i':V'\to 
V'\oplus V''$, $i'':V''\to V'\oplus V''$ the inclusions of the 
direct summands. Consider the two natural maps
\begin{equation}\label{eq1}
\alpha,\beta:FL(V')\times FL(V'')\to FL(V'\oplus V'')
\end{equation}
given by
$$\displaylines{\alpha:\left((F'_1,\ldots,F'_r=V'),(F''_1,\ldots,
F''_s)\right)\mapsto\hfill\cr\hfill
(i'(F'_1),\ldots,i'(F'_r)=i'(V'),i'(V')+i''(F''_1),\ldots,
i'(V')+i''(F''_s)=V'\oplus 
V''),\cr}$$
$$\displaylines{\beta:\left((F'_1,\ldots,F'_r=V'),(F''_1,\ldots,
F''_s)\right)\mapsto\hfill\cr\hfill
(i''(F''_1),\ldots,i''(F''_s)=i''(V''),i''(V'')+i'(F'_1),\ldots,
i''(V'')+i'(F'_r)=V'\oplus 
V'').\cr}$$
\\(A) $\alpha$ and $\beta$ are vertices of a one-simplex of  
$\mathcal{H}om(\FL(V')\times \FL(V''), \FL(V'\oplus V''))$.
\\(B) The maps  $|\FL(V')|\times |\FL(V'')|\to |\FL(V'\oplus V'')|$ induced by $\alpha,\beta$
are homotopic to each other.
\\(C) Let $c:|\FL(V')|\times |\FL(V'')|\to |\FL(V'\oplus V'')|$ denote the map produced by $\alpha$. 
Denote the action of $g\in GL(V'\oplus V'')$ on   $|\FL(V'\oplus V'')|$ by  $|\FL(g)|$. Then 
$c$ and $|\FL(g)|\circ c$ are homotopic to each other, if 
 $g\in \elem(V' \hookrightarrow V' \oplus V'')\subset GL(V'\oplus V'')$.

\end{lemma}
\begin{proof} Part (A).  By the  definition of $\mathcal{H}om(K_1\times K_2,K_3)$ in section 0, we only have to check that $\alpha(\sigma'\times\sigma'')\cup\beta(\sigma'\times\sigma'')$ is a simplex 
of $FL(V'\oplus V'')$  for all simplices $\sigma'$ of $\FL(V')$ and all simplices $\sigma''$ of
 $\FL(V'')$.  Clearly it suffices to prove this for maximal simplices, so we assume that both
 $\sigma'$ and $\sigma''$ are maximal.
 
 Note that if we consider any maximal simplex $\sigma'$ in $\FL(V')$, 
it corresponds to a splitting $\{L'_1,\ldots,L'_r\}\in SPL(V')$. 
Similarly any maximal simplex $\sigma''$ of $\FL(V'')$ corresponds to 
a splitting 
$\{L''_1,\ldots,L''_s\}\in SPL(V'')$. 
This determines the splitting 
$\{i'(L'_1),\ldots,i'(L'_r),i''(L''_1),\ldots,i''(L''_s)\}$ 
of $V'\oplus V''$, giving rise to a maximal simplex $\tau$ of 
$\FL(V'\oplus V'')$, and clearly $\alpha(\sigma'\times\sigma'')$ and 
$\beta(\sigma'\times \sigma'')$ are  both contained in $\tau$.  Thus their union is a simplex.

(B) follows from (A). We now address (C).   We note that $c=g\circ c$ for all 
$g\in id+Hom_A(V'',V')$.
  Denoting by $d$ the map produced by
 $\beta$ we see that $d=g\circ d$ for all $g\in id+Hom_A(V',V'')$. Because $c,d$ are homotopic to each other, we see that $c$ and $g\circ c$ are in the same homotopy class when $g$ is in either of the
 two groups above. These groups generate $\elem(V'\hookrightarrow  V'\oplus V'')$,  and so this 
 proves (C).
\end{proof}
\begin{cor}\label{elementary}
(i) The group $E_{n+1}(A)$ of elementary matrices acts trivially 
on the image of the natural map
\[i_*:H_*(\FL(A^{\oplus n}),\Z)\to H_*(\FL(A^{\oplus n+1}),\Z).
\]
(ii) The action of the group $E(A)$ of elementary matrices on 
$H_*(\FL(A^{\infty}),\Z)$ is trivial. 
\end{cor}
\begin{proof} 
We put $V'=A^n$ and $V''=A$ in the previous lemma. The $c$ in part (B) of the lemma 
is precisely the $i$ being considered here. By (C) of the lemma, $g\circ i$ is homotopic to $i$
 for all $g\in \elem(A^n\hookrightarrow A^{n+1})=E_{n+1}(A)$. This proves (i). 
 The direct limit of the $r$-homology of $|\FL(A^n)|$, taken over all $n$, is the $r$-th homology
 of $|\FL(A^{\infty}|$. Thus (i) implies (ii).

\end{proof}

We will find it useful below to have other ``equivalent models'' of the 
spaces $\FL(V)$, $\FL(A^{\infty})$, by which we mean other simplicial 
complexes, also defined using collections of appropriate $A$-submodules, 
such that there are natural homotopy equivalences between the different 
models of the same homotopy type, compatible with the appropriate group 
actions, etc. 

We apply corollary~\ref{comparison} as follows. Let $V\cong A^n$. We put
  $A=SPL(V), B=FL(V)$ and $Z=\{\alpha,F):F\in[\alpha]\}$.  The simplicial complex
   $S_Z(B)$ of corollary~\ref{comparison} is $\FL(V)$ by its definition. The simplicial
    complex $S_Z(A)$ is our definition of $\SPL(V)$.  The homotopy equivalence of 
    $\SPL(V)$ and $\FL(V)$ follows from this corollary.

We define $SPL(A^{\infty})$ to be the collection of sets $S$ satisfying
\\(a) $L\in S$ implies that $L$ is a free rank one $A$-submodule of $A^{\infty}$,
\\ (b)$\oplus\{L:L\in S\}\to A^{\infty}$ is an isomorphism, and  
\\(c) the symmetric difference of $S$ and the standard collection:
\\ $\{A(1,0,0,...),A(0,1,0,...),\cdots\}$ is a finite set.

Corollary~\ref{comparison} is then applied to the subset $Z\subset SPL(A^{\infty})\times FL(A^{\infty})$
 consisting of the pairs $(S,F)$ so that there is a bijection $h:S\to\N$ so that  
for every $L\in S$, 
\\ $L\subset F_{h(L)}$ and $L\to  gr^F_{h(L)}$ is an isomorphism.   

The above $Z$ defines $\SPL(A^{\infty})$.  The desired homotopy 
equivalence of the geometric realisations of $\SPL(A^{\infty})$ and $\FL(A^{\infty})$ comes from
the same corollary.

We also find it useful to introduce a third model of the homotopy types 
of $\FL(V)$ and $\FL(A^{\infty})$, the ``enriched Tits buildings'' 
$\ET(V)$ and $\ET(A^{\infty})$. The latter is defined in the last remark of this section.

Let $V\cong A^n$ as a left $A$-module. Let 
$\sE(V)$ be the set consisting of ordered pairs 
\[(F,S)=\left((0=F_0\subset F_1\subset\cdots\subset 
F_r=V),(S_1,S_2,\ldots,S_r)\right),\]
where $F$ is a {\em partial flag} in $V$, which means that
$F_i\subset V$ is an $A$-submodule, such that $F_i/F_{i-1}$ is 
a nonzero free module for each $i$, and $S_i\in SPL(F_i/F_{i-1})$ is an 
unordered collection of free $A$-submodules of $F_i/F_{i-1}$ giving 
rise to a direct sum decomposition $\oplus_{L\in S_i}L\cong F_i/F_{i-1}$. 
Thus $S$ is a collection of splittings of the quotients $F_i/F_{i-1}$ 
for each $i$.

We may put a partial order on the set $\sE(V)$ in the following way: 
$(F,S)\leq (F',T)$ if the filtration $F$ is a refinement of $F'$, and 
the data $S$, $T$ of direct sum decompositions of quotients are 
compatible, in the following natural sense --- if 
$F'_{i-1}=F_{j-1}\subset F_j\subset\cdots F_l=F'_i$, 
then $T_i$ must be partitioned into subsets, which map to the sets 
$S_j,S_{j+1},\ldots,S_l$ under the appropriate quotient maps. In 
particular, $(F',T)$ has only finitely many possible predecessors 
$(F,S)$ in the partial order. 

We have a simplicial complex $ET(V):=N\sE(V)$, the nerve of the 
partially ordered set $\sE(V)$ considered as a category, so that 
simplices are just nonempty finite chains of elements of the vertex 
(po)set $\sE(V)$. 

Note that maximal elements of $\sE(V)$ are naturally identified with 
elements of $SPL(V)$, while minimal elements are naturally identified 
with elements of $FL(V)$. Simplices in $\FL(V)$ are nonempty finite 
subsets of $FL(V)$ which have a common upper bound in $\sE(V)$, and 
similarly simplices in $\SPL(V)$ are nonempty finite subsets of 
$SPL(V)$ which have a common lower bound in $\sE(V)$.

We now show that $\ET(V)=B\sE(V)$, the classifying space of the poset 
$\sE(V)$, is another model of the homotopy type of $|\FL(V)|$. 

In a 
similar fashion, we may define a poset $\sE(A^{\infty})$, and a space 
$\ET(A^{\infty})$, giving another model of the homotopy type of 
$|\FL(A^{\infty})|$.

We first have a lemma on classifying spaces of certain posets. 
For any poset $(P,\leq)$, and any $S\subset P$, let 
\[L(S)=\{x\in P|\; x\leq s\;\forall\;s\in S\},\;\;U(S)=\{x\in P|\; 
s\leq x\;\forall\;s\in S\}\]
be the upper and lower sets of $S$ in $P$, respectively. Let ${\mathcal
P}_{min}$ denote the simplicial complex with vertex set $P_{min}$ given by
minimal elements of $P$, and where a nonempty finite
subset $S\subset P_{min}$ is a simplex if $U(S)\neq \emptyset$. 
Let $|{\mathcal P}_{min}|$  denote the geometric realisation of ${\mathcal
P}_{min}$.
\begin{lemma}\label{poset1}
Let $(P,\leq)$ be a poset such that \\
(a) $\forall\;\;s\in P$, the set $L(\{s\})$ is finite\\
(b) if $\emptyset\neq S\subset P$ with $L(S)\neq \emptyset$, then the classifying
space  $BL(S)$ of $L(S)$ (as a subposet) is contractible. \\[2mm]
Then $|{\mathcal P}_{min}|$ is naturally homotopy equivalent to $BP$.
\end{lemma}
\begin{proof} 
We apply Proposition~\ref{covering}. Take 
\[\Lambda=\{L(S)|\; \emptyset \neq S\subset P\mbox{ and }L(S)\neq 
\emptyset\}.\]
This is a poset with respect to inclusion. All $\lambda\in \Lambda$ are finite
subsets of $P$, so $\Lambda$ is Artinian. By assumption, the subsets
$B(\lambda)\subset BP$, for $\lambda\in \Lambda$, are contractible.  
On the other hand, the sets $\lambda\cap P_{min}$ give simplices in $|{\mathcal
P}_{min}|$. Thus, both the spaces $BP$ and $|{\mathcal P}_{min}|$ have
$\Lambda$-admissible coverings, and are thus homotopy equivalent. 
\end{proof}

\begin{remark} 
  If a  poset $P$  has g.c.d. in the sense that $\emptyset\neq S\subset P$ and 
  $\emptyset\neq L(S)$ implies $L(S)=L(t)$ for some  $t\in P$, 
  then condition (b) of the lemma is immediately satisfied.  However $\sE(V)$ does not
 enjoy the latter property.
 
 For example, if $V=A^3$ with basis $e_1,e_2,e_3$, let $s=\{ Ae_1, Ae_2,Ae_3        \}$
 and 
 \\$t= \{Ae_1,A(e_1+e_2), Ae_3\}$ and let $S=\{s,t\}\subset  SPL(V)\subset \sE(V)$. 
 Then $L(S)$  has three minimal elements and two maximal elements.
 In particular, g.c.d. $(s,t)$ does not exist.  In this example, $B(L(S))$
 is an oriented graph in the shape  of the letter M.  
 
\end{remark}
\begin{proposition}\label{mainflags}
If $V$ is a free $A$ module of finite rank, the poset $\sE(V)$ satisfies the
hypotheses of lemma~\ref{poset1}. Thus, $|\FL(V)|$ is naturally homotopy
equivalent to $\ET(V)=B(\sE(V))$.
\end{proposition}
\begin{proof} 
Clearly the condition (a) of lemma~\ref{poset1} holds, so it suffices to prove
(b). 

We now make a series of observations.
\begin{enumerate}
\item[(i)] Regard $SPL(V)$ as the set of maximal elements of the poset 
$\sE(V)$. We observe that for any $s\in \sE(V)$, if $H(s)=SLP(V)\cap 
U(\{s\})$, then we have that $L(\{s\})=L(H(s))$. This is easy to see, 
once one has unravelled the definitions.

Thus, it suffices to show that for sets $S$ of the type 
$\emptyset \neq S\subset SPL(V)\subset \sE(V)$, we have that $B(L(S)$ is 
contractible. We assume henceforth that $S\subset SPL(V)$.

\item[(ii)] Given a submodule $W\subset V$ which determines a 
partial flag $0\subset W\subset V$, we have a natural inclusion of 
posets
\[\sE(W)\times\sE(V/W)\subset 
\sE(V),\]
where on the product, we take the partial order
\[(a_1,a_2)\leq(b_1,b_2)\Leftrightarrow \mbox{$a_1\leq b_1\in \sE(W)$ 
and $a_2\leq b_2\in \sE(V/W)$}.\]

Note that $\alpha\in \sE(V)$ lies in the sub-poset 
$\sE(W)\times\sE(V/W)$ precisely when $W$ is 
one of the terms in the partial flag associated to $\alpha$. Hence, if 
$\alpha$ lies in the sub-poset, so does the entire set $L(\{\alpha\})$.  

\item[(iii)] With notation as above, if $\emptyset \neq S\subset 
SPL(V)\subset \sE(V)$, and $L(S)$ has nonempty intersection with the 
image of $\sE(W)\times\sE(V/W)\subset \sE(V)$,
then clearly there exist nonempty subsets $S'(W)\subset SPL(W)$,
$S''(W)\subset \SPL(V/W)$ such that  
\begin{equation}\label{remark4}
L(S)\cap (\sE(W)\times \sE(V/W))=
 L(S'(W))\times L(S''(W))
\end{equation}

\item[(iv)] If $\emptyset\neq S\subset SPL(V)$, then each $s\in S$ is a 
subset of 
\[\sL(V)=\{L\subset V| L\mbox{ is a free direct summand of rank 1 of 
$V$}\},\]
the set of lines in $V$. Let 
\[T(S)=\displaystyle{\mathop{\cap}}_{s\in S}s=\mbox{ lines common to all 
members of $S$},\]  
so that $T(S)\subset\sL(V)$. Let $M(S)$ denote the direct sum of the 
elements of $T(S)$, so that $M(S)$ is a free $A$-module of finite rank, 
and $0\subset M(S)\subset V$ is a partial flag, in the sense explained 
earlier; further, $T(S)$ may be regarded also as an element of 
$SPL(M(S))\subset \sE(M(S))$.

\item[(v)] We now claim the following: if $\emptyset\neq S\subset 
SPL(V)$ and $b\in L(S)$, then there exists a unique subset $f(b)\subset 
T(S)$ such that if $M(b)$ is the (direct) sum of the lines in 
$f(b)$, then 
\[b=(f(b),b')\in SPL(M(b))\times 
\sE(V/M(b))\subset \sE(M(b))\times\sE(V/M(b))\subset \sE(V).\]

Indeed, if  
\[b =((0=W_0\subset W_1\subset\cdots\subset W_h=V),(t_1,t_2,\ldots,t_h)\}\]
where $t_i\in SPL(W_i/W_{i-1})$, then since $b\in L(S)$, we must have 
that $t_1 \subset s$ for all $s\in S$, which implies that $t_1\subset 
T(S)$. Take $M(b)=W_1$, $t_1=f(b)\in SPL(M(b))$. 

 Let ${\mathcal P}(T(S))$ be the poset of nonempty subsets of $T(S)$, with 
respect to inclusion. Then $b\mapsto f(b)$  
  gives an order-preserving map  
  $f:L(S)\to {\mathcal P}(T(S))$.

\item[(vi)] For any $b\in L(S)$, we have 
\[L(b)=L(\{f(b)\})\times L(b')\subset \sE(M(b))\times\sE(V/M(b)),\]
so that if $b_1\in L(b)$, then $\emptyset\neq f(b_1)\subset f(b)$.
\end{enumerate}

We will now complete the proof of Proposition~\ref{mainflags}. We proceed by
induction on the rank of $V$. Suppose $S\subset SPL(V)$ is nonempty, and
$L(S)\neq \emptyset$. 

If $M(S)=V$, then $S=\{s\}$ for some $s$, and $L(S)=L(\{s\})$ is a cone, hence
contractible. So assume $M(S)\neq V$.

If $T\subset T(S)$ is non-empty, and $M(T)\subset V$ the (direct) sum of 
the lines in $T$, then in the notation of (\ref{remark4}) above, with 
$W=M(T)$, we have $S'(W)=\{T\}$, and so $f^{-1}(T)=\{T\}\times 
L(S''(W)))$ for some $S''(W)\subset SPL(V/W)$. 

Now by induction, we have that $L(S''(W))$ is contractible, provided it 
is non-empty. Hence the non-empty fiber posets of $f$ are contractible. 
If $\emptyset\neq T\subset T'\subset T(S)$, then there is a morphism of 
posets $f^{-1}(T)\to f^{-1}(T')$ given as follows: if $b\in f^{-1}(T)$, 
and 
\[b =((0=W_0\subset W_1\subset\cdots\subset W_h=V),(t_1,t_2,\ldots,t_h)\}\]
where $t_i\in SPL(W_i/W_{i-1})$, then since $b\in f^{-1}(T)$, we must 
have $t_1=T$, $W_1=M(T)$. Now define $b'\in f^{-1}(T')$ using the 
partial flag 
\[0=W'\subset W_1+M(T')\subset W_2+M(T')\subset \cdots\subset 
W_h+M(T')=V\]
and elements $t_i'\in SPL(W_i+M(T')/W_{i-1}+M(T'))$ induced by the 
$t_i$. This is easily seen to be well-defined, and gives a morphism of 
posets $f^{-1}(T)\to f^{-1}(T')$. In particular, if $T\subset T'\subset 
T(S)$ and $f^{-1}(T)$ is non-empty, then so is $f^{-1}(T')$.  

Now take any $b\in L(S)$ and put $T=f(b), T'=T(S)$ in the above to
deduce that $f^{-1}(T(S))\neq\emptyset$. By (iii) above, we see that every 
$f^{-1}X$ is nonempty (and therefore contractible as well)
 for every nonempty 
$X\subset T(S)$. 

 We see that all the fiber posets $f^{-1}(T)$ considered above 
are nonempty.

This makes $f$ pre-cofibered, in the sense of Quillen 
(see \cite{Srinivas}, page 96), with contractible fibers. Hence by 
Quillen's Theorem A, $f$ induces a homotopy equivalence on classifying 
spaces. But ${\mathcal P}(T(S))$ is contractible (for example, since 
$T(S)$ is the unique maximal element).  
\end{proof}
\begin{remark}
Proposition~\ref{covering-equivariant} and the remarks preceding it
apply to the above Proposition. In particular, we obtain homotopy equivalences
$f:\ET(V)\to\FL(V)$ so that the induced maps on homology are $GL(V)$-equivariant.

\end{remark}

\begin{remark} We now define the poset $\sE(A^{\infty})$ and show that
 $\ET(A^{\infty})=B\sE(A^{\infty})$ is homotopy equivalent to $|\FL(A^{\infty}|$. 

We have already observed that a short exact sequence of free modules of finite rank
\begin{center}
$0\to V'\to V\to V''\to 0$
\end{center}
induces a natural inclusion $\sE(V')\times \sE(V'')\hookrightarrow \sE(V)$ of posets. In particular, 
when $V''\cong A$, this yields an inclusion $\sE(V')\hookrightarrow\sE(V)$.

We have $...\subset A^n\subset A^{n+1}\subset...\subset A^{\infty}$ as in the definition of 
$\FL(A^{\infty})$.  From the above, we obtain a direct system of posets 
\begin{center}$...\sE(A^n)\hookrightarrow \sE(A^{n+1})\hookrightarrow...$
\end{center}
and we define $\sE(A^{\infty})$ to be the direct limit of this system of posets. 

We put $P=\sE(A^{\infty})$ in lemma~\ref{poset1}.  We note that $\alpha\leq\beta, \alpha\in
\sE(A^{\infty}),\beta\in\sE(A^n)$ implies that $\alpha\in\sE(A^n)$.  It follows that 
$P=\sE(A^{\infty})$ satisfies the requirements of the lemma because each $\sE(A^n)$ does.
 It is clear that $P_{min}=
 FL(A^{\infty})$, and furthermore that $\sP_{min}=\FL(A^{\infty})$.  This yields the 
 homotopy equivalence of $|\FL(A^{\infty})|$ with $\ET(A^{\infty})$. 
 
 By Proposition~\ref{covering-equivariant}, it follows that $\ET(A^{\infty})//GL(A)$ and 
 $|\FL(A^{\infty})|//GL(A)$ are also homotopy equivalent to each other.

  It has already been 
 remarked that Corollary~\ref{comparison} gives the homotopy equivalence of 
 $|\SPL(A^{\infty}|$ with $|\FL(A^{\infty})|$.  Combined with Proposition~\ref{covering-equivariant},
  this gives the homotopy equivalence of 
  $|\SPL(A^{\infty}|//GL(A)$ with $|\FL(A^{\infty})|//GL(A)$. The remarks 
 preceding that proposition, combined with corollary ~\ref{elementary}, show that 
 the action of $E(A)$ on the homology groups of $|\SPL(A^{\infty})|$ is trivial.

\end{remark}

\section{homology of the Borel constrcuction}

Let $V$ be a free $A$-module of rank $n$.  Fix $\beta\in SPL(V)$ and let $N(\beta)\subset GL(V)$
 be the stabiliser of $\beta$ (when $V=A^n$ and $\beta$ is the standard splitting, then 
 $N(\beta)$ is the subgroup $N_n(A)$ of the last section).  That there is a 
 $GL(V)$-equivariant  $N(\beta)$-torsor on $|\FL(V)|$  has been observed in the previous section.
  In a similar manner, one may construct a  $GL(V)$-equivariant  $N(\beta)$-torsor on $\ET(V)$.
This gives rise to a $N(\beta)$-torsor on $\ET(V)//GL(V)$. Because $BN(\beta)$ is a classifying
 space for such torsors, we obtain a map $\ET(V)//GL(V)\to BN(\beta)$, well defined up to homotopy. 

\begin{sloppypar}
On the other hand,  the inclusion of $\beta$ in $\ET(V)$ gives rise to an inclusion
$BN(\beta)=\{\beta\}//N(\beta)\hookrightarrow \ET(V)//GL(V)$. It is clear that the composite 
$BN(\beta)\hookrightarrow \ET(V)//GL(V)\to BN(\beta)$ is homotopic to the identity.  Thus
 $BN(\beta)$ is a homotopy retract of  $\ET(V)//GL(V)$, but not homotopy equivalent to 
  $\ET(V)//GL(V)$. Nevertheless
  we have the following statement: 
\end{sloppypar}
\begin{proposition}\label{suslin-lemma2}
The map  $BN(\beta)\to \ET(V)//GL(V)$
induces an isomorphism on integral homology, provided $A$ is as in 
theorem~\ref{theorem1}.
\end{proposition}
\begin{proof}
Fix a basis for $V$, identifying $GL(V)$ with $GL_n(A)$. Let $\beta\in SPL(V)$
be the element naturally determined by this basis. Regarded as a vertex of 
$\ET(V)$, let $(\beta,*)\mapsto \overline{\beta}$ under the natural
map
\[\pi:\ET(V)//GL(V)\to \ET(V)/GL(V)\]
from the homotopy quotient to the geometric quotient, where $*\in EGL(V)$ is the
base point (corresponding to the vertex labelled by the
identity element of $GL(V)$).

For any $x\in \ET(V)$, let $\sH(x)\subset GL(V)$ be the isotropy group of $x$
for the $GL(V)$-action on $\ET(V)$. Note that since 
\[
 \ET(V)//GL(V)=\left(\ET(V)\times EGL(V)\right)/GL(V),
\]
the fiber $\pi^{-1}(\pi((x,*))$ may be identified
with $EGL(V)/\sH(x)$, which has the homotopy type of $B\sH(x)$.

In particular, the fiber $\pi^{-1}(\overline{\beta})$ has the homotopy type of
$BN_n(A)$, Further, the principal $N_n(A)$ bundle on $EGL(V)/\sH(\beta)$ is
naturally identified with the universal $N_n(A)$-bundle on $BN_n(A)$ -- its
pullback to $\{\beta\}\times EGL(V)$ is the trivial $N_n(A)$-bundle, regarded as
an $N_n(A)$-equivariant principal bundle, where $N_n(A)$ acts on itself (the
fiber of the trivial bundle) by translation. This means that the composite
\[\pi^{-1}(\overline{\beta})\to \ET(V)/GL(V)\to BN_n(A)\]
is a homotopy equivalence, which is homotopic to the identity, if we identify
$EGL(V)/\sH(\beta)$ with $BN_n(A)$. 

Thus, the lemma amounts to the assertion that 
{\em $\pi^{-1}(\overline{\beta})\to \ET(V)//GL(V)$ induces an isomorphism in
integral homology}.

Fix $\alpha\in FL(V)$ with $\alpha\leq \beta$ in the poset $\sE(V)$. Let 
\[P=\{\lambda\in \sE(V)|\alpha\leq \lambda\leq \beta\}.\]
One sees easily that (i) $BP$ is contractible, and (ii) the map $BP\to 
\ET(V)/GL(V)$ is a homoemorphism. The first assertion is obvious, since $P$
has a maximal (as well as a minimal) element, so that $BP$ is a cone. For the
second assertion, we first note that an element $b\in P\subset \sE(V)$ is
uniquely determined by the ranks of the modules in the partial flag in $V$
associated to $b$. Conversely, given any increasing sequence of numbers
$n_1<\ldots<n_h=\rank V$, there does exist an element of $P$ whose 
partial flag
module ranks are these integers. Given any element $b\in \sE(V)$, there exists
an element $g\in GL(V)$  
so that $g(b)=b'\in P$; the element $b'$ is the unique one determined by the
sequence of ranks associated to $b$. Finally, one observes that if $b\in P$, and
$g\in GL(V)$ such that $g(b)\in P$, then in fact $g(b)=b$: this is a consequence
of the uniqueness of the element of $P$ with a given sequence of ranks. These
observations imply that $BP\to \ET(V)/GL(V)$ is bijective; it is now easy to
see that it is a homeomorphism.
 
We may view $\ET(V)//GL(V)$ as the quotient of $BP\times EGL(V)$ by the
equivalence relation 
\begin{equation}\label{eqreln}
(x,y)\sim (x',y') \Leftrightarrow \mbox{ $x=x'$, and $y'=g(y)$ for some 
$g\in
\sH(x)$}.
\end{equation}

The earlier map $\pi:\ET(V)//GL(V)\to \ET(V)/GL(V)$ may be viewed now as the
map induced by the projection $BP\times EGL(V)\to BP$. We may, with this
identification, also identify $\overline{\beta}$ with $\beta$.

Next, we construct a ``good'' fundamental system of open neighbourhoods of an
arbitrary point $x\in BP$, which we need below. Such a point $x$ lies in the
relative interior of a unique simplex $\sigma(x)$ (called the {\em carrier} of
$x$) 
corresponding  to a chain $\lambda_0<\lambda_1<\cdots<\lambda_r$. Then one sees
that the stabiliser $\sH(x)\subset GL(V)$ is given by
\[\sH(x)=\bigcap_{i=0}^r\sH(\lambda_i),\]
since any element of $GL(V)$ which stabilizes the simplex $\sigma(x)$ must
stabilize each of the vertices (for example, since the $GL(V)$ action preserves
the partial order).

Let $\sta(x)$ be the union of the relative interiors of all simplices in $BP$
containing $\sigma(x)$ (this includes the relative interior of $\sigma(x)$ as
well, so it contains $x$). It is a standard property of simplicial complexes
that $\sta(x)$ is an open neighbourhood of $x$ in $BP$. Then if $z\in\sta(x)$,
clearly $\sigma(z)$ contains $\sigma(x)$, and so $\sH(z)\subset \sH(x)$. 

Next, for such a point $z$, and any $y\in EGL(V)$, it makes sense to consider
the path 
\[t\mapsto (tz+(1-t)x,y)\in \sigma(z)\times EGL(V)\subset BP\times EGL(V)\]
(where we view the expression $tz+(1-t)x$ as a point of $\sigma(z)$, using the
standard barycentric coordinates). In fact this path is contained in
$\sta(x)\times\{y\}$, and gives a continuous map
\[H(x):\sta(x)\times EGL(V)\times I\to \sta(x)\times EGL(V)\]
which exhibits $\{x\}\times EGL(V)$ as a strong deformation retract of
$\sta(x)\times EGL(V)$. Further, this is compatible with the equivalence
relation $\sim$ in (\ref{eqreln}) above, so that we obtain a strong deformation
retraction 
\[\overline{H(x)}:\pi^{-1}(\sta(x))\times I\to\pi^{-1}(\sta(x)).\]

In a similar fashion, we can construct a fundamental sequence of open
neighbourhoods 
$U_n(x)$ of $x$ in $BP$, with $U_1(x)=\sta(x)$, and set 
\[U_n(x)=H(x)(\sta(x)\times EGL(V)\times [0,1/n)).\] 
The same deformation retraction $H$ determines, by reparametrization, a
deformation retraction
\[H_n(x):\pi^{-1}(U_n(x))\times I\to \pi^{-1}(U_n(x))\]
of $\pi^{-1}(U_n(x)$ onto $\pi^{-1}(x)$.

Thus, if $P'=P\setminus{\beta}$, then 
\[\pi^{-1}(\sta(\beta))=\ET(V)//GL(V)\setminus \pi^{-1}(BP'),\]
and from what we have just shown above, the inclusion
\[\pi^{-1}(\beta)\to \pi^{-1}(\sta(\beta))=\ET(V)//GL(V)\setminus
\pi^{-1}(BP')\]
is a homotopy equivalence. To simplify notation, we let $X=\ET(V)//GL(V)$, so
that we  have the map $\pi:X\to BP$, and $X^0=X\setminus \pi^{-1}(BP')$.  
Let $\pi^0=\pi\mid_{X^0}:X^0\to BP$.

We are reduced to showing, with this notation, that the inclusion of the (dense)
open subset 
\[X^0\to X\]
induces an isomorphism in integral homology. Equivalently, it suffices to show
that this inclusion induces an isomorphism on cohomology 
with arbitrary constant coefficients $M$. By the Leray spectral sequence, this
is a consequence of showing that the maps of sheaves 
\[R^i\pi_*M_X\to R^i\pi^0_*M_{X^0}\]
is an isomorphism, which is clear on stalks $x\in BP\setminus BP'$. Now
consider stalks at a point $x\in BP'$. For any point $x'\in \sta(x)$, 
note that $x$ lies in some face of $\sigma(x')$ (the carrier of $x'$). 
We had defined a fundamental system of neighbourhoods $U_n(x)$ of $x$ in $BP$; 
explicitly we have
\[U_n(x)=\{tx'+(1-t)x|0\leq t<1/n\mbox{ and $x'\in\sta(x)$}\}.\]
Here, as before, we make sense of the above expression $tx'+(1-t)x$ using
barycentric coordinates in $\sigma(x')$. 

Define 
\[z_n(x)=\frac{1}{2n}\beta+(1-\frac{1}{2n})x.\]
Note that $z\in BP\setminus BP'=\sta(\beta)$. Further, observe that 
$U_n(x)\cap BP\setminus BP'$ is contractible, contains the point $z$, and for
any $w\in U_n(x)\cap BP\setminus BP'$, contains the line segment 
joining $z$ and $w$
(this makes sense, in terms of barycentric coordinates of any simplex 
containing both $z_n(x)$ and $w$; this simplex is either the carrier of 
$w$, or the cone over it with 
vertex $\beta$, of which $\sigma(w)$ is a face). 

This implies $\sH(w)\subset \sH(z_n(x))=\sH(x)\cap \sH(\beta)$, for all $w\in
U_n(x)$.  A minor modification of the proof (indicated above) that 
$\pi^{-1}(x)\subset \pi^{-1}(U_n(x))$ is a strong deforamtion retract, 
yields
the statement that \[\pi^{-1}(z_n(x))\to \pi^{-1}(U_n(x)\setminus 
BP')\] is a strong deformation retract. Hence, the desired isomorphism  
on stalks follows from:
\begin{equation}\label{suslin-lemma3}
B(\sH(x)\cap \sH(\beta))\to B(\sH(x))\mbox{ induces isomorphisms in integral
homology.}
\end{equation}

We now show how this statement, for the appropriate rings $A$, is reduced to
results of \cite{NS}. 

First, we discuss the structure of the isotropy groups $\sH(x)$ encountered
above. Let $\lambda\in P$, 
given by 
\[\lambda=(F,S)=\left((0=F_0\subset F_1\subset\cdots\subset 
F_r=V),(S_1,S_2,\ldots,S_r)\right),\]
where we also have $\alpha\leq \lambda\leq \beta$ for our chosen elements 
$\alpha\in FL(V)$ and $\beta\in SPL(V)$. We may choose a basis for each of the
lines in the splitting
$\beta$; then $\alpha\in FL(V)$ uniquely determines an order among these basis
elements, and thus a basis for the underlying free $A$-module $V$, such that 
the $i$-th submodule in the full flag $\alpha$ is the submodule generated by 
the first $i$ elements in $\beta$. Now the stabilizer $\sH(\alpha)$ may be viewed
as the group of upper triangular matrices in $GL_n(A)$, while $\sH(\beta)$ is 
the group generated by the diagonal subgroup in $GL_n(A)$ and the group of 
permutation matrices, identified with the permutation group $S_n$. 

In these terms, $\sH(\lambda)$ has the following structure. The 
filtration $F=(0=F_0\subset F_1\subset\cdots\subset F_r=V)$ is a sub-filtration of the 
full flag $\alpha$, and so determines a ``unipotent subgroup'' $U(\lambda)$
of elements fixing the elements of this partial flag, and acting trivially on the 
graded quotients $F_i/F_{i-1}$. These are represented as matrices of the 
form
\[
\left[\begin{array}{ccccc}
I_{n_1} & * & * & \cdots & * \\
0 & I_{n_2} & * & \cdots & * \\
0 & 0 & I_{n_3} & \cdots & *\\   
& \vdots& \ddots &&\vdots \\
0 & \cdots && & I_{n_r}
\end{array}
\right]
\]
where $n_i=\rank(W_i/W_{i-1})$, $I_{n_i}$ is the identity matrix of size 
$n_i$; these are the matrices which are strictly upper triangular with 
respect to a certain ``ladder''. Next, we may consider the group 
$S(\lambda)\subset S_n$ of permutation matrices, supported within the 
corresponding diagonal blocks, of the form 
\[
\left[\begin{array}{ccccc}
A_1 & 0 & 0 & \cdots & 0 \\
0& A_2 & 0 & \cdots & 0 \\
0 & 0 & A_3 & \cdots & 0\\   
& \vdots& \ddots &&\vdots \\
0 & \cdots && & A_r
\end{array}
\right]
\]
where each $A_j$ is a permutation matrix. Finally, we have the diagonal 
matrices $T_n(A)\subset GL_n(A)$, which are contained in $\sH(\lambda)$ 
for any such $\lambda$. In fact $\sH(\lambda)=U(\lambda)T_n(A) 
S(\lambda)$, where the group $T_n(A)S(\lambda)$ normalizes the subgroup 
$U(\lambda)$, making $\sH(\lambda)$ a semidirect product of 
$U(\lambda)$ and $T_n(A)S(\lambda)$. We also have that $S(\lambda)$ normalizes 
$U(\lambda)T_n(A)$.
   
In particular, $\sH(\alpha)$ has trivial associated 
permutation group $S(\alpha)=\{I_n\}$, while $\sH(\beta)$ has trivial 
unipotent group $U(\beta)=\{I_n\}$ associated to it.

Now if $x\in BP$, and $\sigma(x)$ is the simplex associated to  the 
chain $\lambda_0<\cdots<\lambda_r$ in the poset $P$, then it is 
easy to see that $\sH(x)$ is the semidirect product of $U(x):=U(\lambda_r)$ 
and $T_n(A)S(x)$, with $S(x):=S(\lambda_0)$, since as seen earlier, $\sH(x)$ is the 
intersection of the $\sH(\lambda_i)$. In other words, the ``unipotent 
part'' and the ``permutation group'' associated to $\sH(x)$ are each the 
smallest possible ones from among the corresponding groups  attached to 
the vertices of the carrier of $x$. Again we have that $S(x)$ normalizes $U(x)T_n(A)$. 

We return now to the situation in (\ref{suslin-lemma3}). We see that the groups 
$\sH(x)=U(x)T_n(A)S(x)$ and $\sH(x)\cap \sH(\beta)=T_n(A)S(x)$ both have the same associated permutation 
group $S(x)$, which normalizes $U(x)T_n(A)$ as well as $T_n(A)$. By comparing the spectral sequences
\[E^2_{p,q}=H_p(S(x),H_q(U(x)T_n(A),\Z))\Rightarrow H_{p+q}(\sH(x),\Z),\]
\[E^2_{p,q}=H_p(S(x),H_q(T_n(A),\Z))\Rightarrow H_{p+q}(\sH(x)\cap \sH(\beta),\Z)\]
we see that it thus suffices to show that the inclusion 
\begin{equation}\label{suslin-lemma4}
T_n(A)\subset U(x)T_n(A)
\end{equation}
induces an isomorphism on integral homology. 

Now lemma~\ref{suslin-lemma5} below finishes the proof.
\end{proof}

To state lemma~\ref{suslin-lemma5} we use the following notation. Let 
$I=\{i_0=0<i_1<i_2<\cdots<i_r=n\}$ be a subsequence of 
$\{0,1,\ldots,n\}$, so that $I$ determines a partial flag 
\[0\subset A^{i_1}\subset A^{i_2}\subset \cdots\subset A^{i_r}=A^n,\]
where $A^j\subset A^n$ as the submodule generated by the first $j$ basis 
vectors. Let $U(I)$ be the ``unipotent'' subgroup of $GL_n(A)$ 
stabilising this flag, and acting trivially on the associated graded 
$A$-module, and let $G(I)\subset GL_n(A)$ be the subgroup generated by 
$U(I)$ and $T_n(A)=(A^{\times})^n$, the subgroup of diagonal matrices. 
Then $T_n(A)$ normalises $U(I$), and $G(I)$ is the semidirect product of 
$U(I)$ and $T_n(A)$.

\begin{lemma}\label{suslin-lemma5}
Let $A$ be a Nesterenko-Suslin ring. For any $I$ 
as above, the homomorphism $G(I)\to G(I)/U(I)\cong T_n(A)$ induces an 
isomorphism on integral homology $H_*(G(I),\Z)\to H_*(T_n(A),\Z)$.
\end{lemma}
\begin{proof}
We work by induction on $n$, where there is nothing to prove when 
$n=1$, since we must have $G(I)=T_1(A)=A^{\times}=GL_1(A)$.  Next, if 
$n>1$, and $I=\{0<n\}$, then $U(I)$ is the trivial group, so there is 
nothing to prove. 

Hence we may assume $n>1$, $r\geq 2$, and thus 
$0<i_1<n$. There is then a natural homomorphism $G(I)\to G(I')$, where 
$I'=\{0<i_2-i_1<\cdots <i_r-i_1=n-i_1\}$, and $G(I')\subset 
GL_{n-i_1}(A)$. Let $n'=n-i_1$. The induced homomorphism $T_n(A)\to 
T_{n'}(A)$ is naturally split, with kernel $T_{i_1}(A)\subset 
GL_{i_1}(A)\subset GL_n(A)$. 

Let 
\[U_1(I)=\ker \left( U(I)\to U(I')\right) =\ker\left(G(I)\to GL_{i_1}(A)\times GL_{n'}(A)\right).\]
Then $U_1(I)$ is a normal subgroup of $G(I)$, from the last description, 
and 
\[G(I)/U_1(I)\cong U(I')\cdot T_n(A)=T_{i_1}(A)\times G(I').\]

Now $U_1(I)$ may be identified with $M_{i_1,n'}(A)$, the additive group 
of matrices of size $i_1\times n'$ over $A$; this matrix group 
has a natural action of $GL_{i_1}(A)$, and thus of the diagonal matrix 
group $T_{i_1}(A)$, and the resulting semidirect product of $T_{i_1}(A)$ 
with $U_1(A)$ is a subgroup of $G(I)$ (in fact, it is the kernel of 
$G(I)\to  G(I')$). 
This matrix group $M_{i_1,n'}(A)$ is isomorphic, as 
$T_{i_1}(A)$-modules, to the direct sum 
\[\oplus_{i=1}^{i_1}A^{n'}(i),\]
where $A^{n'}(i)$ is the free $A$-module of rank $n'$, with a 
$T_{i_1}(A)$-action given by the the ``$i$-th diagonal entry'' character 
$T_{i_1}(A)\to A^{\times}$. Thus, the semidirect product 
$T_{i_1}(A)U_1(I)$ has a description as a direct product
\[T_{i_1}(A)U_1(i)\cong H\times H\times\cdots \times H=H^{n_i}\]
with $H=A^{n'}\cdot A^{\times}$ equal to the naturally defined semidirect 
product of the free $A$ module $A^{n'}$ with $A^{\times}$, where $A^{\times}$ operates by scalar multiplication.  

Proposition~1.10 and Remark 1.13 in the paper \cite{NS} of Nesterenko and Suslin 
implies immediately that $H\to H/A^{n'}\cong A^{\times}$ induces an 
isomorphism on integral homology. 

We now use the following facts.
\begin{enumerate}
\item[(i)] If $H\subset K\subset G$ are groups, with $H$, $K$ normal in 
$G$, and if $K\to K/H$ induces an isomorphism in integral homology, so 
does $G\to G/H$; this follows at once from a comparison of the two spectral sequences
\[E^2_{r,s}=H_r(G/K,H_s(K,\Z))\Rightarrow H_{r+s}(G,\Z),\]
\[E^2_{r,s}=H_r(G/K,H_s(K/H,\Z))\Rightarrow H_{r+s}(G/H,\Z).\]
\item[(ii)] If $H_i\subset G_i$ are normal subgroups, for 
$i=1,\ldots,n$, such that $G_i\to G_i/H_i$ induce isomorphisms on 
integral homology, then for $G=\prod_{i=1}^nG_i$, $H=\prod_{i=1}^nH_i$, 
the map $G\to G/H$ induces an isomorphism on integral homology. This follows 
from the Kunneth formula.
\end{enumerate}
The fact (ii) implies that $T_{i_1}(A)U_1(I)\to T_{i_1}(A)$ induces an 
isomorphism on integral homology. Then (i) implies that $G(I)\to 
T_{i_1}(A)\times G(I')$ induces an isomorphism on integral homology. By 
induction, we have that $G(I')\to G(I')/U(I')$ induces an isomorphism on 
integral homology. Hence $T_{i_1}(A)\times G(I')\to T_{i_1}\times 
G(I')/U(I')$ also induces an isomorphism on integral homology. Thus, we 
have shown that the composition $G(I)\to G(I)/U(I)=T_n(A)$ induces an 
isomorphism on integral homology. \end{proof}

\section{$\SPL(A^{\infty})^+$ and the groups $L_n(A)$}

We first note that there is a small variation of Quillen's plus construction.

Let $(X,x)$ be a pointed CW complex, $(X_0,x)$ a contractible pointed 
subcomplex, $G$ a group of homeomorphisms of $X$ which acts transitively 
on the path components of $X$, and let $H$ be a perfect subgroup of 
$G$, such that $H$ stabilizes $X_0$.   

Then $X//G$  is clearly path connected, and comes equipped with\\
(i) a natural map $\theta:X//G\to BG=EG/G$, induced by the projection 
$X\times EG\to EG$ \\
(ii) a map $(X_0\times EG)/H\to X//G$, induced by the $H$-stable 
contractible set $X_0\subset X$\\
(iii) a homotopy equivalence $BH\to (X_0\times EG)/H$, such that the 
composition $BH\to X//G\stackrel{\theta}{\to} BG$ is homotopic to the 
natural map $BH\to BG$\\
(iv) a natural map $(X,x)\hookrightarrow  (X//G,x_0)$ determined by the 
base point of $EG$.

Note that, in particular, there is a natural inclusion 
$H\hookrightarrow \pi_1(X//G,x_0)$, which gives a section over 
$H\subset G$ of the surjection $\theta_*:\pi_1(X//G,x_0)\to 
\pi_1(BG,*)=G$.

\begin{lemma}\label{plus}
In the above situation, there is a pointed CW complex $(Y,y)$, together 
with a map $f:(X//G,x_0)\to (Y,y)$ such that\\
(i) the natural composite map 
\[H\hookrightarrow \pi_1(X//G,x_0)\stackrel{f}{\to}\pi_1(Y,y)\]
is trivial\\
(ii) if $g:(X//G,x_0)\to (Z,z)$ such that $H$ is in the kernel of 
\[\pi_1(X//G,x_0)\to \pi_1(Z,z)\]
then $g$ factors through $f$, uniquely upto a pointed homotopy\\
(iii) $f$ induces isomorphisms on integral homology; more 
generally, if $L$ is any local system on $Y$, the map on homology with 
coefficients $H_*(X//G,f^*L)\to H_*(Y,L)$ is an isomorphism\\
(iv) 
$h:(X,x)\to (X',x')$ is a pointed map of such CW complexes with 
$G$-actions, such that $h$ is $G$-equivariant, then there is a map 
$(Y,y)\to (Y',y')$, making  $(X,x)\mapsto (Y,y)$ is functorial (on the 
category of pointed CW complexes with suitable $G$ actions, and 
equivariant maps), and $f$ yields a natural transformation of functors.   
\end{lemma}
The pair $(Y,y)$ is obtained  by applying Quillen's plus construction to $(X//G,x_0)$ 
with respect to the perfect normal subgroup $\tilde{H}$ of 
$\pi_1(X//G,x_0)$ which is generated by $H$. Part (ii) of the lemma is 
in fact the universal property of the plus construction. 
As is well-known, this 
may be done in a functorial way. \emph{We sometimes write 
$(Y,y)=(X//G,x_0)^+$ to denote the above relationship}.

In what follows, the pair $(G,H)$ is invariably $(GL_n(A), A_n)$ for $5\leq n\leq \infty$.  Here 
$A_n$ is the alternating group contained in $N_n(A)$.  The normal subgroups of $GL_n(A)$ 
generated by $A_n$ and $E_n(A)$ coincide with each other. It follows that if we take 
$X=X_0$ to be a point, the $Y$ given by the above lemma is 
just the ``original'' 
$BGL_n(A)^+$.

Recall that there is a natural action of $GL(A)$ on the simplicial 
complex  $\SPL(A^{\infty})$, and hence on its geometric realization 
$|\SPL(A^{\infty})|$.  We apply  
lemma~\ref{plus} with $G=GL(A)$, $H=A_{\infty}$ the 
infinite alternating group,  
$X=|\SPL(A^{\infty})|$, and $X_0=\{x_0\}$ is the vertex of $X$ fixed by $N(A)$ and obtain 
the pointed space 
\[(Y(A),y)=(|\SPL(A^{\infty})|//G), x_0)^+.\]

Taking $X'$ to be a singleton in  (iii) of the above lemma, we get a canonical map 
\[\varphi:(Y(A),y)\to (BGL(A)^+,*)\]
of pointed spaces. 

Let $(\SPL(A^{\infty})^+,z)$ denote the homotopy fibre of $\varphi$.
We define
\[L_n(A)=\pi_n(\SPL(A^{\infty})^+,z)\;\;\forall\;\;n\geq 0.\]
The homotopy sequence of the fibration $\SPL(A^{\infty})^+\to Y(A)\to BGL(A)^+$
 combined with the path-connectedness of $Y(A)$ yields:
\begin{cor}\label{exact}
There is an exact sequence
$$\displaylines{\cdots\to K_{n+1}(A)\to L_n(A)\to \pi_n(Y(A),y)\to 
K_n(A)\cdots\hfill\cr\hfill\cdots\to 
L_1(A)\to \pi_1(Y(A),y)\to K_1(A)\to 
L_0(A)\to 0\cr}$$
where $L_0(A)$ is regarded as a pointed set.
\end{cor}

\begin{lemma}\label{asser(ii)} The natural map 
 $|\SPL(A^{\infty})|\to \SPL(A^{\infty})^+$ 
induces an isomorphism on integral homology.
\end{lemma}
\begin{proof}
We may identify the universal 
covering of $BGL(A)^+$ with $BE(A)^+$, where $BE(A)^+$ is the plus 
construction (see lemma~\ref{plus}) applied to $BE(A)$ with respect to 
the infinite alternating group (or, what is the same thing, with 
respect to $E(A)$ itself). Let $\tilde{\varphi}:\tilde{Y}\to BE(A)^+$ be 
the corresponding pullback map obtained from $\varphi$.
 
We first note that $\SPL(A^{\infty})^+$ is also naturally identified with the homotopy 
fiber of $\tilde{\varphi}$. There is then a homotopy pullback 
$\widehat{\varphi}:\widehat{Y}\to BE(A)$ of $\tilde{\varphi}$ with 
respect 
to $BE(A)\to BE(A)^+$. Thus, our map $\SPL(A^{\infty})\to \SPL(A^{\infty})^+$ may be 
viewed as the natural map on fibers associated to a map
\begin{equation}\label{eqleray}\SPL(A^{\infty})//E(A)\to 
\widehat{Y}
\end{equation}
of Serre fibrations over $BE(A)$. 

From a Leray-Serre spectral sequence argument, we see that since (from 
lemma~\ref{plus}) 
$BE(A)\to BE(A)^+$ induces a isomorphism on integral homology, so does 
$\widehat{Y}\to \tilde{Y}$. Since also $\SPL(A^{\infty})//E(A)\to 
\tilde{Y}$ is a homology isomorphism (from lemma~\ref{plus} again), we 
see that $\SPL(A^{\infty})//E(A)\to \widehat{Y}$ induces an isomorphism 
on integral homology. 

Now we use that the map (\ref{eqleray}) is a map between two 
total spaces of Serre fibrations over a common base, inducing a homology 
isomorphism on these total spaces. We also know that the monodromy 
representation of $\pi_1(BE(A))=E(A)$ on the homology of the fibers is 
trivial, in both cases: for $\widehat{Y}$ this is because it is a 
pullback from a Serre fibration over a simply connected base, while for 
$\SPL(A^{\infty})$, this is one of the key properties we have already 
established (see the finishing sentence of section 2). 
The proof is now complete
modulo the remark below, which is a straightforward consequence of the  
Leray-Serre  spectral sequence of a fibration.
\end{proof}
\begin{remark} Let $p:E\to B$ and $p':E'\to B$ be fibrations  
with fibers $F$ and $F'$ respectively over the base-point $b\in B$.
Let $v:E\to E'$ be a map so that $p'\circ v=p$.
Assume that $B$ is path-connected. Then
$E\to E'$ is a homology isomorphism implies $F\to F'$ is a 
homology isomorphism under the following additional assumption:

$M\neq 0$ implies $H_0(\pi_1(B,b), M)\neq 0$
for every 
$\pi_1(B,b)$-subquotient $M$ of $H_i(F),H_j(F')$ for all $i,j$.
\end{remark}

\section{The H-space structure}

Recall that $BGL(A)^+$ has an H-space structure in a standard way, 
obtained from the direct sum operation on free modules of finte rank; 
this was constructed in \cite{Loday}.

The aim of this section is to prove 
the proposition below.

\begin{proposition}\label{Hspace} The space $Y(A)$ has an H-space 
structure, such that $Y(A)\to BGL(A)^+$ is homotopic to an H-map, for 
the standard H-space structure on $BGL(A)^+$.  
\end{proposition}

\begin{sloppypar}
We first remark that if $V$ is a free $A$-module of finite 
rank, then $|\SPL(V)|//GL(V)$ is homeomorphic to the classifying space 
of the following category $\sSPL(V)$: its objects are simplices in 
$\SPL(V)$ (thus, certain finite nonempty subsets of $SPL(V)$), and 
morphisms $\sigma\to \tau$ are defined to be elements $g\in GL(V)$ such 
that $g(\sigma)\subset \tau$, that is, such that $g(\sigma)$ is a face 
of the simplex $\tau$ of $\SPL(V)$. 

Let $\sAut(V)$ be the category with a single object $*$, with morphisms 
given by elements of $GL(V)$, so that the classifying space $B\sAut(V)$ 
is the standard model for $BGL(V)$. There is a functor $F_V:\sSPL(V)\to 
\sAut(V)$, mapping every obeject $\sigma$ to $*$, and mappng an arrow 
$\sigma\to\tau$ in $\sSPL(V)$ to the corresponding element $g\in 
GLV()$. The fiber $F_V^{-1}(*)$ is the poset of simplices of $\SPL(V)$, 
whose classifying space is thus homeomorphic to $|\SPL(V)|$.

It is fairly straightforward to verify that $B\sSPL(V)$ is homeomorphic 
to  $|\SPL(V)|//GL(V)$ (where we have used the classifying space 
of the translation category of $GL(V)$ as the model for the contractible 
space $E(GL(V))$). One way to think of this is to consider the category 
$\widetilde{\sSPL}(V)$, whose objects are pairs $(\sigma, h)$ with 
$\sigma$ a simplex of $\SPL(V)$, and $h\in GL(V)$, with a unique 
morphism  $(\sigma,h)\to (\tau,g)$  precisely when 
$g^{-1}h(\sigma)\subset \tau$. It is clear that by considering the full 
subcategories of objects of the form $(\sigma,g)$, where $g\in GL(V)$ 
is a fixed element, each of which is naturally equivalent to the poset 
of simplices in $\SPL(V)$, that the classifying space of 
$\widetilde{\sSPL}(V)$  is homeomorphic to $|\SPL(V)|\times E(GL(V)))$. 
Now it is a simple matter to see (e.g., use the criterion of 
Quillen, given in \cite{Srinivas}, lemma~6.1, page~89) that 
$B\widetilde{\sSPL}(V)\to B\sSPL(V)$, given by $(\sigma,h)\mapsto 
h^{-1}(\sigma)$, is a covering space which is a principal 
$GL(V)$-bundle, where the deck transformations are given by the natural 
action of $GL(V)$ on $|\SPL(V)|\times E(GL(V))$. 
\end{sloppypar}

$\sL(V)$ denotes the collection of $A$-submodules $L\subset V$ so that $L$ is free
of rank one and $V/L$ is a free module. 
Now we note that if $V'$, $V''$ are free $A$-modules of finite rank, we note that 
there is a natural inclusion  $\sL(V')\sqcup\sL(V'')\hookrightarrow \sL(V'\oplus V'')$. 
This in turn yields
 a natural map 
\[\varphi_{V',V''}:SPL(V')\times SPL(V'')\to SPL(V'\oplus V''),\]
given by $\varphi_{V',V''}(s,t)=s\sqcup t$.

It follows easily from the definition of $\SPL$ that  the above map on vertices 
induces a simplicial map 
\[\Phi_{V',V''}:\SPL(V')\times \SPL(V'')\to \SPL(V'\oplus V'').\] 
As explained in section 0,  at the level of geometric realisations, this has two
 descriptions. The first description may be used to show that the counterpart of
 lemma~\ref{contiguity}(C) is valid for $\SPL$, namely the homotopy class of the inclusion
 \begin{center}
 $|\SPL(V')|\times |\SPL(V'')|\to |\SPL(V'\oplus V'')|$
 \end{center}
remains unaffected by composition with the action of $g\in \elem(V'\hookrightarrow
 V'\oplus V'')$ on $|\SPL(V'\oplus V'')|$.
\begin{sloppypar}
The second description however is more useful in this context. Let us  abbreviate notation and 
denote the
 (partially ordered) set of simplices of $\SPL(V)$ simply by $\sS(V)$. 
 The desired map  $\sS(V')\times\sS(V'')\to \sS(V'\oplus V'')$ is given simply
  by $(\sigma, \tau)\mapsto \varphi_{V',V''}(\sigma\times\tau)$.   The resulting map 
  $B(\sS(V')\times \sS(V''))\to B\sS(V'\oplus V'')$ is the second description of

 \end{sloppypar}
 \begin{center}
 $|\SPL(V')|\times |\SPL(V'')|\to |\SPL(V'\oplus V'')|$
 \end{center}
for (a) $B\sC'\times  B\sC''
 \cong B(\sC'\times \sC'')$ and (b) $B\sS(V)$ is simply the barycentric subdivision of
  $|\SPL(V)|$. 

\begin{sloppypar}  
  This latter description also allows us to go a step further and define the functor
$\sSPL(V')\times \sSPL(V'')\to 
\sSPL(V'\oplus V'')$, given on objects by 
$(\sigma, \tau)\mapsto \varphi_{V',V''}(\sigma\times\tau)$ as before; on morphisms, it is
given 
 by the natural map $\GL(V')\times \GL(V'')\to GL(V'\oplus 
V'')$. 
Hence on 
classifying spaces, it induces a product
\end{sloppypar}
\[|\SPL(V')|//GL(V')\times |\SPL(V'')|//GL(V'')\to |\SPL(V'\oplus 
V'')|//GL(V'\oplus V'')|.\]
This is clearly compatible with the product
\[BGL(V')\times BGL(V'')\to BGL(V'\oplus V'')\]
under the natural maps induced by the functors $\sSPL\to \sAut$ for the 
three free modules.

One verifies that ${\bf SPL}(A)=\coprod_V \sSPL(V)$, with respect to the 
bifunctor
 \[+:{\bf SPL}(A)\times {\bf SPL}(A)\to {\bf SPL}(A)\]
induced by direct sums on free modules, and the  functors 
$\Phi_{V',V''}$, form a symmetric monoidal category.

 An equivalent 
category, also denoted ${\bf SPL}(A)$ by abuse of notation, is that 
whose objects are pairs $(V,\sigma)$, where $V$ is a free $A$-module 
of finite rank, and $\sigma\in \SPL(V)$ a simplex, and where morphisms 
$(V,\sigma)\to (W,\tau)$ are isomorphisms $f:V\to W$ of $A$-modules such 
that $f(\sigma)$ is a face of $\tau$.  

For the purposes of stabilization, we slightly modify the above to 
consider the related maps 
\[\varphi_{m,n}:SPL(A^m)\times SPL(A^n)\to SPL(A^{\infty})\] 
given by mapping the basis vector $e_i\in A^m$ in the first factor to 
the basis vector $e_{2i-1}\in A^{\infty}$, for each $1\leq i\leq m$, 
and the basis vector $e_j\in A^n$ in the second factor 
to the basis vector $e_{2j}\in A^{\infty}$.  A pair of splittings of 
$A^m$, $A^n$ determine one for the free module spanned by the images of 
the two sets of basis vectors; now one extends this to a splitting of 
$A^{\infty}$ by adjoining the remaining basis vectors of $A^{\infty}$ 
(that is, adjoining those vectors not in the span of the earlier 
images).  If our first two splittings are those given by the basis 
vectors, which correspond to the base points in $|\SPL(A^m)|$ and 
$|\SPL(A^n)|$, the resulting point in $SPL(A^{\infty})$ is again the 
base point of $|\SPL(A)|$.

The corresponding functors 
\[\Phi_{m,n}:\sSPL(A^m)\times \sSPL(A^n)\to \sSPL(A^{\infty})\]
are compatible with similar functors
\[\sAut(A^m)\times \sAut(A^n)\to \sAut(A^{\infty})\]
which, on classifying spaces, yield the diagram of product maps, 
preserving base points, 
\[\begin{array}{ccc}
|\SPL(A^m)//GL_n(A)|\times |\SPL(A^n)|//GL_n(A) & \to & 
|SPL(A^{\infty}|//GL(A)\\
\downarrow&&\downarrow\\
BGL_m(A)\times BGL_n(A)& \to &BGL(A)
\end{array}
\]
where the bottom arrow is the one used in \cite{Loday} to define the 
H-space structure on $BGL(A)^+$. 

As we increase $m$, $n$, the corresponding diagrams are compatible with 
respect to the obvious stabilization maps $|\SPL(A^m)|\hookrightarrow 
|\SPL(A^{m+1})|$, $|\SPL(A^n)\hookrightarrow |\SPL(A^{n+1})|$. Hence we 
obtain on the direct limits a diagram
\[\begin{array}{ccc}
|\SPL(A^{\infty})//GL(A)|\times |\SPL(A^{\infty})|//GL(A) & \to & 
|SPL(A^{\infty}|//GL(A)\\
\downarrow&&\downarrow\\
BGL(A)\times BGL(A)& \to &BGL(A)
\end{array}
\]

From lemma~\ref{plus}, it follows that there is an induced diagram at 
the level of plus constructions
\[\begin{array}{ccc}
Y(A)\times Y(A) & \to & 
Y(A)\\
\downarrow&&\downarrow\\
BGL(A)^+\times BGL(A)^+& \to & BGL(A)^+
\end{array}
\]

It is shown in \cite{Loday} that the bottom arrow defines an H-space 
structure on $BGL(A)$. We claim that, by analogous arguments, the top 
arrow also defines an H-space structure on $Y(A)$. Granting this, the 
map $Y(A)\to BGL(A)^+$ is then an H-map between path connected H-spaces, 
and so the homotopy fiber $Z(A)$ has the homotopy type of an H-space as 
well (and this was what we set out to prove here). 

To show that the product $Y(A)\times Y(A)\to Y(A)$ defines an H-space 
structure, we need to show that left or right translation on $Y(A)$ 
(with respect to this product) by the base point is homotopic to the 
identity. This is also the main point in \cite{Loday}, for the case of 
$BGL(A)^+$.  We first show:
\begin{lemma}\label{loday1}
An arbitrary inclusion $j:\{1,2,\ldots,n\}\hookrightarrow {\mathbb N}$ 
determines an inclusion of $A$-modules $A^n\to A^{\infty}$, given on  
basis vectors by $e_i\mapsto e_{j(i)}$, which induces a map
\[|\SPL(A^n)|//GL_n(A)\to Y(A)\]
which is homotopic (preserving the base point) to the map induced 
by standard inclusion $i_n:A^n\to A^{\infty}$. 
\end{lemma}
\begin{proof}   
We can find an automorphism $g$ of $A^{\infty}$ contained in the 
infinite alternating group $A_{\infty}$, such that $g\circ j=i_n$, 
where $g$ acts on $A^{\infty}$ by permuting the basis vectors (note 
that the induced self-map of $|SPL(A^{\infty})|\times EGL(A)$ fixes the 
base point). Regarding $g$ as an element of 
$\pi_1(|\SPL(A^{\infty})|//GL(A))$, this implies that the maps $(i_n)_*$ 
and $j_*$, considered as elements of the set of pointed homotopy classes of maps
 
\[   \left[ |\SPL(A^n)|//GL_n(A),|\SPL(A^{\infty})|//GL(A) \right],   \] 
are related by $g_*(j_*)=(i_n)_*$, where $g_*$ denotes the action of the 
fundamental  group of the target on the set of pointed homotopy classes 
of maps. However,  $g$ is in the kernel of the map on fundamental groups  
associated to the map 
\[|\SPL(A^{\infty})|//GL(A)\to   (|\SPL(A^{\infty}|//GL(A))^+.\]
Hence the induced maps  
\[|\SPL(A^n)|//GL(A)\to (|\SPL(A^{\infty})|//GL(A))^+\]
determined by $i_n$ and $j$ are homotopic. 
\end{proof}
\begin{cor}\label{loday2}  The map $Y(A)\to Y(A)$ defined by an arbitrary 
injective map $\alpha:{\mathbb N}\hookrightarrow{\mathbb N}$ is 
homotopic, preserving the base point, to the identity.
\end{cor}
\begin{proof}

We first note that if for $n\geq 5$, we let 
$Y_n(A)=(|\SPL(A^n)|//GL_n(A))^+$ be the result of applying 
lemma~\ref{plus} to $|\SPL(A^n)|//GL_n(A)$ for the alternating group 
$A_n$, then there are natural maps $Y_n(A)\to Y(A)$, preserving base 
points, and inducing an isomorphism 
$\limdir{n}\, \pi_*(Y_n(A))=\pi_*(Y(A))$.

We claim that if $\alpha_n:\{1,2\ldots,n\}\hookrightarrow {\mathbb N}$ 
is the inclusion induced by restricting $\alpha$, then the induced map 
$(\alpha_n)_*:Y_n(A)\to Y(A)$ is homotopic, preserving the base points, 
to the natural map $Y_n(A)\to Y(A)$. This follows from 
lemma~\ref{loday1}, combined with the defining universal property of the 
plus construction, given in lemma~\ref{plus}.

This implies that the map $\alpha:Y(A)\to Y(A)$ must then induce 
isomorphisms on homotopy groups, and hence is a homotopy equivalence, 
by Whitehead's theorem.

Thus, we have a map from the set of such injective maps $\alpha$ to the 
group of base-point preserving homotopy classes of self-maps of $Y(A)$. 
This is in fact a homomorphism of monoids, where the operation on the 
injective self-maps of ${\mathbb N}$ is given by composition of maps.
  
Now we use a trick from \cite{Loday}: any homomorphism of monoids from 
the monoid  of injective self-maps of ${\mathbb N}$ to a group is a 
trivial homomorphism, mapping all elements of the monoid to the 
identity. This is left to the reader to verify (or see \cite{Loday}).
\end{proof}

We note that the above monoidal category ${\bf SPL}(A)$ can be used to 
give another, perhaps more insightful construction of the homotopy type 
$Y(A)$, analogous to Quillen's ${\mathcal S}^{-1}{\mathcal S}$ 
construction for $BGL(A)^+$. We sketch the argument below.

We first take ${\bf SPL}_0(A)$ to be the full subcategory of ${\bf 
SPL}(A)$ consisting of pairs $(V,\sigma)$ where $\sigma\in SPL(V)$, 
i.e.,$\sigma$ is a 0-simplex in $\SPL(V)$. This full subcategory is in 
fact a monoidal subcategory, which is a groupoid (all arrow are 
isomorphisms). Also, ${\bf SPL}(A)$ is a symmetric monoidal category, 
in that the sum operation is commutative upto coherent natural 
isomorphisms. Then, using Quillen's results (see Chapter~7 in 
\cite{Srinivas}, particularly Theorem~7.2), one can see that ${\bf 
SPL}_0(A)^{-1}{\bf SPL}(A)$ is a monoidal category whose classifying 
space is a connected H-space, which is naturally homology equivalent to 
$|\SPL(A^{\infty})|//GL(A)$. This then forces this classifying space to 
be homotopy equivalent to $Y(A)$, such that the H-space operations are 
compatible upto homotopy. This is analogous to the identification made 
in Theorem~7.4 in \cite{Srinivas} of ${\mathcal S}^{-1}{\mathcal S}$ 
with $K_0(R)\times BGL(R)^+$ for a ring $R$, and appropriate ${\mathcal 
S}$. (We do not get the factor $K_0$ appearing in our situation since we 
work only with free modules).


\section{Theorem 1 and the groups $\sH_n(A^{\times})$}
\noindent \emph{Proof of Theorem 1}. In view of Proposition~\ref{Hspace},
we see that $\SPL(A^{\infty})^+$, the homotopy fiber of
the H-map $Y(A)\to BGL(A)^+$, 
 is a H-space as well. It follows
 that $L_0(A)=\pi_0(\SPL(A^{\infty})^+)$ is a monoid. Furthermore,
the arrow
 $K_1(A)\to L_0(A)$ in Corollary~\ref{exact} is a monoid
homomorphism. Thus this corollary produces an exact
sequence of Abelian groups.

$\SPL(A^n)$ has a canonical base point fixed under the action of
$N_n(A)$. As in sections 3 and 4, this gives a
 natural 
inclusion $BN_n(A)\to\SPL_n(A)//GL_n(A)$. This is a homology isomorphism
 by lemma~\ref{suslin-lemma2}. Taking direct limits over all $n\in\N$,
 we see that $BN(A)\to \SPL(A^{\infty})//GL(A)$ is a homology isomorphism.

Applying Quillen's plus construction with respect to
the normal subgroup of $N(A)$ generated by the infinite alternating group,
 we obtain a space $BN(A)^+$. That $BN(A)^+$ has a canonical 
H-space structure follows easily by the method of the previous section.
Now the map $BN(A)^+\to Y(A)$ obtained by lemma~\ref{plus}(ii)
 is a homology isomorphism of simple path-connected CW complexes
and is therefore a homotopy equivalence 
(see \cite {Hatcher}) Theorem 4.37, page 371 and Theorem 4.5, page 346).
This gives the isomorphism $\sH_n(A^{\times})\to L_n(A)$. The theorem
now follows from corollary~\ref{exact}.

We now turn to the description of the groups. 
Let $X=B(A^{\times})$. Let $X_+=X\sqcup\{*\}$ be the pointed space 
with $*$ as its base-point. Let $QX_+$ be the direct limit of
$\Omega^n\Sigma^nX_+$ where $\Sigma$ denotes reduced suspension.  
\begin{proposition}\label{stable}

$\sH_n(A^{\times})\cong \pi_n(QX_+) $. 
\end{proposition}
This statement was suggested to us by Proposition 3.6 of \cite{Segal}. 

A complete proof of the proposition was shown us by Peter May.
A condensed version of what we learnt from him is
given below.

Theorem 2.2, page 67 of \cite{May1} asserts that 
$\alpha_{\infty}:C_{\infty}X_+\to QX_+$ is a group completion. This is
proved in pages 50-59, \cite{May3}. The $C_{\infty}$ here is a 
particular case of the construction 2.4, page 13 of \cite{May2}, given
 for any operad. For $C_{\infty}(Y)$, where $Y$ is a pointed space,
the easiest definition to work with is found in 
May's review of \cite{Segal1}. 
It runs as follows. Let $V=\cup_{n=0}^{\infty}\R^n$. Let $C_k(Y)$ be the 
collection of ordered pairs $(c,f)$ where $c\subset V$ has cardinality $k$
 and $f:c\to Y$ is any function. We identify $(c,f)$ with $(c',f')$ if
\\(i) $c'\subset c$, (ii) $f|c'=f'$, and (iii) $f(a)=*$ for all $a\in c,
a\notin c'$. Here $*$ stands for the base-point of $Y$. Then $C_{\infty}Y$
 is the space obtained from the disjoint union of the $C_k(Y), k\geq 0$ by
performing these identifications. This is a H-space. 
In our case, when $Y=X\sqcup\{*\}$, 
it is clear that  $C_{\infty}Y$ is the disjoint union of all 
the $C_k(X)$ 
as a topological space. Thanks to ``infinite codimension'' one gets 
easily the homotopy equivalence of $C_k(X)$ with $X^k//S_k$ where $S_k$ is
the permutation group of $\{1,2,..,k\}$. Now assume that $X$ is 
any path-connected space
equipped with a nondegenerate base-point $x\in X$. This $x$ gives
 an inclusion of $X^n\hookrightarrow X^{n+1}$. Denote by $X^{\infty}$
the direct limit of the $X^n$. Thus $X^{\infty}$ is a pointed space equipped
 with the action of the infinite permutation group $S_{\infty}=\cup n S_n$.
Put $Z= X^{\infty}//S_{\infty}$. As in section 4, we obtain $Z^+$ by
the use of the infinite alternating group. As in section 5, we see that
this is a H-space. It is an easy matter to check that the group
completion of $\sqcup_kC_kX$ is homotopy equivalent to $\Z\times Z^+$.
This shows that $\pi_n(QX_+)\cong \pi_n(Z^+)$ for all $n>0$. 

The proposition is the particular case: $X=B(A^{\times})$.

\section{Polyhedral structure of the enriched Tits building}

From what has been shown so far, we see that it is of interest to determine 
the stable rational homology of the flag complexes $|\FL(A^n)|$ (or equivalently, 
of $|\SPL(A^n)|$, or $\ET(A^n)$). We will construct a spectral sequence that, 
in principle, gives an inductive procedure to do so.

But first we introduce some notation and a definition for posets.
\begin{sloppypar}
Let $P$ be a poset. For $p\in P$, we put $e(p)=BL(p)$ where $L(p)=\{q\in P|q\leq p\}$
 and $\partial e(p)=BL'(p)$ where  $L'(p)=L(p)\setminus\{p\}$. 
 If $\partial e(p)$ is homeomorphic to a sphere for every $p\in P$,  we say the poset $P$ 
  is {\em polyhedral}. We denote by $d(p)$ the dimension of $e(p)$. When $P$ is 
  polyhedral, the space $BP$ gets the structure of a CW complex with 
  $\{e(p):p\in P\}$ as the closed cells. Its $r$-skeleton is 
  $BP_r$ where $P_r=\{p\in P:d(p)\leq r\}$. The homology of $BP$ is then computed 
by the associated complex of cellular chains $Cell_{\bullet}(BP)$, where 
\end{sloppypar}

\[Cell_r(BP)=\bigoplus_{\{p|\dim e(p)=r\}}H_r(e(p),\partial e(p),\Z).\]

\begin{lemma}\label{poly} $\sE(A^n)$ is a polyhedral poset in the above 
sense. Its dimension is $n-1$. 
\end{lemma}
\begin{proof}
First consider the case when $p\in SPL(A^n)$ is a maximal element in $\sE(A^n)$. Then 
$p$ is an unordered collection of $n$ lines in $A^n$ (here, as in \S2, a ``line'' denotes 
a free $A$-submodule of rank 1 which is a direct summand, and the set of lines in $A^n$ is 
denoted by $\sL(A^n)$). Note that the subset $p\subset \sL(A^n)$ of cardinality $n$ determines 
a poset $\tilde{p}$, whose elements are chains $q_{\bullet}=\{q_1\subset 
q_2\subset \cdots \subset q_r=p\}$ of nonempty 
subsets, where $r_{\bullet}\leq q_{\bullet}$  if each $q_i$ is an $r_j$ 
for some $j$, i.e., the 
``filtration'' $r_{\bullet}$ ``refines'' $q_{\bullet}$. We claim that, 
from the definition of the partial order on $\sE(A^n)$, the poset 
$\tilde{p}$ is naturally isomorphic to the poset  
$L(p)$. Indeed, an element  $q\in \sE(A^n)$ consists of a pair, consisting of a partial flag 
\[0=W_0\subset W_1\subset \cdots\subset W_r=A^n\]
such that $W_i/W_{i-1}$ is free, and a sequence $t_1,\ldots,t_r$ with 
$t_i\in SPL(W_i/W_{i-1})$. The condition that this element of $\sE(A^n)$ 
lies in $L(p)$ is that each $W_i$ is a direct sum 
of a subset of the lines in $p$, say $q_i\subset p$, giving the chain of subsets $q_1\subset q_2...\subset q_r=p$; 
the splitting $t_i$ is uniquely determined by the lines in $q_i\setminus q_{i-1}$.  

Let $\Delta(p)$ be the $(n-1)$-simplex with $p$ as its set of vertices. 
Now the chains of non-empty subsets of $p$ correspond to simplices in the barycentric subdivison $sd\Delta(p)$, where 
the barycentre $b$ corresponds to the chain $\{p\}$.  Hence, from the definition of $\tilde{p}$, it is clear that 
it is isomorphic to the poset whose elements are simplices in the barycentric subdivision of $\Delta_n$ with $b$ 
as a vertex, with partial order given by reverse inclusion. Hence $B\tilde{p}$ is naturally identified with the 
subcomplex of the second barycentric subdivision $sd^2\Delta(p)$ which is the union of all simplices containing 
the barycentre. This explicit description implies in particular that $BL'(p)$ is homeomorphic to $S^{n-2}$ 
(with a specific triangulation). 

Before proceeding to the general case,
we set up the relevant notation for orientations. 
For a set $q$ of cardinality $r$, we put $\det(q)=\wedge ^r\Z[q]$ where
$\Z[q]$ denotes the free Abelian group with $q$ as basis. 
we observe that there is a natural isomorphism:
\[H_{n-1}(e(p),\partial e(p))\cong 
H_{n-1} (\Delta(p),\partial\Delta(p))=\det(p).\]

Now let $p\in \sE(A^n)$ be arbitrary, corresponding to a partial flag 
\[0=W_0\subset W_1\subset \cdots\subset W_r=A^n.\]
and splittings $t_i\in SPL(W_i/W_{i-1})$.  Then the natural map
\[\prod_{i=1}^r \sE(W_i/W_{i-1})\to \sE(A^n)\]
is an embedding of posets, where the product has the ordering 
given by $(a_1,\ldots,a_r)\leq (b_1,\ldots,b_r)$ precisely when $a_i\leq 
b_i$ in $\sE(W_i/W_{i-1})$ for each $i$. One sees that, by the 
definition of the partial order in $\sE(A^n)$, the induced map
\[\prod_{i=1}^r L(t_i)\to L(p)\]
is bijective. Hence there is a homeomorphism of pairs
\[(BL(p),BL'(p))=\prod_{i=1}^r(BL(t_i),BL'(t_i)),\]
and so $BL'(p)\cong S^{n-r-1}$, and $BL(p)$ is an $n-r$-cell.
\end{proof}

We now proceed to construct the desired spectral sequence. We use the 
following notation: if $p\in \sE(V)$, where $V$ is a free $A$-module of 
finite rank, and $W_1\subset V$ is the smallest non-zero submodule in 
the partial flag associated to $p$, define $t(p)=\rank W_1-1$. Clearly 
$t:\sE(V)\to \Z$ is monotonic. Hence $F_r\sE(V)=\{p\in\sE(V)|t(p)\leq 
r\}$ is a sub-poset. Define 
\[F_r\ET(V)=BF_r\sE(V)=\cup\{e(p)|t(p)\leq r\},\]
so that
\[F_0\ET(V)\subset F_1\ET(V)\subset \cdots 
F_{n-1}\ET(V)=\ET(V)\]
is an increasing finite filtration of the CW complex $\ET(V)$ by 
subcomplexes. Hence there is an 
associated spectral sequence
\[E^1_{r,s}=H_{r+s}(F_r\ET(V),F_{r-1}\ET(V),\Z)\Rightarrow 
H_{r+s}(\ET(V),\Z).\]
Our objective now is to recognise the above $E^1$ terms. 

It is convenient to use the complexes of cellular chains for these 
sub CW-complexes, which are thus sub-chain complexes of 
$Cell_{\bullet}(\ET(V))$. For simplicity of notation, we write $Cell_{\bullet}(V)$ 
 for  $Cell_{\bullet}(\ET(V))$.
We have the description
\[E^1_{r,s}=H_{r+s}({\rm gr}^F_rCell_{\bullet}(V)).\]

We will now exhibit ${\rm gr}^F_rCell_{\bullet}(V)$ as a direct  sum of 
complexes. Let $W\subset V$ be a submodule such that $W,V/W$ are both 
free, and $\rank W=r+1$. Let $q\in SPL(W)$. The we have an inclusion of 
chain complexes
\[Cell_{\bullet}(e(q))\tensor Cell_{\bullet}(V/W)\subset 
Cell_{\bullet}(W)\tensor Cell_{\bullet}(V/W)\subset F_rCell_{\bullet}(V).\] 
It is clear that
\[{\rm image}\,Cell_{\bullet}(\partial e(q))\tensor 
Cell_{\bullet}(V/W)\subset F_{r-1}Cell_{\bullet}(V),\]
so that we have an induced homomorphism of complexes
\[\left(Cell_{\bullet}(e(q))/Cell_{\bullet}(\partial e(q))\right)\tensor 
Cell_{\bullet}(V/W)\to {\rm gr}^F_rCell_{\bullet}(V).\]
Composing with the natural chain homomorphism
\[H_r(e(q),\partial e(q),\Z)[r]\to 
\left(Cell_{\bullet}(e(q))/Cell_{\bullet}(\partial e(q))\right)\]
for each $q$, we finally obtain a chain map
\[I:\bigoplus_{(W,q\in SPL(W))}
H_r(e(q),\partial e(q),\Z)[r]\tensor Cell_{\bullet}(V/W)\to 
{\rm gr}^F_rCell_{\bullet}(V).\]
Finally, it is fairly straightforward to verify that $I$ is an 
isomorphism of complexes. 

We deduce that the $E^1$ terms have the following description:
\[E^1_{r,s}=\bigoplus_{\scriptstyle{\begin{array}{c}
\rank W=r+1\\
q\in SPL(W)
\end{array}}}
\det(q)\tensor H_s(Cell_{\bullet}(V/W),\Z).\] 
We define $\sL_r(V)$ to be the collection of 
$q\subset \sL(V)$ of cardinality $(r+1)$ for which (a) and (b) below
hold:
\\(a)  $\oplus\{L:L\in q\}\to V$ is injective. Its image will be
denoted by $W(q)$  
\\(b) $V/W(q)$ is free of rank $(n-r-1)$.

Summarising the above, we obtain:
\begin{theorem}\label{ss} There is a spectral sequence with $E^1$ terms
\[E^1_{r,s}=\bigoplus_{\scriptstyle{\begin{array}{c}
q\in\sL_r(V)\\

\end{array}}}
\det(q)\tensor H_s(\ET(V/W(q)),\Z).\] 
that converges to $H_{r+s}(\ET(V))$. We note that $E^1_{r,s}=0$ whenever 
$(r+s)\geq (n-1)$ with one exception: $(r,s)=(n-1,0)$. Here $V\cong A^n$. 
\end{theorem}
\section{Compatible homotopy}
It is true \footnotemark \footnotetext{this only requires the analogue of lemma~\ref{contiguity}(A) 
for the enriched Tits building. More general statements  are contained in lemmas~\ref{silly} and
\ref{comp} .} that $i:\ET(W)\times \ET(V/W)\hookrightarrow\ET(V)$ has
the property that $g\circ i$ is freely homotopic (not preserving base points)
to $i$  whenever $g\in\elem(W\hookrightarrow V)$.

There are several closed subsets of $\ET(A^n)$ with the property that 
homotopy class of  the inclusion
morphism into  $\ET(A^n)$ remains unaffected by composition with the action of
$g\in E_n(A)$.   To prove that the union of a finite collection of such closed subsets
 has the same property, one would require   the homotopies 
 provided for any two members of the collection to agree on their intersection. This is
 the problem we are concerned with in this section.   

We proceed to set up the notation for the problem.

With $q\in \sL_r(V)$ as in theorem~\ref{ss}, we shall define the
subspaces $U(q)\subset \ET(V)$ as follows. Let
$W(q)=\oplus\{L|L\in q\}$. We  regard $q$ as an element of $SPL(W(q))$
and thus obtain the cell $e(q)=BL(q)\subset\ET(W(q))$. This gives 
the inclusion 
\\$\ET'(q)=e(q)\times \ET(V/W(q))\subset \ET(W(q))\times\ET(V/W(q))\subset\ET(V)$.
\\We put $U(q)=\cup\{\ET'(t)|\emptyset\neq t\subset q\}$.
\\ \\
\emph{Main Question:} Let  $i:U(q)\hookrightarrow \ET(V)$ denote the inclusion.  Is it true that
$g\circ i$ is homotopic to $i$ for every $g\in\elem(V,q)$?
\\We focus on the apparently weaker question below.
\\ \emph{Compatible Homotopy Question}: Let $M\subset V$ be a 
submodule complementary to $W(q)$.  Let $g'\in GL((W(q))$ be elementary,
i.e. $g'\in \elem(W(q),q)$. Define $g\in GL(V)$ by $gm=m$ for all $m\in M$ and
 $gw=g'w$ for all $w\in W(q)$.  Is it true that $g\circ i$ is homotopic to $i$? 
   
 Assume that the second question has an affirmative answer
  in all cases. In particular, this holds when $M=0$. Here $V=W(q)$ and $g=g'$ is
  an arbitrary element of $\elem(V,q)$. Let $t$ be a non-empty subset of $q$. 
  Then $U(t)\subset U(q)$.
  We deduce  that $j:U(t)\hookrightarrow \ET(V)$ is homotopic to
 $g\circ j$ for all $g\in\elem(V,q)$.  But $\elem(V,t)=\elem(V,q)$. Thus the Main 
 Question has an affirmative answer for $(q,i)$ replaced by $(t,j)$, which of course,
 up to a change of notation, covers the general case.
\begin{proposition}\label{three} The compatible homotopy question has
an affirmative answer if $q\in\sL_r(V)$ and $r\leq 2$.

\end{proposition} 

The rest of this section is devoted to the proof of this proposition.
To proceed, we will require to introduce the class $\sC$.

This is
our set-up. Let $X$ be a finite set, let $V_x$ be a finitely generated
free module for each $x\in X$ and let $V=\oplus\{V_x:x\in X\}$.  

\begin{sloppypar}
Let $s=\Pi_{x\in X}s(x)\in \Pi_{x\in X}SPL(V_x)$. For each $x\in X$, we regard $s(x)$ as a subset of $\sL(V)$ and 
put $Fs=\cup \{s(x)|x\in X\}$. Thus $Fs\in SPL(V)$.  \emph{The collection of maps  
 $f:\Pi_{x\in X}\ET(V_x)\to \ET(V)$ with the property that  $f(\Pi_{x\in X}BL(s(x)))\subset BL(Fs)$
for all $s\in \Pi_{x\in X}SPL(V_x)$ is denoted by $\sC$}. See lemma~\ref{poset1} and
proposition~\ref{mainflags} and its proof for relevant notation.   
\end{sloppypar}

Every maximal chain $C$ of subsets of $X$ (equivalently every total ordering of $X$)
 gives a member $i(C)\in\sC$. For
instance, if $X=\{1,2,...,n\}$ and $C$ consists of the sets 
$\{1,2,...,k\}$ for $1\leq k\leq n$, we put
\\$D_k=\oplus_{i=1}^kV_i$ and $E=\Pi_{i=1}^n\ET(D_i/D_{i-1})$, denote by
\\$u:\Pi_{x\in X}\ET(V_x)\to E$ and $v:E\to\ET(V)$ the natural
isomorphism and natural inclusion respectively, and put $i(C)=v\circ u$.

\begin{lemma}\label{silly} The  above space $\sC$ is contractible.

\end{lemma}
\begin{proof}
The aim is to realise $\sC$ as the space of
$\Lambda$-compatible maps for a suitable $\Lambda$
 and appeal to Proposition~\ref{covering}.
Let $\Lambda(x)=\{L(S)|\emptyset\neq S\subset SPL(V_x),\emptyset\neq L(S)\}$.  
Proposition~\ref{mainflags} and the proof of lemma~\ref{poset1} combine to
show  that the subspaces 
$\{B\lambda(x):\lambda(x)\in\Lambda(x)\}$ give an 
admissible cover of $\ET(V_x)$.

 For $\lambda=\Pi_{x\in X}\lambda(x)\in  \Lambda=\Pi_{x\in X}\Lambda(x)$,
we put $I(\lambda)=\Pi_{x\in X}B\lambda(x)$ and deduce that
$\{I(\lambda):\lambda\in\Lambda\}$
 gives an admissible cover of
$\Pi_{x\in X}\ET(V_x)$.

We define next a closed $J(\lambda)\subset \ET(V)$ for 
every $\lambda\in\Lambda$ with the properties: 
\\(A): $J(\lambda)\subset J(\mu)$
whenever $\lambda\leq \mu$ and
\\(B):  $J(\lambda)$ is contractible for 
every $\lambda\in\Lambda$.

For each 
$\lambda(x)\in\Lambda(x)$, let $U\lambda(x)$ be its set of upper bounds in $SPL(V_x)$. 
It follows that
that $LU\lambda(x))=\lambda(x)$.  As observed before, we have
\\$F:\Pi_{x\in X} SPL(V_x)\to SPL(V)$.  Thus given $\lambda=\Pi_x\lambda(x)\in \Lambda$,
\\we set $H(\lambda)=F(\Pi_{x\in X}U\lambda(x))\subset  SPL(V)$ and then put 
$J(\lambda)=BLH(\lambda)\subset \ET(V)$. 

The space of $\Lambda$-compatible maps $\Pi_{x\in X}\ET(V_x)\to \ET(V)$ is seen 
to coincide with $\sC$. That the $J(\lambda)$ satisfy property (A)
stated above is straightforward. 

The contractibility of 
$J(\lambda)$ for all $\lambda\in\Lambda$ is guaranteed   
by proposition~\ref{mainflags} once it is checked that these
sets are nonempty. But we have already noted that $\sC$ is nonempty.
Let $f\in\sC$. Now $I(\lambda)\neq \emptyset$ and $f(I(\lambda))\subset 
J(\lambda)$ implies $J(\lambda)\neq \emptyset$. Thus the $J(\lambda)$
 are contractible, and 
as said earlier, an application
of Proposition~\ref{covering} completes the proof of the lemma.

\end{proof}
   We remark that  the class $\sC$ of maps $\Pi_{i=1}^n\ET(W_i)\to\ET(\oplus_{i=1}^n W_i)$ 
 has been defined in general.

 We will continue to employ the notation: $V=\oplus\{V_x:x\in X\}$  all through this section.    
Let  $P$ be a partition of $X$. Each $p\in P$ is a subset of $X$  and we put
\begin{center}$V_p=\oplus \{V_x|x\in p\}$ and $\ET(P)=\Pi\{\ET(V_p)|p\in P\}$.
\end{center}
 When $Q\leq P$ is a partition of $X$ (i.e. $Q$ is finer than $P$), we shall define 
 the contractible collection $\sC(Q,P)$ of maps 
  $f:\ET(Q)\to\ET(P)$
 by demanding  (a) that $f$ is the product of maps $f(p)$ 
 \[f(p):\Pi\{\ET(V_q):q\subset p\,\,\mathrm{and}\,\,q\in Q\}\to \ET(V_p)\]
 \begin{sloppypar}
 and  also (b) each $f(p)$ is in the class $\sC$. For this one should note that
  $V_p=\oplus\{ V_q:q\in Q\,\,\mathrm{and}\,\, q\subset p\}$.
\end{sloppypar} 
 We observe next that there is a distinguished  collection $\sD(Q,P)\subset \sC(Q,P)$. To
 see this, recall that we had the embedding $i(C)$ for every maximal chain $C$ of subsets of $X$ 
 (alternatively, for every total ordering of $X$). Given $Q\leq P$, denote the set 
 of total orderings of 
 $\{q\in Q:q\subset p\}$ by $T(p)$, for every $p\in P$.  The earlier $C\mapsto i(C)$
 now yields, after taking a product over $p\in P$,
  \\$i:\Pi\{T(p):p\in P\}\to \sC(Q,P)$,  and we denote by $\sD(Q,P)\subset \sC(Q,P)$ the
  image of $i$.
  
 The lemma below is immediate from the definitions.  
  \begin{lemma}\label{comp}
 Given partitions $R\leq Q\leq P$ of $X$, if $f$ is in $\sC(R,Q)$ (resp. in $\sD(R,Q)$)and 
 $g$ is in $\sC(Q,P)$ (resp. in $\sD(Q,P)$), then 
 it follows that $g\circ f$ is in $\sC(R,P)$ (resp. $\sD(R,P)$).
\end{lemma}
\bigskip 
We will soon have to focus on the fixed points of certain unipotent $g\in GL(V)$ on
 $\ET(V)$. For instance, if $x,y\in X$ and $x\neq y$, we may consider $g=id_V+h$ where 
$h(V)\subset V_y$ and $h(V_z)=0$ for all $z\in X, z\neq x$.   Let $C$ be a chain
of subsets of $X$, so that $X\in C$. This chain $C$ 
gives rise to a partition $P(C)$ of $X$  and also 
$i(C)\in\sD(P(C),\{X\})$ in a natural manner. Let $C_x=\cap \{S\in C:x\in S\}$. Then $C_x\in C$
because $C$ is a chain. Define $C_y$ in a similar manner.  
\emph{We say the chain $C$ is
 $(x,y)$-compatible if  $C_y\subset C_x$ and $C_x\neq C_y$ }. This condition on $C$
 ensures that the embedding $i(C):\ET(P(C))\to \ET(V)$ has its image within  the fixed points of
 the above $g\in GL(V)$. 
 
 Now let $Q$ be a partition of $X$ so that $q\in Q,x\in q$ implies $y\notin q$.  We shall define
 next the class of $(x,y)$-compatible $\sC$ maps $\ET(Q)\to\ET(V)$ in the
 following manner. Let  $\Lambda$ be the set of chains $C$ of subsets of $X$  so that
 $X\in C$ and $Q\leq P(C)$ (i.e. $Q$ is finer than the partition $P(C)$).  For each
 $C\in\Lambda$,  let $Z(C)$   be the collection of 
 $i(C)\circ f$ where  
 $f\in \sC(Q,P(C))$. Finally, let $Z=\cup\{Z(C):C\in\Lambda\}$. This set $Z$ is defined
 to be the collection of $(x,y)$-compatible maps of class $\sC$ from  $\ET(Q)$ to  $\ET(V)$.  
 Every $z\in Z$ is a map
  $z:\ET(Q)\to\ET(V)$  whose image is contained in the fixed points of 
the above $g$ on $\ET(V)$. Furthermore, in view of lemma~\ref{comp}, this collection of maps is 
contained in $\sC(Q,\{X\})$. 
 \begin{lemma}\label{xy} Let $Q$ be a partition of $X$ that separates $x$ and $y$. 
 Then the collection of $(x,y)$-compatible class  $\sC$ maps $\ET(Q)\to\ET(V)$ is contractible.
 
 \end{lemma}
 \begin{proof} 
In view of the fact that each $\sC(Q,P)$ is contractible, by cor~\ref{cor-covering2}, it follows 
 that the space  of $(x,y)$-compatible chains is homotopy equivalent to
 $B\Lambda$, where  $\Lambda$ is the poset of chains $C$ in the previous
 paragraph.  It remains to show that $B\Lambda$ is
 contractible.
 
 We first consider the  case where $Q$ is the set of all singletons of $X$.   
Let $\sS$ be the collection of subsets $S\subset X$ so that $y\in S$ and $x\notin S$. 
For $S\in\sS$, let $\sF(S)$ be the collection of chains $C$ of subsets of $X$ so that
$S\in C$ and $X\in C$. We see that $\Lambda$ is precisely the union of $\sF(S)$ 
 taken over all $S\in\sS$.   Let $D$ be a finite subset of $\sS$. We see that the
 intersection of the $B\sF(S)$, taken over $S\in D$, is nonempty if and only
 if $D$ is a chain. Furthermore, when $D$ is a chain, this intersection is clearly
 a cone, and therefore contractible. By cor~\ref{cor-covering2}, we see that  $B\Lambda$ 
 has the same homotopy type as the classifying space of the poset of chains
 of $\sS$.   But this is simply the barycentric subdivision of $B\sS$. But the latter is
 a cone as well, with $\{y\}$ as vertex. This completes the proof that  $B\Lambda$ is
 contractible, when $Q$ is the finest possible partition of $X$.
 
 We now come to the general case, when $Q$ is an arbitrary 
partition of $X$ that separates $x,y$.
 So we have $x',y'\in Q$ with $x\in x',y\in y'$ and $x'\neq y'$.  The set $\Lambda$ 
 is identified with the collection of chains $C'$ of subsets of $Q$ so that 
 \\(a) $Q\in C'$, and (b) there is some $L\in C'$ so that $x'\notin L$ and $y'\in L$.  
 
 Thus the general case follows from the case considered first: one replaces $(X,x,y)$
 by $(Q,x',y')$.
  \end{proof}
In a similar manner, we may define, for every ordered $r$-tuple   $(x_1,x_2,...,x_r)$ of  distinct
elements of $X$,  the set of $ (x_1,x_2,...,x_r)     $-compatible chains $C$--we
demand that for each $0<i<r$, there is a member $S$ of the chain so that $x_i\notin S$ and
$x_{i+1}\in S$.  Let $Q$ be a partition of $X$ that separates $x_1,x_2,...,x_r$.  Then the poset of
chains $C$ , compatible with respect to this ordered $r$-tuple, and for which $Q\leq i(P)$,
 is also contractible. One may see this through an inductive version of the proof of the above lemma.
  A corollary
 is that the collection of $(x_1,...,x_r)$-compatible class $\sC$ maps $\ET(Q)\to\ET(V)$ 
 is also contractible.
 We skip the proof. This result is employed in the proof of Proposition~\ref
 {three}
 for $r=2$ (which has already been verified in the above lemma), and for $r=3$, with $\#(Q)\leq 4$. 
  Here it is a simple verification that the poset of chains that arises as above has its 
 classifying space homeomorphic
to a point or a closed interval. 
\\ \\
We are now ready to address the proposition. For this purpose,
we assume that there is $c\in X$ so
that $V_x\cong A$ for all $x\in X\setminus\{c\}$. To obtain  
consistency with the notation of the proposition, we set  $q=X\setminus\{c\}$. 
The closed subset $U(q)\subset \ET(V)$ in the proposition is the union of
$\ET'(t)$ taken over all $\emptyset\neq t\subset q$. For such $t$, we
have $W(t)=V(t)=\oplus\{V_x:x\in t\}$. Recall that $\ET'(t)$ is
the product of the cell $e(t)\subset \ET(V(t))$ with $\ET(V/V(t))$.
To proceed, it will be necessary to 
give a contractible class of maps $D\to\ET(V)$ for certain closed 
subsets $D\subset U(q)$.

The closed subsets $D\subset U(q)$ we consider have the following shape.
For each $\emptyset \neq t\subset q$, we first select a closed subset 
$D(t)\subset e(t)$ and then take $D$ to be the union of the
$D(t)\times \ET(V/V(t))$, taken over all such $t$. This $D$ remains
unaffected if $D(t)$ is replaced by its saturation $sD(t)$. Here
$sD(t)$ is the collection of $a\in e(t)$ for which
$\{a\}\times \ET(V/V(t))$ is contained in $D$.

When $\emptyset \neq t\subset q$, we denote by 
 $p(t)$ the partition of $X$ consisting of all the singletons
contained in $t$, and in addition, the complement $X\setminus t$.
Then there is a canonical identification $j(t):\ET(p(t))\to \ET(V/V(t))$.

A map $f:D\to\ET(V)$ is said to be in class $\sC$ if
for every $\emptyset\neq t\subset q$ and for every $a\in sD(t)$, the
map
$\ET(p(t))\to\ET(V)$ given by $b\mapsto f(a,j(t)b)$ belongs to
$\sC(p(t),\{X\})$. By lemma~\ref{comp}, we see that it suffices to impose 
this condition on all $a\in D(t)$, rather than all $a\in sD(t)$.  

\begin{sloppypar}
We observe that for every $a\in e(t)$, the
map $\ET(p(t))\to\ET(V)$ given by $b\mapsto (a,j(t)b)$ belongs to
$\sC(p(t),\{X\})$. As a consequence, we see that the inclusion $D\hookrightarrow\ET(V)$
is of class $\sC$.
\end{sloppypar}

 When concerned with $(x,y)$-compatible maps, we will assume that
$D(t)=\emptyset$ whenever $t$ and $\{x,y\}$ are disjoint.
 Under this assumption,
a map $f:D\to \ET(V)$ is said to be $(x,y)$-compatible of class $\sC$ 
if  $\ET(p(t))\to\ET(V)$ given by $b\mapsto f(a,j(t)b)$ is 
a $(x,y)$-compatible map in $\sC(p(t),\{X\})$  
 for all pairs $(a,t)$ such that   
 $a\in sD(t)$.  

In a similar manner, we define 
 $(x,y,z)$-compatible maps of class $\sC$ as well. For this, it is
necessary to assume that $D(t)$ is empty whenever the partition $p(t)$
 does not separate $(x,y,z)$, equivalently if $\{x,y,z\}\setminus t$ has
at least two elements.

  \begin{lemma}\label{final} Assume furthermore that $D(t)$ is
a simplicial subcomplex of $e(t)$.  
Then the space of maps  $D\to\ET(V)$ in class $\sC$ is 
  contractible. The same is true of the space of such maps that  are $(x,y)$-compatible,
  or $(x,y,z)$-compatible.
  \end{lemma}
\begin{proof} 
We denote by $d$ the cardinality
of $\{t:D(t)\neq \emptyset\}$. We proceed by induction on $d$, beginning
with $d=0$ where the space of maps is just one point.

We choose $t_0$ of maximum cardinality so that $D(t_0)\neq \emptyset$.
Let $D'$ be the union of $D(t)\times \ET(V/V(t))$ taken over all
$t\neq t_0$. Let $\sC(D')$ and $\sC(D)$ denote the space of 
class $\sC$ maps $D'\to\ET(V)$ and $D\to\ET(V)$ respectively. By
the induction hypothesis, $\sC(D')$ is contractible. We observe that
the intersection of $D'$ and   $e(t_0)\times \ET(V/V(t_0))$  has the form 
$G\times\ET(V/V(t_0)$
 where $G\subset e(t_0)$ is a subcomplex. Furthermore, $G\cup D(t_0)$
 is the saturated set $sD(t_0)$ described earlier.

For a closed subset $H\subset e(t_0)$, denote the space of $\sC$-maps
$H\times\ET(V/V(t_0))\to\ET(V)$ by $A(H)$. Note that 
$A(H)=\Maps(H,\sC(p(t_0),\{X\}))$. By
 lemma~\ref{silly}, the space $\sC(p(t_0),\{X\})$ is itself contractible. 
It follows that $A(H)$ is contractible. In particular, both
$A(G)$ and $A(sD(t_0))$ are contractible. The natural map
$A(sD(t_0))\to A(G)$ is a fibration, because the inclusion $G\hookrightarrow
 sD(t_0)$ is a cofibration. The fibers of  $A(sD(t_0))\to A(G)$ are thus
contractible. 
It follows that 
\[A(sD(t_0))\times_{A(G)}\sC(D')\to \sC(D')
\] 
which is simply $\sC(D)\to\sC(D')$, enjoys the same properties: it is
also a fibration with contractible fibers. Because $\sC(D')$ is contractible,
 we deduce that $\sC(D)$ is itself contractible. This completes the 
proof of the first assertion of the lemma.

 The remaining assertions follow in exactly the same manner by appealing to lemma~\ref{xy}.
\end{proof}
\noindent
 \emph {Proof of Proposition~\ref{three}}.  Choose $x\neq y$ with $x,y\in q$. Let
$g=id_V+\alpha$ where $\alpha(V)\subset V_y$ and $\alpha(V_k)=0$ for
 all $k\neq x\in X$. To prove the proposition, it suffices to show that $g\circ i$ is homotopic to $i$ where
 $i:U(q)\to\ET(V)$ is the given inclusion. \emph{This notation $x,y,\alpha, g$
 will remain fixed throughout the proof}.
\\  \\Case 1. Here $q=\{x,y\}$. 
Now  $x,y$ are separated by the partitions $p(t)$ for every non-empty
$t\subset q$.
 By the second assertion of the above lemma, 
 there exists $f:U(q)\to\ET(V)$ of class $\sC$ and $(x,y)$-compatible. 
 The given inclusion
 $i:U(q)\to\ET(V)$ is also of class $\sC$. By the first assertion of the 
same lemma, $f$ is homotopic to $i$. Now the image of $f$ is contained in the
fixed-points of $g$ and so we get 
$g\circ f=f$. It follows that 
 $g\circ i$ is homotopic to $i$.
This completes the proof of the proposition when $1=r=\#(q)-1$.
 \\  \\ Case 2. Here  $q=\{x,y,z\}$ with $x,y,z$ all distinct. 

We take $Y_1$ to be the union of $e(t)\times \ET(V/V(t))$ taken over
all $t\subset q,t\neq \{z\},t\neq \emptyset$. We put  
$Y_2=\ET(V_z)\times\ET(V/V_z)$ and $Y_3=Y_1\cap Y_2$.
We note that $U(q)=Y_1\cup Y_2$.

The given inclusion $i:U(q)\to\ET(V)$ restricts to $i_k:Y_k\to\ET(V)$ for
$k=1,2,3$.  The required homotopy is a path $\gamma:I\to \Maps(U(q),\ET(V))$ so that 
$\gamma(0)=i$ and
$\gamma(1)=g\circ i$.  Equivalently we require paths 
$\gamma_k$ in $\Maps(Y_k,\ET(V))$
for $k=1,2$ so that
\\(a) $\gamma_k(0)=i_k$ and $\gamma_k(1)=g\circ i_k$  for $k=1,2$ and
\\(b) both $\gamma_1$ and $\gamma_2$ restrict to the same path in 
$\Maps(Y_3,\ET(V))$. 
 
 In view of the fact that $Y_3\hookrightarrow Y_1$ is a cofibration, the weaker conditions $(a')$ and $(b')$ on fundamental groupoids suffice for
the existence of such a $\gamma$:
 \\$(a')$: $\gamma_k\in\pi_1(\Maps(Y_k,\ET(V)); i_k,g\circ i_k)$ for $k=1,2$
 \\$(b')$: both $\gamma_1$ and $\gamma_2$ restrict to the same element of 
$\pi_1((\Maps(Y_3,\ET(V));i_3,g\circ i_3)$

We have the spaces:
$Z_k=\Maps(Y_k,\ET(V))$ for $k=1,2,3$.These spaces come equipped with the
data below:
\\(A)The $GL(V)$-action on
$\ET(V) $ induces a $GL(V)$-action on $Z_k$
\\(B)The maps of class $\sC$ give contractible subspaces $\sC_k\subset Z_k$ for $k=1,2,3$.
\\(C) We have $i_k\in\sC_k$ for $k=1,2,3$.
\\(D) The natural maps $Z_k\to Z_3$ for $k=1,2$ are $GL(V)$-equivariant,
they take $i_k$ to $i_3$ and restrict to maps $\sC_k\to\sC_3$.

Note that the $GL(V)$-action on $Z_k$ turns the disjoint union:
\\$\sG_k=\sqcup\{\pi_1(Z_k;i_k,hi_k)|h\in GL(V)\}$ into a group: given
 ordered pairs $(h_j,v_j)\in\sG_k$, i.e. $h_j\in GL(V)$ and 
$v_j\in\pi_1(Z_k;i_k,h_ji_k)$ for $j=1,2$, we 
get $h_1v_2\in\pi_1(Z_k;h_1i_k,h_1h_2i_k)$ and obtain thereby
$v=(h_1v_2).v_1\in \pi_1(Z_k;i_k,h_1h_2i_k)$ and this produces
the required binary operation $(h_1,v_1)*(h_2,v_2)=(h_1h_2,v)$.

The projection $\sG_k\to GL(V)$ is a group homomorphism. The following
elementary remark will be used in an essential manner when checking
condition $(b')$. 
The data $(H,F,\Delta)$ where 
\\(i) $H\subset GL(V)$ is a subgroup,
\\(ii)$F\in Z_k$ is a fixed-point of $H$, and 
\\(iii)$\Delta\in \pi_1(Z_k;F,i_k)$
\\  
produces the lift $H\to \sG_k$ of the inclusion 
$H\hookrightarrow GL(V)$ by $h\mapsto
(h,(h\Delta).\Delta^{-1})$.
Finally we observe that there are natural homomorphisms $\sG_k\to\sG_3$
induced by $Z_k\to Z_3$ for $k=1,2$.
\\ \emph{Construction of $\gamma_1$. }

 The partitions $p(t)$ for $t\neq\{z\}$ separate $x,y$. 
By lemma~\ref{final}, we have a $(x,y)$-compatible class $\sC$-map
$f:Y_1\to \ET(V)$. Both $i_1$ and $f$ belong to $\sC_1$ and
thus we get $\delta\in\pi_1(\sC_1;f,i_1)$. Now $f$ is fixed by our 
$g\in GL(V)$, so we also get $g\delta\in \pi_1(g\sC_1;f,gi_1)$. The
path $(g\delta).\delta^{-1}$ is the desired $\gamma_1\in\pi_1(Z_1;i_1,gi_1)$.
\\  \emph{Construction of $\gamma_2$. }

Recall that $g=id_V+\alpha$. We choose $m:V_x\to V_z$ and
  $n:V_z\to V_y$ so that $nm(a)=\alpha(a)$ for all $a\in V_x$. We extend $m,n$ by
  zero to nilpotent endomorphisms of $V$, once again denoted by $m,n:V\to V$ and
  put $u=id_V+n,v=id_V+m$ and note that $g=uvu^{-1}v^{-1}$. 
  
Note that the partition $p(\{z\})$ separates 
both the pairs $(x,z)$ and  $(z,y)$. We thus obtain $f',f''\in\sC_2$
so that $f'$ is $(x,z)$-compatible and $f''$ is $(y,z)$-compatible
and also $\delta'\in\pi_1(\sC_2;f',i_2)$ and 
$\delta''\in\pi_1(\sC_2;f'',i_2)$. Noting that $f',f''$ are 
fixed by $v,u$ respectively, we obtain
\\$\epsilon'=(v\delta').\delta'^{-1}\in \pi_1(Z_2;i_2,vi_2)$
and \\$\epsilon''=(u\delta'').\delta''^{-1}\in \pi_1(Z_2;i_2,ui_2)$.

Thus $v'=(v,\epsilon')$ and $u'=(u,\epsilon'')$ both belong to $\sG_2$.
We obtain $\gamma_2$ by 
\\$u'*v'*u'^{-1}*v'^{-1}=(g,\gamma_2)\in\sG_2$
\\ \emph {Checking the validity of $(b')$}.

Let $\gamma_{13},\gamma_{23}\in\pi_1(Z_3;i_3,gi_3)$ be the images 
of $\gamma_1$ and $\gamma_2$ respectively. We have to show that
$\gamma_{13}=\gamma_{23}$.

Consider the spaces $\sH,\sH',\sH''$ consisting of ordered pairs
$(f_3,\delta_3)$, $(f'_3,\delta'_3)$, $(f''_3,\delta''_3)$ respectively,
where $f_3,f_3',f''_3$ are all in $\sC_3$,
\\ $f_3$ is $(x,y)$-compatible, $f_3'$ is $(x,z)$-compatible, and
$f_3''$ is $(y,z)$-compatible, and
\\$\delta_3,\delta_3',\delta_3''$ are all paths in $\sC_3$ that
originate at $f_3,f_3',f_3''$ respectively, and they all
terminate at $i_3$. By lemma~\ref{final}, we see that
the spaces $\sH,\sH',\sH''$ are all contractible.

For $t=\{x,z\},\{y,z\},\{x,y,z\}$, the partition $p(t)$ separates
 $(x,y,z)$.
Note that $Y_3$ is contained in the union of these three $\ET'(t)$.
By lemma~\ref{final},
there is a $(x,y,z)$-compatible $F\in\sC_3$. Let $\Delta$ be a path
in $\sC_3$ that originates at $F$ and terminates at $i_3$. We
see that $(F,\Delta)\in\sH\cap\sH'\cap\sH''$.

Note that $\sH\to \pi_1(Z_3;i_3,gi_3)$ given  by
$(f_3,\delta_3)\mapsto (g\delta_3).\delta_3^{-1}$
 is a constant map 
because $\sH$ is contractible.  
The $(f,\delta)$
employed in the construction of $\gamma_1$ restricts to an element
of $\sH$. Also,  $(F,\Delta)$ belongs to $\sH$. 
 It follows that $\gamma_{13}= (g\Delta).\Delta^{-1}$.

In a similar
manner, we deduce that if $\epsilon'_3,\epsilon''_3
$ denote the
images of $\epsilon',\epsilon''$ in the fundamental groupoid
of $Z_3$, then
\\$\epsilon_3'=(v\Delta).\Delta^{-1}\in \pi_1(Z_3;i_3,vi_3)$
and $\epsilon''_3=(u\Delta).\Delta^{-1}\in \pi_1(Z_3;i_3,ui_3)$

Thus $\sG_2\to\sG_3$ takes $v',u'\in\sG_2$ to 
$(v,(v\Delta).\Delta^{-1}), (u, (u\Delta).\Delta^{-1})\in\sG_3$
 respectively. It follows that their commutator $[u',v']$ maps to 
$(g,\gamma_{23})\in\sG_3$
under this homomorphism. 

We apply the remark preceeding the construction of $\gamma_1$ 
to the subgroup $H$ generated by $u,v$ and  
$F$ and $\Delta$ as above. We conclude that $\gamma_{23}$
equals $(g\Delta).\Delta^{-1}$. That the latter
equals $\gamma_{13}$ has already been shown. Thus $\gamma_{13}
=\gamma_{23}$ and this 
completes the 
proof
of the Proposition. 
\section{Low dimensional stabilisation of homology}
This section contains applications of corollary~\ref{elementary},
 proposition~\ref{three} and 
 Theorem~\ref{ss} to obtain some mild information on the homology
groups of $\ET(V)$.  The notation $\sL(V),\sL_r(V),W(q),
\det(q)$ introduced to state Theorem~\ref{ss}
 will be freely used throughout. The spectral sequence in theorem~\ref{ss} with
 coefficients in an Abelian group $M$ will be denoted by $SS(V;M)$. When
  $V=A^n$, this is further abbreviated to $SS(n;M)$, or even to $SS(n)$ when it
  is clear from the context what $M$ is.

The concept of a commutative ring with \emph{many units} is due to
Van der Kallen. An exposition  of the definition and consequences 
of this term is given in \cite{Mirzaii}. We note that this class of rings includes semilocal
rings with infinite residue fields. The three consequences of this hypothesis on $A$
are listed as I,II,III below. These statements are followed by some elementary
 deductions. \emph{Throughout this section, we will assume that our ring $A$
 has this property.}
\\ \\ \textbf{I}: $SL_n(A)=E_n(A)$.   

This permits a better formulation of Lemma~\ref{contiguity} in many
instances. 
\\ \textbf{Ia}: Let $0\to W\to P\to Q\to 0$ be an exact sequence of free $A$-modules with
of ranks $a,a+b,b$. Let $d=g.c.d.(a,b)$.  Let $H$  be the group of automorphisms of this exact sequence 
 for which the induced automorphisms  on  $W$ and $Q$ are of the type $\alpha.id_W$ and
 $\beta.id_Q$ respectively where $\alpha, \beta$ are arbitrary units of $A$.  We may regard
 $H$ as a subgroup of $GL(P)$. 
This group $H$  acts trivially on the image of the embedding 
$i:\ET(W)\times\ET(Q)\to\ET(P)$. Furthermore 
$\{\det(g)|g\in H\}$  equals  $(A^{\times})^d$. Thus, if
$a,b$ are relatively prime, by lemma~\ref{contiguity}, we see that
$g\circ i$ is freely homotopic to $i$ for all $g\in GL(V)$.

We shall take rank$(W)=1$ in what follows.
 Here $\ET(W)\times\ET(Q)$ is canonically
identified with $\ET(Q)$. The induced $\ET(Q)\to\ET(P)$
 gives rise on homology to   an arrow $H_m(\ET(Q))\to H_m(\ET(P))$  
which has a factoring: 
\begin{sloppypar}
$H_m(\ET(Q))\twoheadrightarrow H_0(PGL(Q), H_m(\ET(Q)))\to H^0(PGL(P), H_m(\ET(P))\hookrightarrow H_m(\ET(P)).$
\end{sloppypar}
The kernel of $H_m(\ET(Q))\to H_m(\ET(P))$ 
does not depend on the choice of the exact sequence. Denoting this kernel by
$KH_m(Q)\subset H_m(\ET(Q))$ 
therefore gives rise to unambiguous
notation.
We abbreviate $H_m(\ET(A^n)),KH_m(A^n)$
to $H_m(n),KH_m(n)$ respectively.
\\ \textbf{Ib}:  In the spectral sequence $SS(n)$, we have:
\begin{enumerate}
\item  $H_0(PGL_{n-1}(A),H_m(n-1))\cong H_0(PGL_n(A),E^{1}_{0,m})$.
\item    $E^{\infty}_{0,m}$ is                                the image
of $H_m(n-1)\to H_m(n)$. 
\item $E^1_{0,m}\to E^{\infty}_{0,m}$  factors as follows: 
\begin{center}$E^1_{0,m}\twoheadrightarrow
 E^2_{0,m}\twoheadrightarrow
H_0(PGL_n(A),E^1_{0,m})
\twoheadrightarrow 
E^{\infty}_{0,m}$.
\end{center}
\item
If $H_0(PGL_n(A),E^p_{p,m+1-p})=0$ for all
$p\geq 2$,  then 
\\the given arrow $H_0(PGL_n(A),E^1_{0,m})\to E^{\infty}_{0,m}$ is an 
isomorphism.
\item Assume that $H_m(n-2)\to H_m(n-1)$ is surjective. Then 
the arrow $E^2_{0,m}\to H_0(PGL_n(A),E^1_{0,m})$ in (3) above is
an isomorphism.
\end{enumerate}
The factoring in part (3) above is a consequence of the factoring of  
$H_m(\ET(Q))\to H_m(\ET(P))$  in part I(a).  

For part (4),
one notes that the composite 
\\$E^2_{2,m-1}\to E^2_{0,m}\to   H_0(PGL_n(A),E^1_{0,m})$ vanishes
because $H_0(PGL_n(A),E^2_{2,m-1} )$ itself vanishes. Thus we obtain
a factoring: $E^2_{0,m}\to  E^3_{0,m}\to   H_0(PGL_n(A),E^1_{0,m})$.
Proceeding inductively, we obtain the factoring:
 \\$E^2_{0,m}\to  E^{\infty}_{0,m}\to   H_0(PGL_n(A),E^1_{0,m})$. In view of (3),
 we see that part (4) follows.

For part (5), it suffices to note that for every $L_0,L_1\in\sL(A^n)$
with $n>1$, there is some $L_2\in\sL(A^n)$ with the property that
both $\{L_0,L_2\}$ and $\{L_1,L_2\}$ belong to $\sL_1(A^n)$. This fact
is contained in consequence III of \emph{many units}.
\\ \\ \textbf{II}: $A$ is a Nesterenko-Suslin ring. 
 
  Let $r>0,p\geq 0$. Put $N=(p+1)!$ and $n=r+p+1$. Let $\sF_r$ be the category with free 
   $A$-modules of rank $r$ as objects; the
   morphisms in $\sF_r$ are $A$-module isomorphisms.  
  Let $F$ be a functor from $\sF_r$  to the category  of $\Z[\frac{1}{N}]$ -modules. Assume that
  $F(a.id_D)=id_{FD}$ for every $a\in A^{\times}$ and for every object $D$ of $\sF_r$.  In other
  words, the natural action of $GL_r(A)$ on $F(A^r)$ factors through the action of
   $PGL_r(A)$.  
  
  For a 
  free $A$-module $V$ of rank $n$, define $\mathit{Ind}'F(V)$ by 
\begin{center}$\mathit{Ind}'F(V)=\oplus\{\det(q)\otimes F(V/W(q)):q\in \sL_p(V)\}$. 
\end{center}   
An alternative description of $\mathit{Ind}'F(V)$ is as follows. Fix some $q\in\sL_p(V)$. Let 
$G(q)$ be the stabiliser of $q$ in $GL(V)$. Then  $\det(q)$ and $F(V/W(q))$ are $G(q)$-modules
 in a natural manner. We have a natural isomorphism of $\Z[GL(V)]$-modules:
 \begin{center}$\mathit{Ind}'F(V)\cong\Z[GL(V)]\tensor_{\Z[G(q)]}[\det(q)\otimes _{\Z}F(V/W(q))]$. 
\end{center}    
 \textbf{IIa}: If $H_i(PGL_r(A),F(A^r))=0$ for all $i<m$, then $H_i(PGL_n(A),\mathit{Ind}'F(A^n))=0$ 
for all $i<p+m$. Furthermore,
\begin{center} $H_{p+m}(PGL_n(A), \mathit{Ind}'F(A^n))\cong H_m(PGL_r(A),F(A^r))\otimes Sym^p(A^{\times}))$
\end{center}

By Shapiro's lemma, the result  of \cite{NS} cited earlier, and the fact that 
the group homology  $H_i(M, C)$ is isomorphic to $ C\otimes\Lambda^i(M)$ for all commutative groups $M$ and
$\Z[1/i!]$-modules $C$ given  trivial $M$-action,
IIa  reduces to the statement: 
\\ Let $\Sigma(q)$ denote the group of permutations of a set $q$ of $(p+1)$ elements.
Let  $M$ be an Abelian group on which $(p+1)!$ acts
invertibly.  Then 
\\ $H_0(\Sigma(q), \det(q)\otimes \Lambda^i(M^q))$ vanishes when $i<p$ and is isomorphic to
$Sym^p(M)$ when $i=p$. 
\\ \\ \textbf{III}: The standard application of the \emph{many units}  hypothesis, see
( \cite{SuslinMilnor} for instance) is that general 
position  is available in the precise sense given below.
Let $V\cong A^n$. 
We denote by $K(V)$ the simplicial complex whose set of vertices is $\sL(V)$. 
A subset $S\subset \sL(V)$ of cardinality $(r+1)$ 
is an $r$-simplex of $K(V)$ if every $T\subset S$ of cardinality $t+1\leq n$
 belongs to $\sL_t(V)$. If $L\subset K(V)$ is a finite simplicial subcomplex, then
there is some $e\in \sL(V)$ with the property that 
$s\cup\{e\}$ is a $(r+1)$-simplex of $K(V)$ for every $r$-simplex $s$
of $L$. This gives an embedding $Cone(L)\hookrightarrow K(V)$ with $e$ as the vertex 
of the cone. Thus $K(V)$ is contractible.  The complex of oriented chains of this simplicial
complex will be denoted by $C_{\bullet}(V)$.  Thus the reduced homologies
$\tilde{H}_i(C_{\bullet}(V))$ vanish for all $i$. 

The group 
 $D(V)=Z_{n-1}C_{\bullet}(V)=B_{n-1}C_{\bullet}(V)$ comes up frequently.
\\ \textbf{IIIa}:

\begin{enumerate}
\item Let $p<n$. Let $M$ be a $\Z[1/N]$-module where $N=(p+1)!$. Then 
$H_j(PGL(V),C_p\otimes M)$ vanishes for $j<p$ and is isomorphic to
$Sym^p(A^{\times})\otimes M$ when $j=p$. 
\item  $H_j(PGL(V),C_n\otimes M)$ vanishes for all $j\geq 0$ and for all 
$\Z[1/(n+1)!]$-modules $M$.
\item  $H_0(PGL(V), D(V)\otimes M)$  for any
 Abelian group $M$ is  isomorphic to $M/2M$ if $n$ is even, 
and vanishes if when $n$ is odd
\item $H_0(PGL(V),B_pC_{\bullet}(V)\otimes M)=0$ for every $\Z[1/2]$-module
$M$ and for every $0\leq p<n$.
\item  $H_1(PGL(V),Z_pC_{\bullet}(V)\otimes M)=0$ for every $\Z[1/(p+2)!]$-
module $M$ and $1\leq p\leq n-2$.
\end{enumerate}

Note that (1) above follows from IIa when $F$ is the constant functor $M$.

 For (2), one observes
that $PGL(V)$ acts transitively on the set of $n$-simplices of the
simplicial complex $K(V)$. The stabiliser of an $n$-simplex is
the permutation group $\Sigma$ on $(n+1)$ letters. The claim now follows from
Shapiro's lemma.

The presentation $C_{p+2}(V)\otimes M\to C_{p+1}(V)\otimes M
\to B_pC_{\bullet}(V)\otimes M\to 0$ and the observation
$H_0(PGL(V),  C_{p+1}(V)\otimes M)\cong M/2M$ whenever $p<n$ suffice
to take care of (3) and (4).

For assertion (5), one applies the long exact sequence of group homology
to the short exact sequence:
\\$0\to Z_{p+1}C_{\bullet}(V)\otimes M\to
C_{p+1}(V)\otimes M\to  Z_{p}C_{\bullet}(V)\otimes M\to 0$. 
\begin{sloppypar}
One
therefore obtains the exact sequence: 
\\
$H_1(PGL(V),C_{p+1}(V)\otimes M)\to  H_1(PGL(V),Z_{p}C_{\bullet}(V)\otimes M)
\to H_0(PGL(V),Z_{p+1}C_{\bullet}(V)\otimes M)$. The end terms here
vanish by (1) and (4).
\end{sloppypar}
\noindent  \textbf{IIIb}:  
\begin{enumerate}
\item     $\ET(V)$ is connected.
\item  $E^2_{0,0}=\Z,E^2_{n-1,0}=D(V),E^2_{m,0}=0$ if $m\neq 0,n-1$ for
the spectral sequence $SS(V)$.
\item  $H_1(\ET(V))\cong \Z/2\Z$ if rank$(V)>2$.

\end{enumerate}
Note that (1) is a consequence of (2).  Part (2) is deduced by induction on  rank$(V)=n$. The induction hypothesis 
 enables the 
identification of  the
$E^1_{m,0}$ terms of the spectral sequence for $V$ (together with differentials) 
 with the $C_m(V)$  (together with  boundary operators) when $m<n$.
Thus (2) follows. 

For part (3), consider the spectral sequence  $SS(3)$. Here $E^2_{2,0}=D(A^3)$ and
$H_0(PGL_3(A),D(A^3))=0$ by IIIa(3). Thus the hypothesis of Ib(4) holds for SS(3). 
Consequently, 
\\$H_1(3)\cong E^{\infty}_{1,0}\cong H_0(PGL_2(A), H_1(2))=H_0(PGL_2(A),D(A^2))\cong\Z/2\Z$,
\\the last isomorphism given by IIIa(3) once again. The isomorphism $H_1(n)\cong\Z/2\Z$
 for $n>3$ is contained in the lemma below for $N=1$.

\begin{lemma}\label{easystab}
Let $M$ be an Abelian group. Let $N\in\N$.  For $0<r<N$, we are given $m(r)\geq 0$ so that 
$H_r(d;M)\to H_r(d+1;M)$ 
\\is a surjection if $d=r+m(r)+1$ , and 
\\an isomorphism if 
$d>r+m(r)+1$. 

 Let $m(N)=\max\{0, m(1)+1,m(2)+1,...,m(N-1)+1\}$. 
\begin{sloppypar}
Then
$H_N(d;M)\to H_N(d+1;M)$ is 
\\(a) an isomorphism if $d>N+m(N)+1$,
\\(b) is a  surjection if $d=N+m(N)+1$.  
\\(c) The surjection in (b) above  factors through an isomorphism
 $H_0(PGL_d(A),H_N(d))\to H_N(d+1;M)$ if 
if $M$ is
a $\Z[1/2]$-module.
\end{sloppypar}
\end{lemma}
\begin{sloppypar}
\begin{proof} Consider the spectral sequence $SS(V;M)$ that computes
the homology of $\ET(V)$ with coefficients in $M$. Here $V$ is free of rank
 $N+h+2$, where $h\geq m(N)$. We make the following claim:
\\ \emph{Claim:}If $E^2_{s,r}\neq 0$ and $0<s$ and $r<N$,
 then $s+r\geq N+1+h-m(N)$.
Furthermore, when equality holds, $H_0(PGL(V),E^2_{s,r}\otimes\Z[1/2])=0$.

We assume the claim and prove the lemma. We take $h=m(N)$.
 All the $E^2_{s,r}$ with $s+r=N$ are zero except possibly for 
$(s,r)=(0,N)$. Part (b) of the lemma now follows from Ib(2). We 
consider next
the $E^2_{s,r}$ with $s+r=N+1$ and $s\geq 2$ (or equivalently with
$r<N$). It follows that $E^s_{s,r}$ is a quotient of $E^2_{s,r}$. The
second assertion of the claim now show that
 $H_0(PGL(V),E^s_{s,r})=0$ if $M$ is a $\Z[1/2]$-module. Part (c)
 of the lemma now follows from Ib(4). 

We take $h>m(N)$ and prove part (a) by induction on $h$. 
The inductive hypothesis implies that $H_N(N+h;M)\to H_N(N+h+1;M)$
 is surjective. By Ib(5), it follows that $H_N(N+1+h;M)\to E^2_{0,N}$
 is an isomorphism.  Now there are no nonzero
$E^2_{s,r}$ with $s+r=N+1$ and $s\geq 2$. Thus $E^2_{0,N}=E^{\infty}_{0,N}$.
It follows that $H_N(N+1+h;M)\to H_N(N+2+h;M)$ is an isomorphism. 

It only remains to prove the claim. We address this matter now.

For $r=0$, both assertions of the claim are valid by IIIb(2)and IIIa(3).

\begin{sloppypar}
So assume now that $0<r<N$. Let  $SH(r)=H_r(d;M)$ for $d=r+m(r)+2$. In view of our hypothesis,
the chain complex 
\\$E^1_{0,r}\leftarrow E^1_{1,r}\leftarrow...\leftarrow E^1_{p-1,r}\leftarrow E^1_{p,r}$ 

for $N+h+1=p+r+m(r)+1$ is identified with 
\\$C_0(V)\otimes SH(r)\leftarrow...\leftarrow C_{p-1}(V)\otimes SH(r)
\leftarrow \oplus\{\det(q)\otimes H_r(\ET(V/W(q));M)|q\in\sL_p(V)\}$.
 
 As in IIIb(2), it follows that $E^2_{s,r}=0$ whenever $0<s<p$. Furthermore, we deduce 
 the following  exact sequence for  $E^2_{p,r}$:
 \\$\oplus \{\det(q)\otimes KH_r(V/W(q);M)|q\in\sL_p\}\to E^2_{p,r}\to Z_pC_{\bullet}\otimes SH(r)\to 0 $.
  By IIIa(4), we see that $H_0(PGL(V), E^2_{p,r})=0$ if $M$ is a $\Z[1/2]$-module.  

Note that  $p+r=N+h-m(r)\geq N+1+h-M(N)$. 
This completes the proof of the claim,
 and therefore, the proof of the lemma as well.
\end{sloppypar}

\end{proof}
\end{sloppypar}
The Proposition below is an application of Proposition~\ref{three}.
The notation here is that of Theorem~\ref{ss}. We regard $B^{r}_{p,q}$ and
$Z^{r}_{p,q}$ as subgroups of $E^1_{p,q}$ for all $r>1$. The notation $KH_m(Q)$
 has been introduced in Ia, the first application of \emph{many units}.
\begin{proposition}\label{kvanish} Let rank$(V)=n$. Let $M$ be a 
$\Z[1/2]$-module. In the spectral sequence $SS(V;M)$, we have:
\begin{enumerate}
\item $\oplus\{\det(q)\otimes KH_m(V/W(q))|q\in\sL_1(V)\}\subset B^{\infty}_{1,m}$ if $n>1$.
\item $E^{\infty}_{1,m}=0$ if $n>2$ and $M$ is a $\Z[1/6]$-module.
\item If, in addition, it is assumed that 
\\$H_{m+1}(n-2;M)\to H_{m+1}(n-1;M)$ is surjective, then 
\\$\oplus\{\det(q)\otimes KH_m(V/W(q))|q\in\sL_2(V)\}\subset B^{\infty}_{2,m}$.

\end{enumerate} 
\end{proposition} 
\begin{proof}  
Let $q\in \sL_r(V)$. We have $U(q)\subset \ET(V)$ as in Proposition~\ref{three}. The
spectral sequence of Theorem~\ref{ss} was constructed from an increasing filtration 
of subspaces of $\ET(V)$. Intersecting this filtration with $U(q)$ we
obtain a spectral sequence that computes the homology of $U(q)$.
Its terms will be denoted by $E^a_{b,c}(q)$. One notes that 
$E^1_{b,m}(q)$ is the direct sum of $\det(u)\otimes H_m(\ET(V/W(u)))$
taken over all $u\subset q$ of cardinality $(b+1)$.

We denote the terms of the spectral sequence in theorem~\ref{ss}
by $E^a_{b,c}(V)$. The given data also provides a homomorphism
$E^a_{b,c}(q)\to E^a_{b,c}(V)$ of $E^1$-spectral
sequences. We assume that $M$ is a $\Z[1/(r+1)!]$-module.

We choose a basis $e_1,e_2,...,e_n$ of $V$ 
so that 
\\$q=\{Ae_i:1\leq i\leq r+1\}$. Let $G\subset GL(V)$ be the subgroup
of $g\in GL(V)$ so that 
\\(A) $g(q)=q$, (B)$g(e_i)=e_i$ for all $i>r+1$,
(C), the matrix entries of $g$ are $0,1,-1$ and (D) $\det(g)=1$.
Now $G$ acts on the pair $U(q)\subset\ET(V)$. Thus the above
homomorphism of spectral sequences is one such in the category 
of $G$-modules.We observe:
\\(a) $G$ is a group of order $2(r+1)!$
\\(b) there are no nonzero $G$-invariants in $E^1_{i,m}(q)$ for $i>0$,
and consequently the same holds for  
all $G$-subquotients, in particular for $E^a_{i,m}(q)$ for all $a>0$ as
well. 
\\ \emph{Proof of part 1}. Take $r=1$.  Proposition~\ref{three} implies that the image of $H_m(U(q))\to H_m(\ET(V))$
 has trivial $G$-action. In view of (b) above, this shows that
$E^{\infty}_{1,m}(q)\to E^{\infty}_{1,m}(V)  $ is zero. But
$E^{\infty}_{1,m}(q)=Z^{\infty}_{1,m}(q)=\det(q)\otimes KH_m(V/W(q))$. It follows that
$\det(q)\otimes KH_m(V/W(q))\subset B^{\infty}_{1,m}(V)$. Part (1)
 follows.
\\ \emph{Proof of part (2)}. We take $r=2$. Here we have 
\\ $Z^{\infty}_{1,m}(q)=Z^{2}_{1,m}(q)$. Appealing to Proposition~\ref{three}
and observation (b) once again, we see that the image of the
homomorphism $Z^2_{1,m}(q)\to Z^2_{1,m}(V)$ is contained in $B^{\infty}_{1,m}$.
Part (2) therefore follows from the claim below.
\\ \emph{Claim:} $\oplus\{Z^2_{1,m}(q)|q\in\sL_2(V)\}\to Z^2_{1,m}(V)$
 is surjective.

Denote the image of $H_m(n-2;M)\to H_m(n-1;M)$ by $I$.
A simple computation produces the exact sequences:
\\$0\to \oplus\{\det(u)\otimes KH_m(V/W(u):u\in\sL_1(V),u\subset q\}\to Z^2_{1,m}(q)
\to \det(q)\otimes I\to 0$, and
\\$0\to \oplus\{\det(u)\otimes KH_m(V/W(u):u\in\sL_1(V)\}\to Z^2_{1,m}(V)
\to Z_1C_{\bullet}(V)\otimes I\to 0$.

The claim now follows from the above description of
 $Z^2_{1,m}(q)$ and $Z^2_{1,m}(V)$.
Thus part (2)
 is proved.
\\ \emph{Proof of part (3)}.We take $r=2$ once again. 
The surjectivity of $H_{m+1}(n-2;M)
\to H_{m+1}(n-1;M)$ implies that $E^2_{0,m+1}(q)$ has trivial $G$-action.
By observation (b), we see that $d^2_{2,m}:E^2_{2,m}(q)\to E^2_{0,m+1}(q)$
 is zero. It follows that $E^{\infty}_{2,m}(q)=Z^2_{2,m}(q)$ here.
 Proposition~\ref{three} and  observation (b) once again show that
 the image of $Z^2_{2,m}(q)\to Z^2_{2,m}(V)$ is contained in $B^{\infty}
_{2,m}(V)$. Because $Z^2_{2,m}(q)=\det(q)\otimes KH_m(V/W(q))$, part (3)
follows. 

This completes the proof of the Proposition.
\end{proof}

\begin{theorem}\label{stab}Let $H_m(n;M)$ denote $H_m(\ET(A^n);M)$ where
$M$ is a $\Z[1/6]$-module. We have:
\item(1) $H_1(n;M)=0$ for all $n>2$,
\item(2)  $H_0(GL_3(A),H_2(3;M))\to H_2(n;M)$ is an isomorphism for all $n\geq 4$,
\item(3)  $H_0(GL_4(A), H_3(4;M))\to H_3(n;M)$ is an isomorphism for all $n\geq 5$,
\item(4)  $H_0(GL_{2m-2}(A), H_m(2m-2;M))\to H_m(n;M)$ is an isomorphism for all $n>2m-2$.
\end{theorem} 
\begin{proof} \begin{sloppypar}
Part (1) has already been proved. 
\\ \emph{Proof of part 2}. For this, 
we study $SS(V;M)$ where rank$(V)=4$. We first note that 
\item(i) $E^2_{3,0}=D(V)$ and therefore $H_0(PGL(V),E^2_{3,0})=0$.
\item(ii) $E^2_{1,1}=E^1_{1,1}=\oplus\{\det(q)\otimes D(V/W(q)):q\in\sL_1(V)\}$, and therefore $H_1(PGL(V),E^2_{1,1})=0$ by IIa.
\item(iii) $E^2_{u,v}=0$ except when $(u,v)=(0,0),(0,2),(1,1),(3,0)$.
We have $E^{\infty}_{1,1}=0$ by
proposition~\ref{kvanish} and thus obtain the short exact 
sequence:
 \\$ 0\to E^3_{3,0}\to           E^2_{3,0}\to E^2_{1,1}\to 0$.

By (i) and (ii) above, we see that $H_0(PGL(V),E^3_{3,0})=0$. By Ib(1,4),
we see that $H_0(PGL_3(A),H_2(3;M))\to H_2(4;M)$ is an isomorphism. 
In particular, $H_2(4;M)$ receives the trivial $PGL_4(A)$-action.
Taking
$N=2$ and $m(1)=0 $ in lemma~\ref{easystab}, we see that
$H_2(4;M)\to H_2(n;M)$ is an isomorphism for all $n\geq 4$. This
proves part (2).
\\ \emph{Proof of part 3}. We inspect $SS(V;M)$ where $V=A^5$.
We note that
\item(1)  $E^{\infty}_{1,2}$ and $E^{\infty}_{2,1}$ both vanish. This
follows from  proposition~\ref{kvanish}, once it is noted that
$KH_1(2;M)=H_1(2;M)$. 
\item(2) $E^2_{0,2}=H_2(4;M)$ has the trival $PGL(V)$-action.
\item(3)  $H_i(PGL(V),E^2_{2,1})=0$ for all $i<3$. This follows
from IIa and IIIa(3) after observing that $H_1(2;M)\cong D(A^2)\otimes M$.
\item(4) From (2) and (3) we see that $d^2_{2,1}=0$.
\item(5) 
We deduce that $E^2_{2,1}\cong E^2_{0,4}/E^3_{0,4}$ and
$E^2_{1,2}=E^3_{1,2}\cong  E^3_{0,4}/E^4_{0,4}$ from observations (1) and (4).
\item(6) $H_1(PGL(V),E^2_{1,2})=0$. 

To see this, first note the 
the short exact sequence:
\\$0\to P\to E^2_{1,2}\to Q\to 0$, where
\\$P=\oplus\{\det(q)\otimes KH_2(V/W(q):q\in\sL_1(V)\}$ and 
$Q=H_2(4;M)\otimes Z_1C_{\bullet}(V)$.

The vanishing of
$H_1(PGL_3(A), Q)$ follows from IIIa(5). By IIa, the vanishing of
$H_1(PGL_3(A), P)$ is reduced to the vanishing of $H_0(PGL_3(A),KH_2(A^3)$.

Now let $I$ be the augmentation ideal of the group algebra
 $R[PGL_3(A)]$ where $R=\Z[1/6]$. In view of the fact that 
$PGL_3(A)_{ab}$ is 3-torsion, we see that $I=I^2$. It follows that 
for all $\Z[1/6]$-modules $N$ equipped with $PGL_3(A)$-action, 
we have $IN=I^2N$, or equivalently, $H_0(PGL_3(A),IN)=0$. We apply
 this remark to $N=H_2(3;M)$.  
By part (2) of the proposition, we see that 
$KH_2(3;M)=IN$. This proves that  $H_0(PGL_3(A),KH_2(A^3))=0$.
We have completed the proof of
observation 6.
\item(7) $H_0(PGL(V),E^4_{4,0})=0$. 

In view of the filtration of (5),
it suffices to check that $H_1(PGL(V), E^2_{a,b})=0$ for 
$(a,b)=(1,2)$ and $(2,1)$ (which has been seen in observations (3) and (6))  
and also that $H_0(PGL(V),E^2_{4,0})=0$ (and this is clear because
$E^2_{4,0}=D(V)$).
\item(8) $H_0(PGL_4(A),H_3(4;M))\to H_3(V;M)$ is an isomorphism. 

That $H_0(PGL_4(A),H_3(4;M))\to E^{\infty}_{0,3}$ is an isomorphism
 follows from 
observation (7) and Ib(4). Now  $E^{\infty}_{a,b}=0$
whenever $a+b=3$ and $(a,b)\neq (3,0)$. This proves (8). 
\\(9) $H_3(5;M)\to H_3(n;M)$
is an isomorphism for all $n\geq 5$.

This follows from lemma~\ref{easystab} by taking $N=3$ and $m(1)=m(2)=0$. 
This finishes the proof of part (3). 

Part (4) now follows from the same lemma
 and induction. 
\end{sloppypar}
\end{proof}

\begin{remark} It can be checked that parts (1,2,4) of the above theorem
are valid for $\Z[1/2]$-modules $M$. In part (3), it is true that
$H_3(n;M)\cong H_3(n+1;M)$ for $n>4$ and also that 
$H_0(PGL_4(A),H_3(4;M)\to H_3(5;M)$ is a surjection. 
\end{remark}
\begin{proposition}\label{generalstab} Assume that the Compatible Homotopy
 Question has an affirmative answer. Then, for all $\Z[1/r!]$-modules $M$
 and for all $d>r+1$, 
\\$H_0(PGL_{r+1}(A),H_r(r+1;M))\to H_r(d;M)$ is an isomorphism.

\end{proposition}
\begin{proof} For $r=1$, this statement has been checked in IIIb(3)
 and lemma~\ref{easystab}. Let $N>1$. 
We assume that the above statement has been proved for
all $r<N$. Let $M$ be a $\Z[1/N!]$-module. In lemma~\ref{easystab},
may may now take $m(1)=m(2)=...=m(N-1)=0$. From this lemma, we obtain:
\[H_0(PGL_{N+2}(A),H_N(N+2;M))\to H_N(N';M)\]
is an isomorphism for all $N'>N+2$. So the proposition is proved once it is
checked that 
\[H_0(PGL_{N+1}(A),H_N(N+1;M))\to H_N(N+2;M)\]
is an isomorphism. To prove this, we consider the spectral sequence 
 $SS(V;M)$ where $V=A^{N+2}$. We will prove:
\item (i) $E^2_{a,b}=0$ or $a=0$ or $a+b=N$ or $(a,b)=(N+1,0)$. Furthermore the action
of $PGL(V)$ on $E^2_{0,b}$ is trivial when $b<N$. 
\item (ii) if $a>0$ and $b>0$, then $ H_i(PGL(V),E^2_{a,b})=0$ for $i=0,1$.  
\item (iii) $E^{\infty}_{a,b}=0$ when $a>0$ and $b>0$. 
\item (iv) 
 $E^2_{a,b}\cong E^{b+1}_{N+1,0}/E^{b+2}_{N+1,0}$ whenever $a>0,b>0$ and $a+b=N$. 

We first observe that (iv) is true for any spectral sequence of $PGL(V)$-modules where 
(i),(ii) and (iii) hold.  Next note that (ii) and (iv) imply that 
$H_0(PGL(V),E^{N+1}_{N+1,0})$ is contained in $H_0(PGL(V), E^2_{N+1,0})$.  And
since the latter is zero, we see that the former also vanishes.

We deduce that both arrows 
 $H_0(PGL(V), E^1_{0,N})\to E^{\infty}_{0,N}\to H_N(N+2;M)$ are isomorphisms exactly as 
 in earlier proofs. Thus it only remains to prove (i), (ii) and (iii).
\\  \emph{Proof of (i)}. This is contained in the proof of lemma~\ref{easystab}.
\\  \emph{Proof of (ii)}.  $0\to P\to E^2_{a,b}\to Q\to 0$ is exact, where
\\$P=\oplus\{\det(q)\otimes KH_b(V/W(q)|q\in\sL_a(V)\}$, and $Q=Z_aC_{\bullet}\otimes H_b(b+2;M)$
 as in the proof of the lemma~\ref{easystab}.   The required vanishing of $H_i(PGL(V), T)$ for
 $i=0,1$ holds for  $T=Q$ by IIIa(5). For $T=P$ and $a>1$,
the required vanishing follows from II(a). For $T=P$ and $a=1$, this
 is deduced from the vanishing of $H_0(PGL_N(A),KH_{N-1})$ (see
the proof
 of observation (6) in the proof of theorem~\ref{stab}).
\\ \emph{Proof of (iii)}. We follow the steps of the proof of Proposition~\ref{kvanish}. We first choose $q\in\sL_{N-1}(V)$ and consider the inclusion $U(q)\hookrightarrow \ET(V)$.  As in that proof we get a homomorphism
of $E^1$ spectral sequences of $G$-modules with $G\subset SL(V)$
 as  given there. The terms of  these spectral sequences are denoted by 
$E^a_{b,c}(q)$ and $E^a_{b,c}(V)$ respectively. From the inductive hypothesis,
we deduce:

\item (i$'$) $E^2_{a,b}(q)=0$ or $a=0$ or $a+b=N$ or $(a,b)=(N+1,0)$. Furthermore the action
of $G$ on $E^2_{0,b}(q)$ is trivial when $b<N$.

\item (ii$'$) if $a>0,h>0$, then $ H_0(G,E^h_{a,b}(q))=0$.

These observations together imply 
\item (iii$'$)  $Z^{\infty}_{a,b}(q)=Z^2_{a,b}(q)$ when $a>0$
 and $b>0$.

For $a>0,b>0$, we obtain $E^{\infty}_{a,b}(q)\to E^{\infty}_{a,b}(V)$ is zero, from the affirmative answer to 
the Compatible Homotopy Question. For such $(a,b)$, 
\\the image of 
$x(q):Z^2_{a,b}(q)\to Z^2_{a,b}(V)$ is thus contained in $B^{\infty}_{a,b}(V)$.
As in the proof of proposition~\ref{kvanish}, we see that the 
sum of the images of
 $x(q)$, taken over all $q\in\sL_{N-1}(V)$, is all of $Z^2_{a,b}(V)$. It
follows that $Z^2_{a,b}(V)=B^{\infty}_{a,b}(V)$ and thus
$E^{\infty}_{a,b}(V)=0$ whenever $a>0,b>0$. This proves assertion (iii)
 and this completes the proof of the Proposition.

\end{proof}
\section {a double complex}
We will continue to assume that $A$ is a commutative ring with many units.
The paper \cite{BMS} of Beilinson, Macpherson and Schechtman 
introduces a Grassmann complex, intersection and
projection maps, and a torus action. The terms of the double-complex constructed below may be obtained from the quotients by the torus action of the objects of \cite{BMS}. The arrows
of the double-complex are signed sums of 
their intersection and projection maps.
 
$D(V),C_{\bullet}(V)$ etc. are as in the previous section.
When 
rank$(V)=n$, we have the resolution:
\\$0\leftarrow D(V)\leftarrow C_n(V)\leftarrow C_{n+1}(V)
...$.

We put $\oC_r(V)=H_0(PGL(V),C_r(V))$ when $r\geq n$ and define $\oC_r(V)$
 to be zero otherwise. We put $\oC_r(A^n)=\oC_r(n)$.
 We observe that the above resolution of $D(V)$ tensored with the rationals
 is a projective resolution in the category of $\Q[PGL(V)]$-modules.
 It follows that  $H_i(PGL_n(A),D(A^n))\otimes\Q\cong H_{n+i}(\oC(n)_{\bullet})
\otimes\Q$. We denote by $\partial':\oC_r(n)\to\oC_{r-1}(n)$ the
boundary operator of $\oC(n)_{\bullet}$. We will now define
$\partial'':\oC_{r}(n)\to\oC_r(n-1)$. 

Let $V\cong A^n$. Let $(L_0,L_1,...,L_r)$ be an ordered $(r+1)$-tuple in $\sL(V)$ 
that gives rise to a $r$-simplex of $K(V)$ (see consequence III 
of \emph{many units} for notation).
We define $\partial_i(L_0,L_1,...,L_r)\in C_{r-1}(V/L_i)$ 
by  $\partial_i(L_0,L_1,...,L_r)=(\oL_0,..,\oL_{i-1},\oL_{i+1},...,\oL_r)$ where 
$\oL_j=L_j+L_i/L_i\in\sL(V/L_i)$ whenever $j\neq i$. Now
let 
\[g_r(L_0,L_1,...,L_r)=\Sigma_{i=0}^r(-1)^i\partial(L_0,L_1,...,L_r)\in
\oplus\{C_{r-1}(V/L):L\in \sL(V)\}.\]
The above $g_r:C_r(V)\to \oplus\{C_{r-1}(V/L):L\in \sL(V)\}$ anti-commutes
with the boundary operator. The functor $M\to H_0(PGL(V),M)$ takes
$g_r$ to $\partial'':\oC_r(n)\to\oC_{r-1}(n-1)$. This defines $\partial''$.

We put $F_r(A)=\oplus\{\oC_r(n):n\geq 1\}$ and define $\partial:F_r(A)
\to F_{r-1}(A)$ by $\partial=\partial'+\partial''$. The exact relation
 between the homology of $F_{\bullet}(A)$ and groups $L_n(A)$ is as
yet unclear. However, we do have:
\begin{lemma}$H_3(F_{\bullet}(A))\otimes\Q\cong L_2(A)\otimes\Q\cong
H_3(\oC_{\bullet}(2))\otimes\Q$.
\end{lemma}
We sketch a proof. In view of the H-space structure, 
$L_i(A)\otimes\Q$ is the primitive homology of $\ET(A^n)$
with $\Q$ coefficients for $n$ large. The vanishing of
$H_1(n;\Q)$ for $n>2$ implies that the primitive homology 
is all of $H_i(n;\Q)$ for $i=2,3$ and $n$ large. By theorem~\ref{stab},
 we get $L_i(A)\otimes\Q\cong H_0(PGL_{i+1}(A),H_i(i+1;\Q))$
 for $i=2,3$. For the computation
 of $H_0(PGL(V),H_2(\ET(V);\Q))$ where $V=A^3$, we recall the exact
sequence obtained from $SS(V;\Q)$:
\\$0\to H_2(\ET(V);\Q)\to D(V)\otimes\Q\to D_2(V)=\oplus\{D(V/L):L\in\sL(V)\}
\to 0$.

This identifies $L_2(A)\otimes\Q$ with the 
cokernel of $H_1(PGL(V),D(V))\otimes\Q\to H_1(PGL(V),D_2(V))\otimes\Q$.
In view of IIa and the above remarks, this is readily identified
with $H_3(F_{\bullet}(A))\otimes\Q$.That gives the first isomorphism
 of the lemma. For the second isomorphism, what one needs is:
\\ \emph{Claim}: The arrow $H_4(\oC_{\bullet}(3))\to H_3(\oC_{\bullet}(2))$
 induced by $\partial''$ is zero.

The proof of this claim, which we address now, was already known to Spencer Bloch. Let $V=A^3$.
Given an ordered 5-tuple $(L_0,...,L_4)$ with the $L_i\in\sL(V)$ as
vertices of a 4-simplex in $K(V)$ (i.e.
in general position), they belong to a conic $C$ and the projection from
the points $L_i$ induces an isomorphism $p_i:C\to \PP(V/L_i)$.
We put $(M_0,...,M_4)=(p_0L_0,p_0L_1,...,p_0L_4)$. 
 Let $q_i=p_i\circ p_0^{-1}$. With the $\partial_i$ as in the definition
of $g_4$, we see that $\partial_i(L_0,...,L_4)\in C_3(V/L_i)$ and
 $q_i\partial_i(L_0,...,L_4)\in C_3(V/L_0)$ both give rise to the
same element of $\oC_3(2)$. It follows that 
 $\partial(M_0,M_1,...,M_4)\mapsto    \partial''(L_0,...,L_4)$
under the map $C_3(V/L_0)\to\oC_3(2)$.      
Thus $\partial''(L_0,...,L_4)\mapsto
 0\in H_3(\oC_{\bullet}(2)$. This proves the claim and the lemma.

Thus we have shown that 
\\$L_2(A)\tensor\Q\cong \coker(\oC_4(A^2)\to\oC_3(A^2))$.

The Bloch group tensored with $\Q$ is the homology of
\[\oC_4(A^2)\to\oC_3(A^2)\to \Lambda^2(A^{\times})\tensor\Q. \]
Thus this discussion amounts to a proof of Suslin's theorem on the Bloch
group. 

It remains to obtain a closed form for $L_3(A)\tensor\Q$
by this method. 

\textbf{Acknowledgements}: This paper owes, most of all, to Sasha Beilinson who introduced us
to the notion of a homotopy point. The formulation in terms of contractible spaces
 of maps occurs in many places in this manuscript. Rob de Jeu lent a
patient ear while we groped towards the definition of $\FL(V)$.
Peter May showed us a proof of Proposition~\ref{stable}. 
We thank them, Kottwitz and Gopal Prasad for interesting discussions, and
also  Spencer Bloch for his support throughout the project. Finally we thank the referee 
who brought our attention to several careless slips 
 in the first draft.


\Addresses
\end{document}